\documentclass[a4paper,11pt]{amsart}

\usepackage{amsfonts,amssymb,amsmath,amsthm,abstract,color}
\usepackage[mathscr]{euscript}
\usepackage[ps,all,arc,rotate]{xy}
\usepackage{enumerate}
\usepackage[lmargin=1in,rmargin=1in,tmargin=1in,bmargin=1in]{geometry}
\usepackage{fancyhdr}
\usepackage{pb-diagram}
\usepackage{graphicx}
\usepackage{multirow,enumitem,framed}
\usepackage{hyperref}
\usepackage{comment}
\usepackage{yfonts}
\usepackage{caption}
\usepackage{cleveref}
\usepackage{tikz-cd}

\usepackage{xcolor}

\newtheorem{prop}{Proposition}[section]
\newtheorem{lemma}[prop]{Lemma}
\newtheorem{thm}[prop]{Theorem}
\newtheorem{cor}[prop]{Corollary}

\theoremstyle{definition}
\newtheorem{defn}[prop]{Definition}

\newtheorem{rmk}[prop]{Remark}
\newtheorem{ex}[prop]{Example}
\newtheorem{ass}[prop]{Assumption}

\newtheorem{const}[prop]{Construction}

\DeclareMathOperator{\Proj}{Proj}
\DeclareMathOperator{\quot}{Quot}

\DeclareMathOperator{\diag}{diag}
\DeclareMathOperator{\GL}{GL}
\DeclareMathOperator{\SL}{SL}
\DeclareMathOperator{\PGL}{PGL}
\DeclareMathOperator{\SU}{SU}

\DeclareMathOperator{\U}{U}
\DeclareMathOperator{\proj}{Proj}
\DeclareMathOperator{\Hom}{Hom}
\DeclareMathOperator{\dom}{dom}

\DeclareMathOperator{\sym}{Sym}
\DeclareMathOperator{\conv}{Conv}
\DeclareMathOperator{\stab}{Stab}

\DeclareMathOperator{\id}{id}

\DeclareMathOperator{\Aut}{Aut}
\DeclareMathOperator{\End}{End}
\DeclareMathOperator{\Ext}{Ext}
\DeclareMathOperator{\Lie}{Lie}

\DeclareMathOperator{\gr}{gr}

\DeclareMathOperator{\Mat}{Mat}

\newcommand{\bexe}{\begin{exe}}
	\newcommand{\eexe}{\end{exe}}

\newcommand{\brf}{\begin{reference}}
	\newcommand{\erf}{\end{reference}}
\newcommand{\bnt}{\begin{nt} \normalfont}
	\newcommand{\ent}{\end{nt}}
\newcommand{\bntt}{\begin{ntt}}
	\newcommand{\entt}{\end{ntt}}
\newcommand{\bcl}{\begin{cl} \normalfont }
	\newcommand{\ecl}{\end{cl}}
\newcommand{\bqn}{\begin{qn} \normalfont  \begin{bf} }
		\newcommand{\eqn}{\end{bf} \end{qn}}
\newcommand{\bid}{\begin{idea}}
	\newcommand{\eid}{\end{idea}}

\newcommand{\bas}{\begin{assertion} \normalfont \begin{bf} }
		\newcommand{\eas}{\end{bf} \end{assertion}}
\newcommand{\bcr}{\begin{cor}}
	\newcommand{\ecr}{\end{cor}}
\newcommand{\bex}{\begin{example} \normalfont}
	\newcommand{\eex}{\end{example}}
\newcommand{\blm}{\begin{lemma}}
	\newcommand{\elm}{\end{lemma}}
\newcommand{\bthm}{\begin{thm}}
	\newcommand{\ethm}{\end{thm}}
\newcommand{\bcd}{\begin{tikzcd}}
	\newcommand{\ecd}{\end{tikzcd}}
\newcommand{\bdf}{\begin{defn} \normalfont}
	\newcommand{\edf}{\end{defn}}
\newcommand{\bpp}{\begin{prop}}
	\newcommand{\bam}{\begin{aim}}
		\newcommand{\eam}{\end{aim}}
	
	\newcommand{\rar}{\rightarrow}
	\newcommand{\brm}{\begin{rmk} \normalfont }
		\newcommand{\erm}{\end{rmk}}
	\newcommand{\epp}{\end{prop}}
\newcommand{\bpf}{\begin{proof}}

	\newcommand{\epf}{\end{proof}}

\newcommand{\kom}{$\mathcal{\kom}$}
\newcommand{\CC}{\mathbb{C}}

\newcommand{\rd}{\textcolor{red}}

\newcommand{\overbar}[1]{\mkern 1.5mu\overline{\mkern-1.5mu#1\mkern-1.5mu}\mkern 1.5mu}

\newcommand{\RR}{\mathbb{R}}
\newcommand{\bnu}{\begin{enumerate}}
	\newcommand{\enu}{\end{enumerate}}

\newcommand{\PP}{\mathbb{P}}
\newcommand{\ZZ}{\mathbb{Z}}
\newcommand{\NN}{\mathbb{N}}
\newcommand{\QQ}{\mathbb{Q}}
\newcommand{\GG}{\mathbb{G}}

\newcommand{\bdoc}{\begin{document}}
	\newcommand{\edoc}{\end{document}}
\newcommand{\bpm}{\begin{pmatrix}}
	\newcommand{\epm}{\end{pmatrix}}
\newcommand{\eps}{\varepsilon}
\newcommand{\bfct}{\begin{fact}}
	\newcommand{\efct}{\end{fact}}
\newcommand{\bslg}{\begin{slg}}
	\newcommand{\eslg}{\end{slg}}
\newcommand{\Uh}{\hat{U}}
\newcommand{\git}{\mathbin{
		\mathchoice{\mkern-5mu/\mkern-6mu/\mkern-4mu}
		{\mkern-5mu/\mkern-6mu/\mkern-4mu}
		{\mkern-5mu/\mkern-6mu/\mkern-4mu}
		{\mkern-5mu/\mkern-6mu/\mkern-4mu}}}

\newcommand{\rqq}{\rd{[??]}}
\newcommand{\re}{\rd{[...]}}

\newcommand{\diagentry}[1]{\mathmakebox[1.8em]{#1}}
\newcommand{\xddots}{%
	\raise 4pt \hbox {.}
	\mkern 6mu
	\raise 1pt \hbox {.}
	\mkern 6mu
	\raise -2pt \hbox {.}
}
\newcommand{\bgs}{\begin{gss}}
	\newcommand{\egs}{\end{gss}}
\def\acts{\curvearrowright}

\newcommand{\bit}{\begin{itemize}}
	\newcommand{\eit}{\end{itemize}}

\newcommand{\ra}{\rightarrow}     
\newcommand{\la}{\leftarrow} 
\newcommand{\surj}{\twoheadrightarrow}

\def\cA{\mathcal A}\def\cB{\mathcal B}\def\cC{\mathcal C}\def\cD{\mathcal D}
\def\cE{\mathcal E}\def\cF{\mathcal F}\def\cG{\mathcal G}\def\cH{\mathcal H}
\def\cI{\mathcal I}\def\cJ{\mathcal J}\def\cK{\mathcal K}\def\cL{\mathcal L}
\def\cM{\mathcal M}\def\cN{\mathcal N}\def\cO{\mathcal O}\def\cP{\mathcal P}
\def\cQ{\mathcal Q}\def\cR{\mathcal R}\def\cS{\mathcal S}\def\cT{\mathcal T}
\def\cU{\mathcal U}\def\cV{\mathcal V}\def\cW{\mathcal W}\def\cX{\mathcal X}
\def\cY{\mathcal Y}\def\cZ{\mathcal Z}

\def\AA{\mathbb A}\def\BB{\mathbb B}\def\CC{\mathbb C}\def\DD{\mathbb D}
\def\EE{\mathbb E}\def\FF{\mathbb F}\def\GG{\mathbb G}\def\HH{\mathbb H}
\def\II{\mathbb I}\def\JJ{\mathbb J}\def\KK{\mathbb K}\def\LL{\mathbb L}
\def\MM{\mathbb M}\def\NN{\mathbb N}\def\OO{\mathbb O}\def\PP{\mathbb P}
\def\QQ{\mathbb Q}\def\RR{\mathbb R}\def\SS{\mathbb S}\def\TT{\mathbb T}
\def\UU{\mathbb U}\def\VV{\mathbb V}\def\WW{\mathbb W}\def\XX{\mathbb X}
\def\YY{\mathbb Y}\def\ZZ{\mathbb Z}

\def\fa{\mathfrak a}\def\fb{\mathfrak b}\def\fc{\mathfrak c}\def\fd{\mathfrak d}
\def\fe{\mathfrak e}\def\ff{\mathfrak f}\def\fg{\mathfrak g}\def\fh{\mathfrak h}
\def\fj{\mathfrak j}\def\fk{\mathfrak k}\def\fl{\mathfrak l}\def\fm{\mathfrak m}
\def\fn{\mathfrak n}\def\fo{\mathfrak o}\def\fp{\mathfrak p}\def\fq{\mathfrak q}
\def\fr{\mathfrak r}\def\fs{\mathfrak s}\def\ft{\mathfrak t}\def\fu{\mathfrak u}
\def\fv{\mathfrak v}\def\fw{\mathfrak w}\def\fy{\mathfrak y}\def\fz{\mathfrak z}

\def\fA{\mathfrak A}\def\fB{\mathfrak B}\def\fC{\mathfrak C}\def\fD{\mathfrak D}
\def\fE{\mathfrak E}\def\fF{\mathfrak F}\def\fG{\mathfrak G}\def\fH{\mathfrak H}
\def\fI{\mathfrak I}\def\fJ{\mathfrak J}\def\fK{\mathfrak K}\def\fL{\mathfrak L}
\def\fM{\mathfrak M}\def\fN{\mathfrak N}\def\fO{\mathfrak O}\def\fP{\mathfrak P}
\def\fQ{\mathfrak Q}\def\fR{\mathfrak R}\def\fS{\mathfrak S}\def\fT{\mathfrak T}
\def\fU{\mathfrak U}\def\fV{\mathfrak V}\def\fW{\mathfrak W}\def\fX{\mathfrak X}
\def\fY{\mathfrak Y}\def\fZ{\mathfrak Z}

\def\SU{\mathrm{SU}} \def\U{\mathrm{U}}\def\GL{\mathrm{GL}} \def\SL{\mathrm{SL}}
\def\PGL{\mathrm{PGL}}

\newcommand{\env}{\!
	\mathbin{\text{\rotatebox[origin=c]{70}{\scalebox{1.2}{$\approx$}}}} \!}
\newcommand{\wenv}{\, \widehat{\env} \,}
\newcommand{\tenv}{\, \widetilde{\env} \,}
\newcommand{\inenv}{\dblslash \!_{\circ}}

\newcommand{\genv}{\dblslash \! _{\mathcal{C}}}
\newcommand{\ten}{\otimes}
\newcommand{\mc}{\mathcal}
\newcommand{\mf}{\mathfrak}
\newcommand{\mb}{\mathbb}
\newcommand{\mbf}{\mathbf}
\newcommand{\ol}{\overline}
\newcommand{\pt}{\mathrm{pt}}
\newcommand{\symdot}{\sym^{\bullet}}
\newcommand{\nss}{\mathrm{nss}}
\newcommand{\ssfg}{\mathrm{ss,fg}}
\newcommand{\ssUfg}{\mathrm{ss,} H_u-\mathrm{fg}}
\newcommand{\rms}{s}
\newcommand{\rmss}{\mathrm{ss}}
\newcommand{\ssC}{\mathrm{ss},\mc{C}}
\newcommand{\aff}{\mathrm{aff}}
\newcommand{\dblslash}{/\! \!/}
\newcommand{\act}{\curvearrowright}
\newcommand{\umax}{{U_{\max}}}
\newcommand{\weight}{\omega}
\newcommand{\weightlie}{\hat{\omega}}
\newcommand{\reg}{\mathrm{reg}}
\newcommand{\sss}{s_1,\ldots, s_n}

\newcommand{\syms}{\mathrm{sym}}
\newcommand{\symk}{\sym^{\mathbf{\weight}}\CC^n}
\newcommand{\wsymk}[1]{\wedge^{#1}(\symk)}
\newcommand{\grass}{\mathrm{Grass}}
\newcommand{\bg}{\mathbf{\gamma}}
\newcommand{\bX}{\overline{X}}
\newcommand{\bY}{\overline{Y}}
\newcommand{\bZ}{\overline{Z}}
\newcommand{\bR}{\overline{R}}
\newcommand{\bz}{\mathbf{z}}
\newcommand{\ta}{\tilde{\a}}
\newcommand{\ts}{\tilde{s}}
\newcommand{\tk}{\tilde{k}}
\newcommand{\bm}{\mathbf{m}}
\newcommand{\bbg}{\mathbb{G}}
\newcommand{\br}{\mathbf{r}}
\newcommand{\bk}{\mathbf{\omega}}
\newcommand{\bi}{\mathbf{i}}
\newcommand{\bj}{\mathbf{j}}
\newcommand{\bv}{\mathbf{v}}
\newcommand{\bu}{\mathbf{u}}
\newcommand{\bl}{\mathbf{\lambda}}
\newcommand{\zdis}{\mathfrak{p}}
\newcommand{\qdis}{\mathfrak{q}}

\newcommand{\ba}{\mathbf{\alpha}}
\newcommand{\bs}{\mathbf{s}}

\newcommand{\Fl}{\mathcal{F}}
\newcommand{\be}{\mathbf{e}}
\newcommand{\imm}{\mathrm{im}}
\newcommand{\fnk}{\mathfrak{Fl}_k(n)}
\newcommand{\hofi}{\Hom^\triangle(\CC^k,\symdot)}
\newcommand{\phit}{\bar{\phi}}
\newcommand{\bmnd}{\calm(n,d)}
\newcommand{\liek}{{\mathfrak k}}
\newcommand{\lieu}{{\mathfrak u}}
\newcommand{\kk}{\Bbbk}

\newcommand{\hU}{\widehat{U}}
\newcommand{\hH}{\widehat{H}}
\newcommand{\hX}{\widehat{X}}
\newcommand{\hZ}{\widehat{Z}}
\newcommand{\hp}{\widehat{p}}
\newcommand{\hq}{\widehat{q}}
\newcommand{\hs}{\hat{s}}
\newcommand{\Gm}{\mathbb{G}_m}
\newcommand{\ct}{ CITE }

\title[quotients by parabolic groups and moduli spaces of unstable objects]{Quotients by Parabolic Groups and Moduli Spaces of Unstable Objects}

\author{Victoria Hoskins}
\address{V.\ Hoskins \\IMAPP, Radboud University \\ PO Box 9010, 6500GL Nijmegen \\ The Netherlands}
\email{v.hoskins@math.ru.nl}

\author{Joshua Jackson}
	\address{J.\ Jackson  \\ Mathematics Department, Imperial College \\London \\ SW7 2AZ \\UK  \vspace{-5pt}}
	\address{Heilbronn Institute for Mathematical Research\\ Bristol \\ UK}
\email{j.jackson@imperial.ac.uk}

\begin{document}
	
\maketitle
	
	\begin{abstract}
	
Motivated by constructing moduli spaces of unstable objects, we use new ideas in non-reductive GIT to construct quotients by parabolic group actions. For moduli problems with semistable moduli spaces constructed by reductive GIT, we consider associated instability (or HKKN) stratifications, which are often closed related to Harder--Narasimhan stratifications, and construct quotients of the unstable strata under various stabiliser assumptions by further developing ideas of non-reductive GIT. Our approach is to construct parabolic quotients in stages, in order for the required stabiliser assumptions to be more readily verified.   
To illustrate these ideas, we construct moduli spaces for certain sheaves of fixed Harder--Narasimhan type on a projective scheme in cases where our stabiliser assumptions can be verified.
\end{abstract}

\setcounter{tocdepth}{1} 	
\renewcommand{\tocsection}[3]{	\indentlabel{{\bfseries\ignorespaces#1 #2\quad}}\bfseries#3} %makes sections bold in toc
\tableofcontents

	\section{Introduction}

Many algebro-geometric moduli spaces are constructed using Mumford's Geometric Invariant Theory (GIT) \cite{Mumford}: given a reductive group $G$ acting linearly on a projective scheme $Y$, one obtains a quotient of an open semistable locus $Y^{ss} \subset Y$. Semistability is combinatorially determined via the \emph{Hilbert-Mumford criterion} by calculating the weights of 1-parameters subgroups (1-PS) of $G$. In good cases, this enables a moduli-theoretic interpretation of GIT semistability from which one obtains a moduli space of semistable objects; the paradigmatic example being slope semistability for vector bundles on a smooth projective curve.

In this paper, we view the above as merely the beginning of the story and aim to additionally construct moduli spaces of unstable objects. Instability in reductive GIT is a structured phenomenon: for each unstable point of $Y$, there is a \emph{maximally destabilising} 1-PS of $G$, unique up to conjugation under a parabolic subgroup, which is most responsible\footnote{This notion depends on a choice of norm on 1-PSs of $G$.} for violating the Hilbert--Mumford criterion \cite{Kempf}; this is akin to the Harder--Narasimhan filtration \cite{HN} of a vector bundle being a maximally destabilising filtration. These maximally destabilising 1-PSs gives rise to an instability stratification (or HKKN stratification following \cite{Hesselink,Kempf,Kirwan_thesis,Ness}) of $Y$ into finitely many $G$-invariant locally closed subschemes such that the lowest stratum $S_0$ is the semistable set and the higher (unstable) strata $S_\beta$ are indexed by conjugacy classes $\beta$ of rational 1-PSs of $G$ which are maximally destabilising for the points in $S_\beta$. Recently, Halpern-Leistner greatly generalised this idea to $\Theta$-stratifications of stacks \cite{HL_Theta}. The HKKN stratification suggests that the instability types $\beta$ should be viewed as discrete invariants and one should then look to construct moduli spaces for unstable objects of instability type $\beta$. For the GIT construction of moduli of vector bundles on a curve,  GIT instability types correspond to Harder--Narasimhan (HN) types \cite{GSZ,HK}.

Our primary motivating examples come from moduli of objects in a linear abelian category, where moduli spaces of semistable objects are constructed by reductive GIT and the HKKN stratification is closely related to a Harder-Narasimhan (or Shatz \cite{Shatz}) stratification; for example, this is the case for moduli of bundles or sheaves \cite{HK,Hoskins}, Higgs sheaves \cite{Eloise} and quiver representations \cite{Hoskins_quivers}. In these cases, constructing quotients of the unstable strata $S_\beta$ would give rise to moduli spaces of objects of fixed HN type; this has been described for certain length $2$ HN filtrations in \cite{BrambilaPaz2013,Eloise,Josh_length2}.

\subsection{Moduli of Unstable Objects}\label{sec intro mod unstable}

The initial technical obstruction to constructing quotients of unstable strata is clear enough: the strata consist by definition of GIT-unstable points, so their quotients are empty. Since GIT semistability depends on a choice of linearisation, a natural approach is to vary the linearisation $\cL$ by twisting by a character to change the semistable locus, as is described by variation of GIT \cite{DH,Thaddeus}. However, there may not exist a character $\chi_\beta$ of $G$ which we can use to make all of $S_\beta$ semistable (for example, if $G=\SL_n$, there are no non-trivial characters). Fortunately, if we pick a representative 1-PS $\lambda_\beta$ of the conjugacy class $\beta$, this determines a parabolic group $P_\beta:=P(\lambda_\beta)< G$ and locally closed subscheme $Y_\beta^{ss} \subset S_\beta$ such that $S_\beta$ is the $G$-sweep of $Y_\beta^{ss}$, and such that flowing under $\lambda_\beta$ induces a retraction $p_\beta : Y_\beta^{ss} \ra Z_\beta^{ss}$ \cite{Kirwan_thesis}. Moreover, a categorical quotient of $G$ acting on $S_\beta$ is equivalent to a categorical quotient of $P_\beta$ acting on $Y_\beta^{ss}$. The parabolic subgroup $P_\beta$ has larger character group than $G$ and in particular, there is a character $\chi_\beta$ such that the twisted linearisation $\cL_\beta$ (which we refer to as a borderline linearisation below) produces a categorical $P_\beta$-quotient of $Y_\beta^{ss}$ \cite{HK}. However, this quotient factors by the retraction $p_\beta$ followed by the categorical quotient of $Z_\beta^{ss}$ under the action of the Levi $L_\beta < P_\beta$. In particular, this categorical quotient is far from being an orbit space. Often the moduli theoretic interpretation of the retraction $p_\beta$ is given by taking the associated graded of a Harder--Narasimhan filtration and so the obtained categorical quotient is just a product of moduli spaces for the successive semistable quotients in the HN filtration. 

To avoid these identifications, one would like to remove the locus $Z_\beta^{ss} \subset Y_\beta^{ss}$, which in many moduli theoretic examples parametrises objects with completely degenerate HN filtrations. To do this, we would like to further perturb the twisted linearisation $\cL_\beta$, so that $Z_\beta^{ss}$ becomes unstable. The $P_\beta$-linearisation $\cL_\beta$ on $Y_\beta^{ss}$ also extends to an ample $P_\beta$-linearisation on the closure of  $Y_\beta^{ss}$ in $Y$, but it does not extend to an ample $G$-linearisation on a natural compactification of $S_\beta$. For this reason, one must work with the non-reductive group $P_\beta$ and so classical reductive GIT cannot be applied.

Recently there has been great progress on developing non-reductive GIT (NRGIT) for groups whose unipotent radical is \lq graded' by a one-parameter subgroup \cite{BDHK2,BDHK_handbook,BDHK,BDK_grad_lin,BHK}. The best results are obtained under certain strong stabiliser assumptions by appropriately twisting the linearisation to make it \lq adapted'. Even when these stabiliser assumptions fail, one can perform a blow-up sequence to obtain these assumptions on the blow-up. Although there are important applications in which the blow-up process is tractable (for example, see\ \cite{BK_hyp}), in general explicitly determining the open set one obtains a quotient of is often quite involved (see $\S$\ref{sec non red GIT}). This theory has applications to hyperbolicity questions \cite{BK_hyp} and the construction of moduli of hypersurfaces in toric varieties \cite{dom}, as well as providing tools for cohomological computation via a symplectic description of the quotient \cite{BK_cohomolgy}; but in this paper we will focus on its application to the construction of moduli spaces of unstable objects, motivated by the desire to construct moduli spaces of objects of fixed Harder--Narasimhan type. More precisely, given an unstable stratum $S_\beta$ associated to a reductive GIT set-up, our aim is to construct a non-reductive quotient of the action of the parabolic group $P_\beta$ on $Y_\beta^{ss}$ (or rather its closure). This is an ideal situation to apply non-reductive GIT, as the unipotent radical of $P_\beta$ is graded by the 1-PS $\lambda_\beta$. Consequently, we can apply NRGIT to 
construct explicit quotients of unstable HKKN strata subject to certain stabiliser assumptions denoted \ref{starU}, \ref{starURss} and \ref{starR0} (see Definition \ref{def star cond}). The first main result of this paper is the following theorem concerning quotients of unstable HKKN strata.

\begin{thm}\label{thm nred quot unstable strata}\label{mainthm1}
	For a reductive group $G$ acting linearly on a projective scheme $Y$, let $\beta$ be a HKKN index corresponding to an unstable stratum $S_\beta \subset Y$. Consider the non-reductive group $P_\beta = U_\beta \rtimes L_\beta$ with unipotent radical $U_\beta$ internally graded by $\lambda_\beta$ acting on $X_\beta:= \overline{Y_\beta}$ with respect to the ample adapted linearisation $\cL_\beta^{\text{per}}$. Then the following statements hold.
	\begin{enumerate}[label=\roman*)]
		\item Assume \ref{starU} holds for the action of $\hU = \hU_\beta$; then the open stable set $X_\beta^{\hU_\beta-s}:=Y_\beta - U_\beta  Z_\beta$ admits a geometric and projective $\hU_\beta$-quotient $\overline{Y_\beta}^{\hU_\beta-s}/\hU_\beta$, given by the projective spectrum of the $\hU$-invariants with respect to a well-adapted twist of $\cL_{\beta}^{\text{per}}$.
		\item Assume \ref{starURss} and \ref{starR0} hold for the action of $P_\beta$; then the open stable set \[X_\beta^{P_\beta-s} := Y_\beta^{ss} - UZ_\beta^{ss}=Y_\beta^{s} - UZ_\beta^{s}\] admits a geometric and projective $P_\beta$-quotient, given by the projective spectrum of the $P_\beta$-invariants with respect to a well-adapted twist of $\cL_{\beta}^{\text{per}}$.
		\item Assuming \ref{starU} holds, $U_\beta  Z_\beta \subset Y_\beta$ is a closed subscheme, the map $p_\beta : U_\beta Z_\beta \ra Z_\beta$ is a geometric $\hU_\beta$-quotient and the reductive quotient $Z_\beta /\!/_{\cL_\beta} L_\beta$ is a good $P_\beta$-quotient of $U_\beta Z_\beta^{ss}$ (and also a geometric quotient if \ref{starR0} holds).
	\end{enumerate}
\end{thm}

To apply this theorem, one needs to check whether these stabiliser assumptions are satisfied. The assumptions \ref{starU} and \ref{starURss} concern unipotent stabilisers, whereas \ref{starR0} concerns stabilisers for the reductive Levi factor and is used to obtain an explicit description of the open \lq semistable' set of which we obtain a quotient. If these stabiliser assumptions fail, then NRGIT uses a blow-up sequence to produce a (only quasi-projective) quotient of an open set, which in practice can be hard to determine.

Unfortunately, for many unstable strata, the reductive stabiliser condition \ref{starR0} fails. More precisely, \ref{starR0} asks for semistability to coincide with stability for the action of the reductive Levi factor $L_\beta$ on $Z_\beta$ (which is the minimal weight space in the $\lambda_\beta$-fixed locus). Often in the case when our reductive GIT set-up comes from moduli of objects in a linear abelian category, the subscheme $Z_\beta$ parametrises the associated graded of Harder--Narasimhan filtrations and the stabiliser group of such an object contains a torus of dimension equal to the length of the HN filtration; this results in a torus $T_\beta$ in the Levi $L_\beta$ acting trivially on $Z_\beta$ (but not on $Y_\beta$) so that \ref{starR0} fails except in low lengths. More precisely, in the case of sheaves of HN type of length $l$, one obtains an $(l-1)$-dimensional\footnote{Note that the dimension drops by one, as we are working inside $\SL$ and not $\GL$.} torus $T_\beta <  L_\beta$ acting trivially on $Z_\beta$ (see $\S$\ref{sec sheaves nonred git setup}). Hence, \ref{starR0} fails for HN filtrations of length $l > 2$; this happens for moduli of (Higgs) sheaves and quiver representations. For HN filtrations of length $l=2$, the problem of constructing moduli spaces of fixed HN type for sheaves and Higgs bundles are considered in \cite{Josh_length2} and \cite{Eloise}.

\subsection{Quotienting-in-Stages}\label{sec intro QiS}

The main result of this paper is motivated by the issue outlined above, namely that in many applications of interest (both within \lq moduli spaces of unstable objects', and beyond) the stabiliser conditions of Theorem \ref{mainthm1} are not met and so we cannot obtain an explicit satisfactory quotient using the current methods of NRGIT. 

To overcome this, we prove a new theorem in NRGIT that is better suited to this family of applications than Theorem \ref{mainthm1}.  More precisely, we consider a projective scheme $X$ with a linear action of a parabolic subgroup $P =  P(\lambda) = U \rtimes L < \SL_N$, defined by a one-parameter subgroup $\lambda: \Gm \rar \SL_N$. The basic idea is to break the quotienting process up into stages, one for each step in the row filtration of $P$ (one should consider these as the steps in the HN filtration).  At each stage, we take a NRGIT quotient of a maximal unipotent group inside $U$ graded by a 1-PS coming from the central torus $T = Z(L)$. The key challenges are determining an explicit description of the open locus we get a quotient of at all these stages and determining the relevant stabiliser conditions which are required for this procedure. To determine this open set explicitly, we impose the so-called \emph{Quotienting-in-Stages Assumption} \hyperlink{QiS}{(QiS)}, which are modelled on the behaviour of objects of fixed HN type in a linear abelian category and are used to determine the relevant minimal weight spaces at each stage in the quotient procedure. In place of the strong stabiliser assumptions of Theorem \ref{mainthm1} we introduce conditions on stabiliser groups that are more likely to hold in this setting, which we call the \emph{Weak Upstairs Unipotent stabiliser Assumption} \hyperlink{WUU}{(WUU)} and the \emph{Upstairs Stabiliser Assumption} \hyperlink{U}{(U)}. Under these assumptions, we obtain an explicit quotient as follows.

\begin{thm}\label{mainthm2}
	Let $P < \SL_N$ be a parabolic group acting linearly on an irreducible projective scheme $X$ satisfying \hyperlink{QiS}{(QiS)}. If \hyperlink{WUU}{(WUU)} holds, then there is an open set $X^{P-qs}$ with an explicit Hilbert-Mumford style description that admits a geometric $P$-quotient given by a quasi-projective scheme with a natural projective completion. If moreover, \hyperlink{U}{(U)} holds, then this quotient of $X^{P-qs}$ is projective.
\end{thm}	

The Upstairs Stabiliser Assumption \hyperlink{U}{(U)} is chosen to induce the relevant stabiliser conditions at each stage of the quotienting process, so that the stronger results of NRGIT can inductively be applied at each stage. Crucially, the reductive stabiliser assumptions in \hyperlink{U}{(U)} do not concern the full Levi $L$, but only its semisimple part and so the fact that the central torus $T=Z(L)$ acts trivially is no longer an issue.
By performing a blow-up sequence, we can weaken this to \hyperlink{WUU}{(WUU)}, which ensures we understand enough about the blow-up procedure to explicitly describe the obtained quotient.

\subsection{Applications}

As mentioned above, the motivating family of applications for our results is the construction of moduli spaces of objects of fixed HN type in a linear abelian category for which a reductive GIT construction of the moduli space of semistable objects is known. These applications, including sheaves, Higgs sheaves, and quiver representations, will be dealt with in future papers \cite{BHJK,HHJ}.

For ease of exposition, we illustrate our ideas by briefly analysing the problem of constructing moduli spaces for sheaves of fixed HN type, concluding with Theorem \ref{thm sheaves} which allows the construction of moduli spaces for certain unstable sheaves in special situations. The open stable set $X^{P-qs}$ obtained by Theorem \ref{mainthm2} in the case of sheaves of HN type $\tau$ admits a moduli-theoretic interpretation which we call $\tau$-stability, see Definition \ref{def tau stable}: for a sheaf of HN type $\tau$ to lie in $X^{P-qs}$, each inclusion in its HN filtration must be non-split and the successive HN quotients must be stable rather than just semistable. In addition, one requires certain unipotent stabiliser assumptions, which admit interpretations in terms of certain filtered endomorphism groups of the sheaf. The semistability notion given by NRGIT for sheaves of fixed HN type naturally fits with previous work describing moduli of rank 2 bundles of fixed HN type \cite{BrambilaPaz2013} using a projectivised Extension bundle associated to the universal families for the two subquotients: this forces the HN filtration to be non-split and, for this Extension sheaf to be a vector bundle, one requires constant dimensional filtered endomorphism groups. More generally, for a length $l$ HN type, provided one has universal families for the semistable subquotients, one can iteratively take projectivised Extension bundles to obtain an analogous construction. However, without the existence of such universal families, one can instead apply the quotienting-in-stages procedure we develop.

However, the assumptions of Theorem \ref{mainthm2} are very restrictive, due to the necessity of condition \hyperlink{WUU}{(WUU)}, which at present we can only verify for sheaves in slightly contrived examples (see $\S$ \ref{sec applications}). A construction better adapted to the case of sheaves in particular will be considered in \cite{BHJK}. On the other hand, condition \hyperlink{WUU}{(WUU)} is satisfied in far greater generality for Higgs sheaves, and so the results of this paper are well suited to the construction of moduli spaces of unstable Higgs sheaves; we will pursue this application in \cite{HHJ}.

\subsection{Overview of the paper}
In $\S$\ref{sec GIT}, we recall the basic properties of reductive GIT and instability stratifications. Moreover, in \S\ref{sec non red GIT} we collect the results of NRGIT we shall need and prove some slightly stronger results. We then apply these results in $\S$\ref{sec quot HKKN strata}, leading to the proof of our first main result, Theorem \ref{mainthm1} on quotienting unstable strata. 

The remainder of the paper is devoted to proving Theorem \ref{mainthm2}. Let us sketch the structure; for the detailed road map, see $\S$\ref{sec roadmap}. In Section \ref{sec quotient stages par}, we explain the set-up and notation, including defining the quotienting-in-stages stable locus $X^{P-qs}$ in Definition \ref{defn quotienting-in-stages stable locus}. We outline the quotienting-in-stages process in $\S$\ref{sec outline quotient in stages}, which we can perform under so-called \emph{Downstairs Stabiliser Assumption} \hyperlink{D}{(D)} conditions. From $\S$\ref{subsec qis and u assumptions} onwards, we suppose the \emph{Quotienting-in-Stages Assumptions} \hyperlink{QiS}{(QiS)} hold. We then introduce certain upstairs stabiliser assumptions \hyperlink{U}{(U)} and \hyperlink{WUU}{(WUU)}, before moving onto the proof in Section \ref{sec proof}. Assuming \hyperlink{D}{(D)}, we first describe the open locus we obtain a quotient of in $\S$\ref{sec describing domain q with downstairs stab}, and we show that \hyperlink{D}{(D)} can be deduced from \hyperlink{U}{(U)} in $\S$\ref{subsec from downstairs to upstairs}. To deduce the final version of Theorem \ref{mainthm2} we perform certain blow-up procedures. We conclude with a brief discussion of applications in $\S$\ref{sec applications}.

\subsection*{Notation and Conventions}
By a scheme, we mean a scheme of finite type over an algebraically closed field of characteristic 0. By a sheaf over a scheme $X$, we mean a coherent sheaf of $\cO_X$-modules. For a linear algebraic group $G$ and a subgroup $H \subset G$ acting on a quasi-projective scheme $Y$, the Borel construction $G \times^H Y$ denotes the geometric quotient of $H$ acting on $G \times Y$ diagonally.

\subsection*{Acknowledgements}
We are very grateful to F. Kirwan, G. B\'{e}rczi, D. Bunnett and E. Hamilton for many helpful discussions concerning non-reductive GIT and the ideas in this paper. The second named author is funded by the Heilbronn Institute for Mathematical Research.

\section{Geometric invariant theory}\label{sec GIT}
	
	For an action of a linear algebraic group on a scheme, the aim of geometric invariant theory (GIT) is to provide a quotient of this action (potentially on an invariant open subscheme) in the category of schemes. The theory splits into two branches, according to whether the group acting is reductive or non-reductive. The former was established in the 1960s by Mumford \cite{Mumford}, building on the results of 19th Century invariant theory. The latter has been studied by several authors; in particular, the recently developed theory for non-reductive groups with so-called graded unipotent radicals has many of the favourable properties of reductive GIT \cite{BDHK2}. In $\S$\ref{sec red GIT}, we survey the key results of reductive GIT, including the Hilbert--Mumford criterion, and the partial desingularisation blow-up procedure \cite{K2}; in $\S$\ref{sec HKKN strat} we review the instability (HKKN) stratifications arising from work of Hesselink \cite{Hesselink} Kempf \cite{Kempf}, Kirwan \cite{Kirwan_thesis} and Ness \cite{Ness}. In $\S$\ref{sec non red GIT}, we summarise the key results from \cite{BDHK2} which enable constructions of GIT quotients for actions of linear algebraic groups with graded unipotent radicals. 
	
	\subsection{Reductive GIT quotients}\label{sec red GIT}
	
Let $G$ be a reductive group acting on a projective scheme $Y$ with an ample linearisation of the action $\cL$. One of the most well-known properties of reductive groups actions is that the ring of $G$-invariants is finitely generated; this is contrast to the counterexamples to Hilbert's 14th Problem \cite{Nagata} in the non-reductive case. In addition taking $G$-invariants is exact, which makes the theory behave functorially with respect to $G$-equivariant closed immersions. These properties were utilised by Hilbert and Mumford in order to construct reductive GIT quotients. A standard reference is \cite{Newstead}.

The reductive GIT quotient $\pi : Y^{ss} \ra Y/\!/_{\!\cL} G$ is a categorical and good quotient of the open subset $Y^{ss}:=Y^{G-ss}(\cL)$ of GIT semistable points in $Y$ (with respect to the $G$-linearisation $\cL$) and was constructed by Mumford \cite{Mumford}. The quotient $Y/\!/_{\!\cL} G$ is a projective scheme equal to the projective spectrum of the graded ring of $G$-invariant sections of non-negative powers of $\cL$. In general, the GIT quotient is not an orbit space: its closed points correspond to S-equivalence classes of semistable orbits, where two semistable orbits are S-equivalent if their orbit closures meet in the semistable set. However, the GIT quotient restricts to a geometric quotient, which in particular is an orbit space, on the open subset $Y^{s}(\cL)\subset Y^{ss}$ of stable points. 

Stable points\footnote{Note that stable points in modern terminology were referred to as properly stable points by Mumford \cite{Mumford}.} can be described as semistable points whose stabiliser groups are finite and whose orbits are closed in the semistable set. Often when the linearisation and group is fixed, we will drop the notation $\cL$ and $G$ indicating the dependence of the semistable set on this linearisation and $G$. Variation of GIT describes the birational transformations between GIT quotients with respect to different $G$-linearisations \cite{DH,Thaddeus}.
	
	The (semi)stable loci are defined in \cite{Mumford} using $G$-invariant sections of positive powers of $\cL$. Critical to the use of reductive GIT in practice, though, is the ability to determine these loci without computing ring of invariants, afforded by the Hilbert--Mumford criterion. This criterion states that a closed point $y \in Y$ is semistable if and only if it is semistable for each 1-parameter subgroup (1-PS) $\lambda$ of $G$. To state this as a numerical criterion, we define the Hilbert--Mumford weight $\mu(y,\lambda)$ to be minus the weight of the $\GG_m$-action induced by $\lambda$ on the fibre of $\cL$ over the fixed point $\lim_{t \ra 0} \lambda(t) \cdot y$. 
	
	\begin{prop}[Hilbert--Mumford criterion]
		A closed point $y \in Y$ is semistable if and only if $\mu(y,\lambda) \geq 0$ for all 1-PS $\lambda: \GG_m \ra G$.
	\end{prop}
	
	Alternatively, by using the linearisation $\cL$ to give an equivariant embedding $Y \hookrightarrow \PP^n$, the Hilbert--Mumford criterion can be stated combinatorially using the weights of a maximal torus $T < G$ on $\PP^n$ (or rather the affine cone over this projective space). A point is semistable for the $T$-action if its convex hull of $T$-weights contains the origin \cite[$\S$9.4]{Dolgachev}. Moreover, a point is semistable for the $G$-action if all points in its $G$-orbit are $T$-semistable; thus we have \[Y^{G-ss} = \cap_{g \in G} \: gY^{T-ss}.\]
	
	\subsubsection{Partial desingularisation via equivariant blow-ups}\label{sec red part desing}
	
	In this section, we summarise a partial desingularisation of the reductive GIT quotient $Y/\!/G$ that is constructed by performing a canonical sequence of blow-ups along closed invariant subschemes determined by studying the dimensions of reductive stabiliser groups of semistable points \cite{K2}. In \cite{K2} this procedure was described when the stable locus was non-empty, in which case one has that semistability and stability coincide on the final blow-up. If one starts with a smooth projective scheme $Y$, the resulting blow-up $\widetilde{Y}$ is also smooth and projective and its GIT quotient $\widetilde{Y}/\!/G$ has only orbifold singularities and thus provides a partial desingularisation of $Y/\!/G$. One can also perform this blow-up procedure when then stable set is empty but the semistable set is non-empty; in this case, the procedure can terminate in several different ways, as we describe below. 
	
To describe the blow-up sequence in \cite{K2} for irreducible $Y$ and also later describe the non-reductive GIT blow-up sequences, we introduce the following notation for the locus in a subscheme where the stabilisers of a certain group have fixed dimension.

\begin{defn}\label{def centres of blowups}
Given a scheme $Y$ with an action by a linear algebraic group $G$, for any subscheme $Y' \subset Y$ and subgroup $G' < G$ and any $d \in \NN$, we define the locally closed subscheme 
\[ C(Y',G',d) := \{y \in Y' : \dim \stab_{G'}(y) =d \} \subset Y' \]
to be the locus in $Y'$ where the $G'$-stabilisers are $d$-dimensional. We also write
\[ d_{\max}(Y',G') :=\max\{ \dim \stab_{G'}(y) : y \in Y' \}  \quad \text{and} \quad d_{\min}(Y',G') :=\min\{ \dim \stab_{G'}(y) : y \in Y' \}\]
and $C(Y',G'):= C(Y',G',d_{\max}(Y',G'))$ for the locus with maximal dimensional stabilisers.
\end{defn}

Often for simplicity we will assume that $Y$ is irreducible (in the case when $Y'$ is not irreducible, for example if it is not connected, it may be better to think of $d_{\max}(Y',G')$ as a tuple of dimensions indexed by the irreducible components). Since the dimension of the stabiliser group is upper semi-continuous, $C(Y',G') $ is a closed subscheme of $Y'$. If $Y'$ is $G'$-invariant, these subschemes are $G'$-invariant. If $Y'$ is $G$-invariant and $G' \triangleleft G$ is normal then these subschemes are $G$-invariant, so that the blow up of $Y$ along $Y'$ is $G$-equivariant, and the blown up space has a natural $G$-action. 

Returning to the case where $G$ is a reductive group acting linearly on a projective scheme $Y$,	if $\emptyset \neq  Y^{s} \subsetneq Y^{ss}$, then $C(Y^{ss},G)$ is a proper closed $G$-invariant subscheme of $Y^{ss}$ which is disjoint from $Y^{s}$. Assuming that $Y$ is irreducible for simplicity, the first step\footnote{In fact, this first step in \cite{K2} is performed as a sequence of blow-ups: one chooses representatives of all conjugacy classes of connected reductive subgroups of $G$ of dimension $d_{\max}(Y^{ss},G)$ and for each such subgroup $R$ one blows-up the closure of the $G$-sweep of the locus of semistable points which are fixed by $R$.} in the partial desingularisation procedure in \cite{K2} is to blow-up $X$ along the closure of $C(Y^{ss},G)$ to obtain $\varphi_1: \widetilde{Y}_{1} \ra Y$. The $G$-action on $Y$ naturally lifts to $\widetilde{Y}_{1}$ and by pulling back the ample linearisation $\cL \ra Y$ and twisting by a small multiple of the exceptional divisor one obtains an ample $G$-linearisation $\cL_1 \ra \widetilde{Y}_1$ such that
	\[ \varphi_1^{-1}(Y^s) \subset \widetilde{Y}_1^s \subset \widetilde{Y}_1^{ss} \subset \varphi_1^{-1} (Y^{ss})\]
and also
\[ d_{\max}(\widetilde{Y}^{ss}_1,G) < d_{\max}(Y^{ss},G).\]
By inductively blowing-up $\widetilde{Y}_i$ along the closure of $C(\widetilde{Y}^{ss}_i,G)$, we eventually obtain a projective scheme $\widetilde{Y}$ for which semistability and stability coincide. 

If the stable set is empty, the partial desingularisation procedure can terminate in different ways; however, we will not be interested in these cases in this paper, so we refer to \cite{BHK} for an overview of the possible ways in which this can terminate.
	
	\subsubsection{Example: moduli of semistable sheaves via reductive GIT}\label{sec red GIT sheaves}
	
Let us recall the construction of moduli spaces of (Gieseker) semistable sheaves with fixed Hilbert polynomial $P$ over a projective scheme $B$ with fixed ample line bundle $\cO_B(1)$ following Simpson \cite{Simpson}. As the set of semistable sheaves with Hilbert polynomial $P$ is bounded, one can take $n$ sufficiently large so they are all $n$-regular; thus for $\cE$ semistable we have a surjective evaluation map
	\[  H^0(\cE(n)) \otimes \cO_B(-n) \ra \cE\]
	which with a choice of isomorphism $\CC^{P(n)} \cong H^0(\cE(n))$ determines a point in the Quot scheme 
	\[Q_n:=Q_n(P):=\quot(\cO_B(-n)^{\oplus P(n)},P)\] 
	parametrising quotient sheaves of $\cO_B(-n)^{\oplus P(n)}$ with Hilbert polynomial $P$. Let $Q_n^{ss}:=Q^{ss}_n(P)$ be the open subscheme of $Q_n$ that parametrises quotient sheaves $q: \cO_B(-n)^{\oplus P(n)} \surj \cE$ such that $\cE$ is semistable and $H^0(q(n))$ is an isomorphism; then the isomorphism classes of semistable sheaves over $B$ with Hilbert polynomial $P$ are in bijective correspondence with the orbits of the natural $\GL_{P(n)}$-action on $Q^{ss}_n$. The moduli space of semistable sheaves is constructed as a reductive GIT quotient of $\SL_{P(n)}$ acting on the closure $R_n:=R_n(P)$ of $Q^{ss}_n$ in $Q_n$ with respect to the linearisation $\cL_{n,m}$ for $m >\!> n$ corresponding to Grothendieck's embedding of the Quot scheme $Q_n$ into a Grassmanninan of $P(m)$-dimensional quotients of $H^0(\cO_B(m-n))^{\oplus P(n)}$.	
	
	\subsection{Instability stratifications}\label{sec HKKN strat}

	Associated to a linearised action of the reductive group $G$ on $Y$, there is a GIT instability, or HKKN, stratification of $Y$ due to work of Kempf \cite{Kempf} and Hesselink \cite{Hesselink}. It also agrees with a symplectic Morse-theoretic stratification associated to the norm square of the moment map for the action of the maximal compact subgroup of $G$ \cite{Kirwan_thesis, Ness}.
	
	To define this stratification, we must also choose an invariant norm $||-||$ on the Lie algebra of $G$. More precisely, we mean a norm on $\Lie G$ that is invariant under the adjoint action and such that for a (or equivalently every) maximal torus $T \subset G$, the norm is integer valued on the co-character lattice of $T$. In particular, this norm gives an identification of the co-character lattice and character lattice of $T$ and of the Lie algebra and co-Lie algebra of $T$. 
	
	The HKKN stratification of $Y$ is then a finite stratification by $G$-invariant locally closed subschemes $S_\beta$ such that the lowest stratum $S_0$ is $Y^{ss}$ and the unstable strata are indexed by conjugacy classes of rational 1-PSs of $G$. The idea is to associate to each unstable point $y \in Y - Y^{ss}$ the conjugacy class $\beta=[\lambda_\beta]$ of rational 1-PSs of $G$ that are \lq most responsible' for the instability of this point (or rather its orbit), where most responsible means that a representative minimises the normalised Hilbert--Mumford weight $\mu(y,-)/|| - ||$. 
	
	The unstable strata $S_\beta$ can be inductively described following \cite{Kirwan_thesis}: first, pick a representative $\lambda_\beta$ of each conjugacy class $\beta$ (to do this canonically, fix a maximal torus $T$ and positive Weyl chamber $\ft_+$ in the Lie algebra $\ft$ of $T$, and then pick the representatives in $\ft_+$). For any 1-PS $\lambda$ of $G$, there is a parabolic subgroup
	\begin{equation}\label{eq parabolic}
	P(\lambda) : = \{ g \in G : \lim_{t \ra 0} \lambda(t) \: g \: \lambda(t)^{-1} \mathrm{\: exists \: in \: } G \},
	\end{equation}
	which is a semi-direct product of a reductive Levi group $L(\lambda)$ with a unipotent group $U(\lambda)$. If $\lambda_\beta$ is the chosen representative of the conjugacy class $\beta$, then we write $P_\beta = P(\lambda_\beta) =U_\beta \rtimes L_\beta$ and let $q_\beta : P_\beta \ra L_\beta$ denote the retraction onto the Levi factor.
	
	\begin{defn}
		Let $Z_\beta$ be the components of the fixed locus $Y^{\lambda_\beta(\GG_m)}$ on which the Hilbert--Mumford weight $\mu(-,\lambda_\beta)$ takes the value $-|| \lambda_\beta||^2$. Let $Y_\beta$ be the locally closed subscheme of $Y$ whose closed points are those $y \in Y$ such that $\lim_{t \ra 0} \lambda_\beta(t) \cdot y \in Z_\beta$.
	\end{defn}
	
	There is a natural retraction $p_\beta : Y_\beta \ra Z_\beta$. Moreover, $Y_\beta$ is invariant under the $P_\beta$-action, $Z_\beta$ is invariant under the $L_\beta$-action and $p_\beta$ is equivariant with respect to $q_\beta  : P_\beta \ra L_\beta$. Under the identification between (rational) characters and co-characters given by $|| - ||$, we let $\chi_{-\beta}: L_\beta \ra \GG_m$ denote the rational character corresponding to the rational 1-PS $\lambda_\beta^{-1} : \GG_m \ra Z(L_\beta)$.
	
	The subschemes $Y_\beta$ and $Z_\beta$ can be more explicitly defined in terms of weights for the action of a maximal torus $T < G$ as follows. The linearisation $\cL$ (or some positive power) determines a $G$-equivariant embedding $Y \hookrightarrow \PP^n$ and we can choose coordinates on $\PP^n$ such that the action of the fixed maximal torus $T \subset G$ is given by $t \mapsto \diag(\alpha_0(t), \dots , \alpha_n(t))$ for characters $\alpha_i : T \ra \GG_m$; then by \cite{Kirwan_thesis}, we have
	\begin{align}\label{weight desc Ybeta}
	\overline{Y_\beta} &  =  Y \times_{\PP^n} \{ [p_0: \dots :p_n] \in \PP^n : p_i = 0 \text{ if } \alpha_i \cdot \beta < ||\beta||^2 \} \nonumber \\
	{Y_\beta} & =  \overline{Y_\beta} \times_{\PP^n} \{ [p_0: \dots :p_n] \in \PP^n :  \exists i   \text{ such that } p_i \neq 0 \text{ and } \alpha_i \cdot \beta = || \beta||^2 \} \\
	Z_\beta &  =  Y \times_{\PP^n} \{ [p_0: \dots :p_n] \in \PP^n : p_i = 0 \text{ if } \alpha_i \cdot \beta \neq ||\beta||^2 \} \nonumber.
	\end{align}
	
	\begin{defn}
		Let $\cL_\beta$ denote the $L_\beta$-linearisation on $Z_\beta$ given by restricting the linearisation $\cL$ to $Z_\beta$ and twisting by the character $\chi_{-\beta}$. We define $Z_\beta^{ss} \subset Z_\beta$ to be the open subset of GIT semistable points for the $L_\beta$-action on $Z_\beta$ with respect to $\cL_\beta$. Let $Y_\beta^{ss}:= p_\beta^{-1}(Z_\beta^{ss})$.
	\end{defn}
	
	Then the following theorem collects the key results concerning this instability stratification, which we refer to as the Hesselink-Kempf-Kirwan-Ness stratification (or HKKN stratification).
	
	\begin{thm}
		Let $G$ be a reductive group acting on a projective scheme $Y$ with respect to an ample linearisation $\cL$. For a choice of norm $|| - ||$ as above, there is an associated stratification
		\[ Y = \bigsqcup_{\beta \in \cB} S_\beta \]
		into finitely many $G$-invariant locally closed subschemes $S_\beta$ with the following properties.
		\begin{enumerate}[label=\roman*)]
			\item If $Y^{ss}(\cL) \neq \emptyset$, then $\cB$ has a minimal element denoted $0$ such that $S_0 = Y^{ss}(\cL)$.
			\item The closure of a stratum $S_\beta$ is contained in the union $ S_\beta \sqcup\bigsqcup_{\beta' > \beta} S_{\beta'}$.
			\item For $\beta \in \cB$, we have $ S_\beta := G \cdot Y_\beta^{ss} \cong G \times^{P_\beta} Y_\beta^{ss}$.
		\end{enumerate}  
	\end{thm}
	
	\begin{rmk}\label{sch str}
		The description of the unstable strata in \cite{Kirwan_thesis} was originally only described for varieties. For a scheme $Y$, one uses the linearisation $\cL$ to $G$-equivariantly embed $Y$ into $\PP^n$ and then constructs a scheme theoretic stratification of $Y$ by taking the fibre product inside $\PP^n$ of $Y$ with the strata in $\PP^n$, which are locally closed subvarieties of $\PP^n$ (see \cite[$\S$4]{HK}). 
	\end{rmk}
	
\subsubsection{Example: HKKN stratification for sheaves}\label{sec instab sheaves}

For a projective polarised scheme $(B,\cO_B(1))$, recall from $\S$\ref{sec red GIT sheaves} that the moduli space of semistable sheaves on $B$ with Hilbert polynomial $P$ is constructed as a GIT quotient of the $\SL_{P(n)}$-action on the closed subscheme $R_n$ of the Quot scheme $Q_n:=\quot(\cO_B(-n)^{\oplus P(n)},P)$ linearised by $\cL_{n,m}$, provided $m >\!>n >\!> 0$.

Associated to the universal quotient sheaf $\cU_n \ra B \times Q_n$, there is a schematic Harder--Narasimhan stratification \cite{Nitsure}
	\begin{equation} \label{HN strat Quot}
	Q_n = \bigsqcup_{\tau} Q^\tau_n. 
	\end{equation}
In \cite{Hoskins,HK}, this HN stratification was compared with the HKKN
	\[ Q_n = \bigsqcup_{\beta \in \cB_{n,m}} S_\beta \]
associated to this $\SL_{P(n)}$-action on $Q_n$ linearised by $\cL_{n,m}$. The first step is to associate to a HN type $\tau$ a candidate HKKN index $\beta$ (i.e. the conjugacy class of a rational 1-PS).

\begin{defn}\label{def hess index from HN type} 
For a tuple of Hilbert polynomials $\tau = (P_1, \dots , P_l)$ with $\sum_{i=1}^l P_i =P$ and $(n,m) \in \NN^2$, we let $\beta_{n,m}(\tau)$ denote the conjugacy class of the following rational 1-PS of $\SL_{P(n)}$
		\[ \lambda_{\beta_{n,m}(\tau)}(t) = \diag (t^{\beta_1}I_{P_1(n)}, \dots , t^{\beta_l} I_{P_l(n)}) \: \: \text{ where } \:\beta_i:= \frac{P(m)}{P(n)} - \frac{P_i(m)}{P_i(n)}.\]
	\end{defn}

Let $P_\tau := P(\lambda_{\beta_{n,m}(\tau)}) = U_\tau \rtimes L_\tau < \SL_{P(n)}$ be the parabolic subgroup determined by $\lambda_{\beta_{n,m}(\tau)}$. 
	
\begin{thm}[\cite{Hoskins,HK}]\label{thm comp strat}
		Let $\tau=(P_1, \dots , P_l)$ be a HN type for sheaves on $(B,\cO_B(1))$ with Hilbert polynomial $P$; then, for $m > \!> n > \! > 0$, the following statements hold for $\beta:=[\beta_{n,m}(\tau)]$.
		\begin{enumerate}[label=\roman*)]
			\item All sheaves over $B$ with HN type $\tau$ are $n$-regular.
			\item If $S_\tau$ denotes the open subscheme of $Q_n^\tau$ consisting of quotients $q: \cO_B(-n)^{\oplus P(n)} \surj \cE$ such that $H^0(q(n))$ is an isomorphism, then the $SL_{P(n)}$-orbits of closed points in $S_\tau$ are in bijection with the isomorphism classes of sheaves over $B$ with HN type $\tau$.
			\item $\beta \in \cB$ is a HKKN index for the $\SL_{P(n)}$-action on $Q_n$ with respect to $\cL_{n,m}$ and $S_\tau$ is a closed subscheme of the corresponding HKKN stratum $S_{\beta}$.
			\item There is a closed $L_\beta$-invariant subscheme $Z_\tau^{ss} \subset Z_\beta^{ss}$ such that if $Y_\tau^{ss} := p_\beta^{-1}(Z_\tau^{ss})$; then
			\[ S_\tau = \SL_{P(n)} \cdot Y_\tau^{ss} \cong \SL_{P(n)} \times^{P_\tau} Y_\tau^{ss}. \]
		\end{enumerate}
	\end{thm}
	
	Moreover, by \cite{Hoskins,HK}, we have an isomorphism
	\begin{equation}\label{Ztauss}
	Z_\tau^{ss} \cong Q^{ss}_n(P_1) \times \cdots \times Q^{ss}_n(P_l),
	\end{equation}
	where $Q^{ss}_n(P_i)$ is the open subscheme of $Q_n(P_i):=\quot(\cO_B(-n)^{\oplus P_i(n)},P_i)$ consisting of quotient sheaves $q_i : \cO_B(-n)^{\oplus P_i(n)} \surj \cE_i$ such that $\cE_i$ is semistable and $H^0(q_i(n))$ is an isomorphism. In fact, if $R_n(P_i)$ denotes the closure of $Q^{ss}_n(P_i)$ in $Q_n(P_i)$; then $Z_\tau^{ss}$ is the open subscheme of
	\begin{equation}\label{Ztau}
	Z_\tau \cong R_n(P_1) \times \cdots \times R_n(P_l)
	\end{equation}
	equal to the GIT semistable set for the action of the Levi subgroup $L_{\tau} < P_\tau$ with respect to the canonical linearisation $\cL_\tau:=\cL_\beta$. The closed points of $Y_\tau^{ss}$ are quotients $q : \cO_B(-n)^{\oplus P(n)} \surj \cE$ such that $H^0(q(n))$ is an isomorphism, $\cE$ has HN type $\tau$ and  the filtration on $\CC^{P(n)}$ given by the 1-PS $\lambda_{\beta_{n,m}(\tau)}$ induces the HN filtration of $\cE$. The retraction $p_\tau : Y_\tau^{ss} \ra Z_\tau^{ss}$ sends a quotient sheaf $q : \cO_B(-n)^{\oplus P(n)} \surj \cE$ to the associated graded sheaf for its HN filtration.
	
	\begin{rmk}
If $B$ is a smooth projective curve, then for $n >\!> 0$, the Quot scheme $Q^{ss}_n(P)$ is smooth (see \cite[Proposition 2.2.8]{HL}). Consequently, for a HN type $\tau$ of a vector bundle on a smooth projective curve, the HN stratum $S_\tau$ is also smooth for $n >\!> 0$, as $Z_\tau^{ss}$ is isomorphic to a product of smooth schemes by (\ref{Ztauss}).
	\end{rmk}

	\subsection{Non-reductive GIT quotients}\label{sec non red GIT}
	
	The results of non-reductive GIT are strongest for linear algebraic groups with so-called graded unipotent radical; for linear actions of such groups on projective schemes, there are open (semi)stable sets with Hilbert--Mumford type descriptions that admit quasi-projective good quotients with natural projective completions \cite{BDHK2,BDHK}. Thus the key features and computational flexibility of reductive GIT carry over to this non-reductive setting. In this subsection, we summarise and slightly expand on these results.
	
	\subsubsection{Linear algebraic groups with graded unipotent radical}
	
	Throughout this section, we let $H = U \rtimes L$ be a linear algebraic group, where $U$ is the unipotent radical of $H$ and $L$ is a reductive Levi subgroup (which is unique up to conjugation). We only consider linear algebraic groups whose unipotent radical is graded in the following sense.
	
	\begin{defn}\label{def graded unipotent radical}
		An \emph{internal grading of the unipotent radical} $U$ of $H$ is a 1-PS $\lambda : \GG_m \ra Z(L)$ such that $\lambda(\GG_m)$ acts via conjugation on the Lie algebra of $U$ with strictly positive weights. 
		If there is a central 1-PS $\lambda : \GG_m \ra Z(L)$ which grades the unipotent radical $U$ of $H$ via its conjugation action, we say that $\lambda$ \emph{internally grades the unipotent radical} $U$ of $H$, and define the normal subgroup \[\hU  :=  U \rtimes \lambda(\Gm).\] Thus $\hU <H$ with quotient $\overline{L}:=L/\lambda(\GG_m)$.
	\end{defn}
	
For us, the key example of a group with internally graded unipotent radical is a parabolic subgroup $P = U \rtimes L$ of a reductive group (see Lemma \ref{lemma parabolic graded}). 

Note that the grading 1-PS is not necessarily unique. We will see that a choice of grading 1-PS together with a choice of linearisation determines an open subset admitting a quotient. The dependence of this open set on the grading and linearisation is analogous to reductive Variation of GIT \cite{DH} \cite{Thaddeus}, and is explored in \cite{Berczi2018a}.
	
\begin{rmk} \label{external grading}
In this paper, we will be primarily concerned with groups with internally graded unipotent radical and so we describe the results of \cite{BDHK} in this setting. There is however a more general notion of groups with \emph{externally graded unipotent radical}, as follows. An \emph{ external grading of the unipotent radical} $U$ of $H$ is a 1-PS $\lambda : \GG_m \ra \Aut(H)$ not coming from a 1-PS of $H$, which commutes with the action of $L$ and acts on $U$ with strictly positive weights. One can then construct the associated \emph{externally graded group} $\hH$ as the semi-direct product of $H$ and this 1-PS $\lambda$
\[ \hH := H \rtimes \lambda(\GG_m) = U \rtimes (L \times \lambda(\GG_m)) = \hU \rtimes L \: \text{ where } \hU:= U \rtimes \lambda(\GG_m). \]

 The general procedure for taking quotients using external gradings is described in \cite{BDHK,BDHK2}.

\end{rmk}

Often in this paper we will choose a subgroup $L' < L$ such that the composition $L' \ra L \ra \overline{L} = L/\lambda(\GG_m)$ is surjective with finite kernel, and then think of $\lambda$ as externally grading the subgroup $H' = U \rtimes L' < H = U \rtimes L$; then the quotients of $Y$ by the internally graded group $H$ and externally graded group $\hH'$ coincide.

	\subsubsection{Set-up and notation}
	
	\begin{ass}
		Throughout $\S$\ref{sec non red GIT}, we assume that $H= U \rtimes L$ is a linear algebraic group which is internally graded by a fixed 1-PS $\lambda: \GG_m \ra Z(L)$ and that $H$ acts on a projective scheme $Y$ with respect to a very ample $H$-linearisation $\cL \ra Y$. We write $\hat{U}: = U \rtimes \lambda(\GG_m) < H$ and $\overline{L}:= L/\lambda(\GG_m)$ for the quotient of $L$ by this central $\GG_m$.  
	\end{ass}
	
	Following \cite{BDHK}, we will first construct a quotient of the action of the graded unipotent group $\hU$, and then use reductive GIT to construct the residual quotient by $\overline{L}$. In $\S$\ref{sec star cond} and $\S$\ref{sec nrGIT blowups}, we explain the central results in \cite{BDHK,BDHK2}, which give the construction of a quasi-projective geometric quotient of the $H$-action on an open stable set of $Y$ admitting a Hilbert--Mumford type description together with a natural projective completion. The first step is to show that one can construct a projective quotient of the $H$-action if one assumes the stabilisers are reasonably well-behaved (see Definition \ref{def star cond} below, which one can think of as \lq semistability coinciding with stability' condition); we describe this case in $\S$\ref{sec star cond} below (\textit{cf.}\ Theorems \ref{mainthmBDHK} and \ref{mainthm nrGIT quotient}). To prove the more general statement without assuming this condition on the stabilisers, one performs a sequence of equivariant blow-ups to arrange for this condition to hold on the resulting blow-up and then we construct a quotient of an open set in $Y$ as an open subset of the projective quotient of the blow-up; this is outlined in $\S$\ref{sec nrGIT blowups}.
	
	To describe the stable loci and the relevant stabiliser conditions, we need to introduce some important notation.
	
	\begin{defn}\label{def Zmin}
		For the 1-PS $\lambda$ grading $U$, we let $\weight_{\min} < \weight_{\min +1} < \cdots < \weight_{\max} $ denote the weights of the $\lambda(\GG_m)$-action on $V:= H^0(Y,\cL)^*$ and we write $V_{\min}$ for the weight space of $\weight_{\min}$. 
	\begin{enumerate}
		\item The \emph{$\lambda$-minimal weight space}, denoted $Z_{\min}= Z(Y,\lambda)_{\min}$ is the closed subscheme \[
		Z(Y,\lambda)_{\min}:=Y \cap \PP(V_{\min})=\left\{ y \in Y^{\lambda(\GG_m)} : \lambda(\GG_m) \text{ acts on } \cL^*|_y \text{ with weight } \weight_{\min} \right\},
		\]
		\item The \emph{$\lambda$-attracting open set}, denoted\footnote{In \cite{BDHK2}, the open scheme $Y_{\min}$ is denoted $Y_{\min}^0$.} $Y_{\min}= Y(\lambda)_{\min}$, is the open subscheme
		\[
		Y(\lambda)_{\min}:=\{y\in Y \mid p(y)  \in Z_{\min}\}  \quad \mbox{and } \quad  p(y) :=  \lim_{\substack{ t \to 0\\ t \in \GG_m }} \lambda(t) \cdot y \quad \mbox{ for } y \in Y,
		\]   
		\item The \emph{$\lambda$-retraction} is $p : Y_{\min} \ra Z_{\min}$, given by $p(y):= \lim_{t \ra 0} \lambda(t) \cdot y$,
		\item The \emph{(semi)stable minimal weight space} is $Z_{\min}^{(s)s}:=Z_{\min}^{\overline{L}-(s)s}$ and we let $Y_{\min}^{(s)s} := p^{-1}(Z_{\min}^{(s)s})$.	
		\end{enumerate}	
	\end{defn}
	
	The subscheme $Z_{\min} \subset Y$ is closed and is contained in the fixed locus $Y^{\lambda(\GG_m)}$. The subscheme $Y_{\min}$ is the open stratum in the Bia{\l}ynicki-Birula decomposition for the downward flow of the $\lambda(\GG_m)$-action on $Y$. The retraction $p : Y_{\min} \ra Z_{\min}$ is a Zariski locally trivial fibration in affine varieties (or even affine spaces if $Y$ is smooth) \cite{BB}.  We assume that there are at least two distinct $\lambda(\GG_m)$-weights on $V$, as otherwise the $U$-action on $Y$ is trivial, and in this case we could consider GIT for the action of the reductive group $L = H/U$.
	
	\subsubsection{Adapted linearisations}
	
	The choice of linearisation in NRGIT is very important. In order to achieve the best results, we will later require that the linearisation satisfies the following \emph{adaptedness} condition.

\begin{defn}\label{def well adapted}
	A very ample $H$-linearisation $\cL$ over $Y$ is \emph{adapted} for the action of $\hU < H$ if the $\lambda(\GG_m)$-weights on $V:= H^0(Y,\cL)^*$ satisfy the following inequalities
	\[ \weight_{\min} < 0 <  \weight_{\min +1} < \cdots < \weight_{\max}. \]
\end{defn}

By twisting a $H$-linearisation $\cL$ over $Y$ by a (rational) character $\chi : H \ra \GG_m$, one effectively shifts the $\lambda(\GG_m)$-weights on $V$ above, and thus one can arrange for the twisted linearisation to be adapted. We note that the subschemes $Z_{\min}$ and $Y_{\min}$ do not change under this twisting.

\brm In \cite{BDHK2}, there is also a notion of \emph{well-adapted} linearisation which one can think of as requiring $0 \in (\omega_{\min},\omega_{\min} + \epsilon)$ (see \cite[$\S$2]{BDHK2} for the full definition). \erm 

For our major results we will always assume that the linearisation in question is adapted. However, the ancillary notion of \emph{borderline} linearisations\footnote{Referred to as a canonical linearisation in \cite[$\S$3.2]{HK}.} is also useful.
	
	\begin{defn}\label{def borderline adapted} 
		A very ample $H$-linearisation $\cL$ over $Y$ is \emph{borderline} for the action of $\hU < H$ if the $\lambda(\GG_m)$-weights on $V:= H^0(Y,\cL)^*$ satisfy the following inequalities
		\[ \weight_{\min} = 0 <  \weight_{\min +1} < \cdots < \weight_{\max}. \]
	\end{defn}
	
	 The significance of a borderline linearisation is that it allows one to obtain non-reductive GIT quotients without needing to impose any of the conditions on stabiliser groups discussed in \S\ref{subsubsec various stab assumptions} below.  Unfortunately, these quotients are very far from being geometric, and so are of limited practical use: if one works with such a borderline linearisation, one can obtain a categorical $H$-quotient of a larger semistable set, but this $H$-quotient factors via the retraction $p : Y_{\min} \ra Z_{\min}$ given by flowing under the 1-PS $\lambda$ and so contracts many orbits. More concretely, we have the following result, whose proof follows analogously to the argument in \cite[$\S$3.2]{HK}.
	
	\begin{prop}\label{prop canonical linearisation}
		Let $H = U \rtimes L$ be a linear algebraic group whose unipotent radical $U$ is internally graded by $\lambda: \GG_m \ra L$, acting on a projective scheme $Y$ with respect to a very ample borderline linearisation $\cL_{0}$. Then
		\[  \bigoplus_{r \geq 0} H^0(Y,\cL_{0}^{\otimes r} )^H \cong \bigoplus_{r \geq 0} H^0(Z_{\min},\cL_{0}^{\otimes r} )^L \]
		and the associated projective scheme $Z_{\min}/\!/_{\cL_0} L$ is a categorical quotient of 
		\begin{enumerate}
			\item the $L$-action on $Z_{\min}^{ss}$,
			\item the $L$-action on $Y_{\min}^{ss}$,
			\item the $H$-action on $Y_{\min}^{ss}$.
		\end{enumerate}
		In particular, we have that
		\[ Y_{\min}^{ss} = Y^{H-ss}(\cL_0) :=\{ y \in Y : \exists \sigma \in H^0(Y,\cL_0^{\otimes r})^H \text{ for } r > 0 \text{ with } \sigma(y) \neq 0 \} \]
		and the morphism $Y_{\min}^{ss} \stackrel{p}{\ra} Z_{\min}^{ss} \ra Z_{\min}/\!/_{\cL_0} L$ is a surjective categorical $H$-quotient that identifies every $y \in Y_{\min}^{ss}$ with its limit $p(y) = \lim_{t \ra 0} \lambda(t) \cdot y \in Z_{\min}^{ss}$. 
	\end{prop}	
	\begin{proof} 
		This follows by adapting the arguments in \cite[Lemma 3.1 and Proposition 3.2]{HK} for the special case when $Y_{\min} = Y_\beta$ and $Z_{\min} = Z_{\beta}$ to the retraction $p : Y_{\min}^{ss} \ra Z_{\min}^{ss}$ which is equivariant with respect to the surjection $q : H \twoheadrightarrow L$. 
	\end{proof}
	
	Another use of borderline linearisations is that they enable us to obtain information about adapated linearisations, via variation of GIT style arguments. A well-adapted linearisation is a small perturbation of a borderline linearisation $\cL_0$; typically they are constructed by twisting $\cL_0$ by a small rational character. If the associated ring of invariant sections for the perturbed linearisation are still finitely generated, then the following result compares the associated \lq quotient morphisms' (which, for the perturbed linearisation, is not necessarily surjective even when the invariant sections are finitely generated).
	
	\begin{prop}\label{prop ss locus of perturbation of canonical linearisation}
		Let $H= U \rtimes L$ be a linear algebraic group with unipotent radical $U$ internally graded by $\lambda: \GG_m \ra L$ acting on a projective scheme $Y$ with respect to a very ample borderline linearisation $\cL_{0}$. Given a (rational) character $\chi : H \ra \GG_m$, we let $\cL_\chi$ denote the twist of $\cL_0$ by the rational character $\epsilon \chi$ for $\epsilon \in \QQ_{>0}$ sufficiently small and let $R(Y,\cL_\chi) := \bigoplus_{r \geq 0} H^0(Y,\cL_{\chi}^{\otimes r} )$. If $R(Y,\cL_\chi)^H$ is finitely generated, then there is an induced, not necessarily surjective, $H$-invariant morphism 
		\[ q_{H,\chi} : Y^{H-ss}(\cL_\chi) \longrightarrow  Y/\!/_{\cL_\chi} H:= \proj (\bigoplus_{r \geq 0} H^0(Y,\cL_{\chi}^{\otimes r} )^H), \]
		where we define  $Y^{H-ss}(\cM) := \{ y \in Y : \exists \sigma \in H^0(Y,\cM^{\otimes r})^H \text{ for } r > 0 \text{ with } \sigma(y) \neq 0 \}$ for any linearisation, and we have \[Y^{H-ss}(\cL_\chi) \subset Y^{H-ss}(\cL_0) = Y_{\min}^{ss}.\] Furthermore, if $q_{H,\chi}$ is a categorical $H$-quotient, then there is an induced morphism between the $H$-quotients with respect to $\cL_\chi$ and $\cL_0$:
		\begin{equation}\label{VGIT type map to borderline linearisation}
			Y/\!/_{\cL_\chi} H \longrightarrow Y/\!/_{\cL_0} H.
		\end{equation} 
	\end{prop}
	\begin{proof}
		The inclusions of invariant rings, which are by assumption finitely generated,
		\[ R(Y,\cL_\chi) \supset R(Y,\cL_\chi)^L \supset R(Y,\cL_\chi)^H \]
		induce rational morphisms of projective schemes
		\[ q_{H,\chi} : Y \dashrightarrow Y/\!/_{\cL_{\chi}} L \dashrightarrow Y/\!/_{\cL_{\chi}} H. \]
		Note that there is not a residual action of $U$ on $Y/\!/_{\cL_{\chi}} L$. Nevertheless, we obtain that
		\[ Y^{H-ss}(\cL_\chi) \subset Y^{L-ss}(\cL_\chi) \subset Y^{L-ss}(\cL_0) = Y^{H-ss}(\cL_0) = Y^{ss}_{\min} \]
		where the first inclusion follows from the above factorisation of $q_{H,\chi}$, the second inclusion follows by variation of reductive GIT quotients (as $\cL_\chi$ is a small perturbation of $\cL_0$), and the remaining equalities follow from Proposition \ref{prop canonical linearisation}. Since the composition $Y^{H-ss}(\cL_\chi) \hookrightarrow Y^{H-ss}(\cL_0) \ra Y/\!/_{\cL_0} H$ is $H$-invariant, provided $q_{H,\chi}$ is a categorical $H$-quotient, we obtain via its universal property a morphism as in Equation \eqref{VGIT type map to borderline linearisation}.
	\end{proof}
	
	\begin{rmk}
		In fact, in the above situation, we have $ Y^{H-ss}(\cL_\chi) \subset Y^{ss}_{\min} \setminus UZ^{ss}_{\min}$ provided the perturbation $\cL_\chi$ is non-trivial, as $ Y^{H-ss}(\cL_\chi) \subset Y^{\lambda(\GG_m)-ss}(\cL_\chi)$. The latter is either equal to $Y_{\min} \setminus Z_{\min}$, if  $\cL_\chi$ is adapted, or empty, if the borderline linearisation is perturbed in the other direction. 
	\end{rmk}

	\subsubsection{Stabiliser assumptions} \label{subsubsec various stab assumptions}

	Let us introduce some assumptions on the stabiliser groups, considered as analogous to \lq semistability coinciding with (Mumford) stability for some group'. For any $y \in Y$ and any subgroup $G <H$, we write $\stab_G(y)$ for the stabiliser group of $y$.
	
	\begin{defn} \label{def star cond}
		For the $H$-action on $Y$ with respect to an ample adapted linearisation $\cL$, we define the following conditions:
		\begin{enumerate}
			\item \emph{Semistability coincides with stability for $\hU$} 
			\begin{equation}\label{starU0}
			\dim \stab_U(z) = 0 \: \text{ for all } z \in Z_{\min}, \tag*{$[\hU]_0$}
			\end{equation}
			and \emph{semistability coincides with stability for $\hU$ on the reductive semistable locus} if
			\begin{equation}\label{starU0Rss}
			\dim \stab_U(z) = 0 \: \text{ for all } z \in Z_{\min}^{ss}. \tag*{$[\hU; \overline{L}\text{-ss}]_0$} 
			\end{equation}
			\item  If $U$ is abelian, we say \emph{semistability coincides with Mumford stability for $\hU$} if 
			\begin{equation}\label{starU}
			\dim \stab_{U}(-) \text{ is constant on } Y_{\min}, \tag*{$[\hU]$}
			\end{equation}
			and \emph{semistability coincides with Mumford stability for $\hU$ on the reductive semistable locus} if
			\begin{equation}\label{starURss}
			\dim \stab_{U}(-) \text{ is constant on } Y_{\min}^{ss}. \tag*{$[\hU;\overline{L}\text{-ss}]$}
			\end{equation}
				\item \emph{Semistability coincides with stability for $\overline{L}$}
			\begin{equation}\label{starR0}
			\dim \stab_{\overline{L}}(z) = 0 \: \text{ for all } z \in Z_{\min}^{ss}. \tag*{$[\overline{L}]_0$}
			\end{equation}
			\item If $U$ is abelian, we say \emph{the minimal weight space has minimal dimensional $U$-stabilisers} if
			\begin{equation}\label{CZminNE}
				 \dim \stab_{U}(z) = d_{\min}:=d_{\min}(Y_{\min},U) \text{ for some}\footnote{By semicontinuity of stabiliser dimensions, if this holds for some $z$, it holds generically on $Z_{\min}$ (and $Z_{\min}^{ss}$).} \: z\in Z_{\min}, 
				 \tag*{$[C(Z_{\min},U,d_{\min}) \neq \emptyset]$}
			\end{equation}
			where $d_{\min}(Y_{\min},U)$ is as in Definition \ref{def centres of blowups}; this is equivalent to $C(Z_{\min},U,d_{\min})\neq \emptyset$.
\end{enumerate}			
If $U$ is not abelian then (2) and (4) above are modified by choosing a series of subgroups $ \{e\}=U^{(0)} \leqslant U^{(1)} \cdots \leqslant U^{(r)}= U$ normalised by $H$ whose successive quotients are abelian and asking for these conditions to hold for each $U^{(i)}$.
\end{defn}
			
\begin{rmk} 
Definitions (1) -- (2) on the unipotent stabilisers appear in \cite{BDHK2,BDHK}, whereas (3) is an additional condition that we use in this paper to determine the stable locus explicitly. Condition (4) is a weakened version of (2): when (2) fails, we will perform a blow-up sequence to arrange for (2) to hold, and (4) ensures that we never blow up all of $Z_{\min}$ and so we can easily determine the minimal weight space in the blow-up and obtain an explicit expression for the stable set in this case.
\end{rmk}

Our primary motivation for constructing quotients by parabolic group actions in stages is that conditions such as these may hold for the relevant groups at each stage in the quotienting process, even when they fail to hold directly for the initial $H$-action on $Y$ graded by $\lambda$. For example, when looking at sheaves of fixed Harder--Narasimhan type of length $l > 2$, Condition \ref{starR0} fails directly (see $\S$\ref{sec sheaves nonred git setup}) and our approach via Quotienting-in-stages circumvents this problem by enabling us to set up a sequence of quotients such that these types conditions may hold at each stage.

\begin{rmk}\
		\begin{enumerate}
			\item Condition \ref{starU0} is referred to as \lq semistability coincides with stability' in \cite{BDHK} and is denoted by $(ss = s \neq \emptyset[\hU])$ in \cite{BDHK2}. As a unipotent group has no non-trivial finite subgroups, \ref{starU0} is equivalent to asking that all $U$-stabilisers on $Z_{\min}$ are trivial. Furthermore, as the $\GG_m$-action normalises the $U$-action, \ref{starU0} is equivalent to $U$ acting freely on $Y_{\min}$ (see \cite[Remark 5.7]{BDHK2}). This condition depends on the grading $\hU$ of $U$ by $\lambda$, as the 1-PS $\lambda$ determines $Z_{\min}$.
			\item Condition \ref{starR0} says that for $\overline{L}$ acting on $Z_{\min}$, the GIT stable set coincides with the GIT semistable set. 
		
			\item Conditions \ref{starU} and \ref{starURss} allow for the presence of positive dimensional unipotent stabilisers, provided they are of constant dimension for each unipotent subgroup in the chosen series. See \cite[Remark 7.12]{BDHK2} and the proof of Theorem \ref{mainthmBDHK} below.
			
		\end{enumerate}
	\end{rmk}
	
Let us finally include a result about stabilisers from \cite{BDHK2} that we will need later on.
	
\begin{lemma}\cite[Lemma 5.2]{BDHK2}\label{lem second small lemma}
For the $H$-action on $Y$ with respect to an ample adapted linearisation $\cL$, we have for $z \in Z_{\min}$ and $u\in U$ that $uz \in Z_{\min}$ if and only if $u\in \stab_U(z)$.  Furthermore, \[\stab_H(z) = \stab_U(z)\rtimes \stab_L(z)\] and, in particular, $\dim\stab_H(z) = \dim\stab_U(z)+\dim\stab_L(z)$.
\end{lemma}
\begin{proof}
We extend the argument of \cite[Lemma 5.2]{BDHK2}. For the first statement, as $p$ is $U$-invariant and restricts to the identity on $Z_{\min}$, if $uz \in Z_{\min}$, we see $p(uz) = uz$ and also $p(uz) = z$, and so $u\in \stab_U(z)$. For the second statement, any $h \in H$ can be written as $h=ul$ for $u\in U$ and $l\in L$. If $hz= z$ then $uz = l^{-1}z \in Z_{\min}$ since $L$ preserves $Z_{\min}$ as the grading 1-PS $\lambda$ is central in $L$. Thus any $h\in \stab_H(z)$ can be written as a product of $ul$ with both elements fixing $z$. Since $\stab_U(z)\cap\stab_L(z)= \{e\}$ the result follows.
\end{proof}
	
	\subsubsection{The construction of quotients under optimal stabiliser assumptions}\label{sec star cond}
	
	In this section, we assume we have an adapted linear action of $H = U \rtimes L$ on a projective scheme $Y$ with internally graded unipotent radical. Let us define the stable loci in the case when the weaker assumption \ref{starU} on unipotent stabilisers holds.
	
	\begin{defn}\label{def nonred ss}
		For a $H$-action on $Y$ with respect to a very ample adapted linearisation $\cL$ satisfying \ref{starU}, we define the following stable sets.
		\begin{enumerate}
			\item The $\hU$-stable set is defined to be 
			\[ Y^{\hU-s}:= Y_{\min} \setminus UZ_{\min} = \bigcap_{u \in U} u Y^{\lambda(\GG_m)-s}.\]
			\item The $H$-stable set is defined to be
			\[ Y^{H-s}:= Y^{s}_{\min} \setminus U Z_{\min}^s. \]
		\end{enumerate}
	\end{defn}
	
	\begin{rmk}\label{rmk nrGIT ss def}\
		\begin{enumerate}
			\item The equality in the first definition holds as $Y^{\lambda(\GG_m)-s} = Y^{\lambda(\GG_m)-ss}= Y_{\min} \setminus Z_{\min}$ by the Hilbert--Mumford criterion for $\lambda(\GG_m)$ and the assumption that the linearisation is adapted. This is an open set in $Y$, as the assumption \ref{starU0} implies that $U Z_{\min} \subset Y_{\min}^0$ is a closed subscheme (\textit{cf.}\ \cite[Lemma 3.4]{BDHK2}). 
			\item In \cite{BDHK}, the $\hU$-stable set is denoted $Y^{\hU-s}_{\min +}$ rather than $Y^{\hU-s}$, but we have simplified the notation as we will not introduce any other stable sets (see \cite{BDHK,BDHK2}\ for a discussion of various other notions of stability defined using invariant sections).
			\item Our definition of the $H$-stable set differs slightly from that given in \cite{BDHK2}. However, under the assumption \ref{starR0}, both these notions coincide (see Theorem \ref{thm nrGIT HM type inclusions} below).			
			\item These definitions are functorial with respect to equivariant closed immersions $Z \hookrightarrow Y$; that is, $Z^{H-s} = Z \times_Y Y^{H-s}$. 
		\end{enumerate}
	\end{rmk}
	
Let us state the main result concerning $\hU$-quotients from	\cite{BDHK}, whose assumptions are weakened in \cite{BDHK2}.
	
\begin{thm}\label{mainthmBDHK} \cite[Theorem 0.2 and Corollary 0.6]{BDHK}
Let $H= U \rtimes L$ be a linear algebraic group with unipotent radical $U$ internally graded by $\lambda: \GG_m \ra L$ acting on a projective scheme $Y$ with respect to a very ample adapted linearisation $\cL$. Assuming condition \ref{starU} holds, the following statements hold.
\begin{enumerate}[label=\emph{(\roman*)}]
\item The $\hU$-action on $Y^{\hU-s}= Y_{\min} \setminus U Z_{\min}$ has a projective geometric $\hU$-quotient \[q_{\hU}: Y^{\hU-s} \ra Y/\!/\hU = \Proj R(Y,\cL)^{\hU}.\]   
\item By taking a reductive GIT quotient of the induced $\overline{L}$-action on the projective geometric $\hU$-quotient $Y/\!/\hU$
\[ \xymatrixcolsep{4pc} \xymatrix{ q: Y^{\hU-s} =Y_{\min} \setminus UZ_{\min}  \ar@{->>}[r]^{\quad \quad \quad q_{\hU}} & Y/\!/\hU \ar@{-->}[r]^{q_{\overline{L}}\quad \quad \quad} & Y/\!/H := (Y/\!/\hU)/\!/\overline{L}}\]
one obtains a projective good $H$-quotient of an open set of $Y^{\hU-s}$, \[\dom(q) \twoheadrightarrow Y/\!/H = \Proj R(Y,\cL)^H.\]

\end{enumerate}
\end{thm}
\begin{proof}
Under the stronger assumption \ref{starU0} this is \cite[Theorem 0.2 and Corollary 0.6]{BDHK}; let us describe how to weaken this assumption to \ref{starU} as outlined in \cite[Remarks 2.8 and 7.12]{BDHK2}. The idea is to iteratively construct quotients by the abelian groups $U^{(j)}/U^{(j-1)}$ for $j =1 , \dots ,r$ by locally choosing complementary subgroups to the stabiliser groups. To illustrate the argument, we assume for simplicity that $U$ is abelian; then by \ref{starU}, there exists $d$ such that $\dim \stab_U(y) =d$ for all $y \in Y_{\min}$. As $U$ is an abelian unipotent group, for each stabiliser group $\stab_U(y)$ of $y \in Y_{\min}$, we can choose a complementary subgroup $U'$ with $\stab_U(y)U' = U$ and $\stab_U(y)\cap U = \{1\}$, by choosing a complementary subspace in the Lie algebra. Taking the Lie algebra of the $U$-stabiliser gives a $\hU$-invariant function $Y_{\min} \ra \mathrm{Gr}(d,\Lie(U))$ and being complementary to a given $U' < U$ is an open condition in this Grassmannian. Therefore, we get a corresponding open $\hU$-invariant set of $Y_{\min}$ on which all stabilisers have complementary subgroup $U'$ and thus quotienting by $U$ is equivalent to quotienting by $U'$. This open set of $Y_{\min}$ admits a geometric $U'$-quotient (and thus $U$-quotient) by \cite[Proposition 7.4]{BDHK2}. Then by varying $U'$ we can cover $Y_{\min}$ by open sets admitting geometric $U$-quotients in order to obtain a geometric $U$-quotient of $Y_{\min}$. From there one proceeds as in the proof of the case when \ref{starU0} holds.
\end{proof}

\begin{rmk}\label{rmk Uhat part of proj completion constr}
Let us make a few comments about the construction of the quotients in the above theorem.
\begin{enumerate}
\item One first shows that there is a geometric $U$-quotient $q_U: Y_{\min} \ra Y_{\min}/U$ by gluing together affine geometric quotients; this uses the grading 1-PS $\lambda$ and the fact that \ref{starU0} holds (see \cite[Proposition 7.4]{BDHK2}).
\item Then the grading 1-PS $\lambda$ is also used to give a projective $\hU$-quotient as follows. First one constructs a $\lambda(\GG_m)$-equivariant embedding of $Y_{\min}/U$ in a projective space and we let $\overline{Y_{\min}/U}$ denote its projective completion (see \cite[Lemma 7.6]{BDHK2}). Then the projective $\hU$-quotient is constructed by taking a reductive quotient of $\overline{Y_{\min}/U}$ by $\lambda(\GG_m)$, where one slightly perturbs the adapted linearisation by twisting by a rational character to obtain a so-called well-adapted linearisation (\textit{cf.}\ \cite[$\S$2]{BDHK2}). 
\item If \ref{starU0} holds, then by \cite[Theorem 0.2 (ii)]{BDHK}, $Y^{\hU-s}$ coincides with the \emph{locally trivial stable locus} (see \cite[Definition 1.5)]{BDHK}, and it is shown that this geometric $\hU$-quotient is projective and coincides with the \emph{enveloping quotient} (see \cite[Definition 1.5]{BDHK}). Hence, the ring of $\hU$-invariant sections of a sufficiently divisible tensor power of this well-adapted linearisation is finitely generated and its projective spectrum is $Y^{\hU-s}/\hU=Y/\!/\hU$. In particular, $q_{\hU}$ can be constructed by taking invariant sections for a slight perturbation of the original linearisation.
\item \label{rmk when UZ equals Ymin}  In general we will assume $U Z_{\min} \neq Y_{\min}$. If this does not hold then the $\hU$-stable set $Y^{\hU-s}$ is empty and in this case, we should define $Y_{\min}$ as a replacement stable set which admits a the geometric $p: Y_{\min} \ra Z_{\min}$ (see \cite[Remark 2.5]{BDHK2}). Moreover, there is a good $H$-quotient $UZ_{\min}^{ss} \stackrel{p}{\ra} Z_{\min}^{ss} \ra Z_{\min} /\!/_{\! \cL} \overline{L}$.
\end{enumerate}
\end{rmk}

Determining the domain of the quotient $q$ in terms of the $H$-action on $Y$ is quite delicate, even though we can determine the domain of $q_{\overline{L}}$ using the Hilbert--Mumford criterion for the reductive group $\overline{L}$. If $T$ is a maximal torus of $L$ containing $\lambda(\GG_m)$, one would like a description of $\dom(q)$ in terms of $T$-weights on $Y$, but passing to the $U$-quotient involves deleting some of these weights. Nevertheless, we can sandwich the domain of definition between sets admitting Hilbert--Mumford type descriptions in terms of weights as in Theorem \ref{thm nrGIT HM type inclusions} \emph{\ref{HM type inclusions}} below. Moreover, under the assumption \ref{starR0}, we can explicitly determine the domain of $q$. In \cite[Lemma 7.8]{BDHK2}, a Hilbert--Mumford type statement is also stated without assuming \ref{starR0}.

\begin{thm}\label{thm nrGIT HM type inclusions}
Let $H = U \rtimes L$ be a linear algebraic group with unipotent radical $U$ internally graded by $\lambda: \GG_m \ra L$ acting on a projective scheme $Y$ with respect to a very ample adapted linearisation $\cL$ such that condition \ref{starU} holds. Let $q : Y \dashrightarrow Y/\!/H$ be the non-reductive GIT quotient given by Theorem \ref{mainthmBDHK}.
\begin{enumerate}[label=\emph{(\roman*)}]
\item\label{nrGIT HM inclusion 1}  If $y\in Y^s_{\min}\setminus U Z^s_{\min}$, then $q_{U}(y) \in Y_{\min}/U$ is stable for the reductive group $L$ acting on $\overbar{Y_{\min}/U}$. 
\item\label{HM type inclusions} The domain of $q$ sits between the open sets with Hilbert--Mumford type descriptions
\[  Y^{H-s} = Y_{\min}^s \setminus UZ_{\min}^s  \subset \dom(q) \subset \bigcap_{h \in H} hY^{T-ss}  =  \bigcap_{u \in U} uY^{L-ss}  \subset Y_{\min}^{ss} \setminus UZ_{\min}^{ss}. \]
\item The open set $Y^{H-s}$ admits a quasi-projective geometric $H$-quotient by restricting $q$
\[  q|_{Y^{H-s}} : Y^{H-s} = Y_{\min}^s \setminus UZ_{\min}^s \twoheadrightarrow q(Y^{H-s} ) \subset Y/\!/H,\] 
with a natural projective completion $Y/\!/H$.
\end{enumerate}
In particular, if one also assumes \ref{starR0}, then all the above inclusions in \ref{HM type inclusions} are equalities and so $q : Y^{H-s} \ra Y/\!/H$ is a projective geometric $H$-quotient.
\end{thm}
\begin{proof}
After replacing the linearisation with a sufficiently high tensor power (which does not change the notion of semistability), the quotient map $q_U$ is constructed from a linear projection $\PP(V) \dashrightarrow \PP(W)$ where $V= H^0(Y,\cL)^*$ and $W = (H^0(Y,\cL)^U)^*$. More precisely, we obtain a locally closed immersion $Y_{\min}/U \hookrightarrow \PP(W)$ and we let $\overbar{Y_{\min}/U} \subset \PP(W)$ denote the projective completion. There is an induced $L$-action on both of these projective schemes, whose (semi)stable loci can be described in terms of the (reductive GIT) Hilbert--Mumford criterion. 

Let $T$ be a maximal torus in $L$ containing $\lambda(\GG_m)$ and consider the $T$-weights on $V$ and $W$. For $y \in Y$ (or strictly speaking a lift to the affine cone $V$), we need to compare the convex hull of $T$-weights of $y$ and $q_U(y)$; since $q_U$ is induced by a linear projection, it forgets some weights:
\[ \conv_T(y) \supset \conv_T(q_U(y)). \]
However, the sections in $H^0(Y,\cL)$ with maximal $\lambda$-weights (dually corresponding to the locus $V_{\min} \subset V =H^0(Y,\cL)^*$ with minimal $\lambda$-weights) are all $U$-invariant (see \cite[$\S$7]{BDHK2}) and thus these weights are preserved by $q_U$. Let $H_{\min} \subset \mathfrak{t}^*$ be the hyperplane containing the $T$-weights on $V$ which are $\lambda$-minimal (this hyperplane is perpendicular to the ray corresponding to $\lambda$ when we identify $\mathfrak{t} \cong \mathfrak{t}^*$ using our given norm); then we have
\[ \conv_T(y) \cap H_{\min} = \conv_T(q_U(y)) \cap H_{\min}.\]
All the $T$-weights are contained in a closed half space bounded by this hyperplane and, as the linearisation is twisted to be well-adapted, the origin lies very close to this hyperplane. Points in $Y_{\min}$ have at least one weight on $H_{\min}$ and points in $Z_{\min}$ have all their weights on this hyperplane. The map $p : Y_{\min} \ra Z_{\min}$ is the retraction to the minimal $\lambda$-weight space $\PP(V_{\min})$, i.e.\ it throws away all the $T$-weights not lying on the hyperplane $H_{\min}$. If $y \in Y_{\min} \setminus UZ_{\min}$, then the weight polytope $\conv_T(q_U(y))$ contains a weight outside of $H_{\min}$. If $y\in Y^s_{\min}$, then $p(y) \in Z^s_{\min}$ is $\overline{L}$-stable; that is, the projection of the origin to $H_{\min}$ is contained inside the interior of $\conv_T(y) \cap H_{\min}$. 

For the proof of (i), if $y\in Y^s_{\min}\setminus U Z^s_{\min}$, then we claim that the origin is contained in the interior of $\conv_T(q_U(y))$. Indeed, the origin lies very close to $H_{\min}$, and since $y\in Y^s_{\min}$ this weight polytope meets $H_{\min}$ in a codimension 1 polytope $\conv_T(y) \cap H_{\min}$ containing the projection of the origin. Since $y \notin UZ_{\min}$, the weight polytope also contains one weight away from $H_{\min}$, so it is full dimensional. Since $Y^s_{\min}\setminus U Z^s_{\min}$ is $H$-invariant, the same is true for all points in the $H$-orbit of $y$, whose image under $q_U$ is the $L$-orbit of $y$; hence, by the Hilbert--Mumford criterion for the $L$-action on $\PP(W)$, we conclude that $q_U(y)$ is $L$-stable.

For (ii), the first inclusion follows from (i), as $\dom(q) = q_U^{-1}(\overbar{Y_{\min}/U}^{L-ss})$. The second inclusion follows as we can pullback invariant sections along $q_U$. The middle equality is just the reductive Hilbert--Mumford criterion and the right inclusion holds by a Hilbert--Mumford type argument involving $T$-weights similar to that given in (i).

Finally (iii) follows from (i) and (ii), and the final statement follows as then $Z_{\min}^s = Z_{\min}^{ss}$.
\end{proof}

Finally, let us explain how one can weaken \ref{starU} to \ref{starURss} following  \cite[Remark 2.8]{BDHK2}. In this case, one minor difference is that although the $H$-quotient is projective, the $\hU$-quotient is only quasi-projective. 

\begin{thm}\label{mainthm nrGIT quotient}
For a $P$-action on $X$ with an adapted linearisation $\cL$ satisfying \ref{starURss}, there is a map
\[ \xymatrixcolsep{4pc} \xymatrix{ q: Y^{ss}_{\min} \setminus UZ^{ss}_{\min} \ar@{->>}[r]^{ q_{\hU} \quad} & (Y^{ss}_{\min} \setminus UZ^{ss}_{\min})/\hU \ar@{-->}[r]^{\quad \quad \quad q_{\overline{L}}} & Y/\!/H}, \]
giving a projective good $H$-quotient $\dom(q) \twoheadrightarrow Y/\!/H$ of an open set of $Y^{ss}_{\min} \setminus UZ^{ss}_{\min}$, which restricts to a quasi-projective geometric $H$-quotient of $Y^{H-s} = Y_{\min}^s \setminus UZ_{\min}^s$. If additionally,  \ref{starR0} holds, then $\dom(q) = Y^{H-s}$ and this has a projective geometric quotient $Y/\!/H$ given by taking the projective spectrum of the ring of $H$-invariant sections for a twist of the linearisation which is well-adapted.
\end{thm}
\begin{proof}
The argument is a minor modification of \cite[\S 7]{BDHK2}. For simplicity, let us outline the proof in the case when $U$ is abelian and acts with zero dimensional stabilisers on $Z^{ss}_{\min}$; if $U$ is not abelian, one needs to successively take quotients by the abelian groups appearing as quotients in a normal series and, if the stabilisers are positive dimensional, one needs to choose complementary subgroups as in the proof of Theorem \ref{mainthmBDHK}.

First one constructs a geometric $U$-quotient of $Y^{ss}_{\min}$ similarly to the affine local construction of \cite[Proposition 7.4]{BDHK2}; however, one instead takes a basis $\{\sigma_i\}_{i=1}^r$ of the $\overbar{L}$-invariant sections of maximal $\lambda(\GG_m)$-weight $H^0(Y,\cL)^{\overbar{L}}_{\max}$ (recall that these sections are also $U$-invariant), then $Y_{\min}^{ss}$ is covered by the open affine sets $Y_{\sigma_i}$ and each of these admits a trivial geometric $U$-quotient given by taking rings of $U$-invariants (\textit{cf.} \cite[Lemmas 7.2 and 7.3]{BDHK2}). Since these invariant rings are finitely generated, one can pick $s$ so they are all generated in degree $s$. Let $W:=(H^0(Y,\cL^{\otimes s})^U)^*$; then we have a natural $U$-invariant morphism $\phi: Y^{ss}_{\min} \ra \PP:=\PP(W)$, which thus factors via the quotient $q_U : Y^{ss}_{\min} \ra Y^{ss}_{\min}/U$ and a locally closed immersion $\overline{\phi} :Y^{ss}_{\min}/U \ra \PP$.  Moreover, there is an induced $L$-action on $\PP$ such that $\phi$ is $L$-equivariant; we denote the corresponding (semistable) $\lambda$-minimal weight space $Z^{(ss)}_{\min}(\PP)$ and the corresponding open by $\PP_{\min}^{(ss)}$. After appropriately twisting the linearisation on $\PP$ so it is well-adapted, we have $\PP^{\lambda-(s)s} = \PP_{\min} \setminus Z_{\min}(\PP)$ by the Hilbert--Mumford criterion.

The same argument as in \cite[Lemma 7.6]{BDHK2} shows that the image of $\overline{\phi}$ is contained in $\PP_{\min}^{ss}$ as a closed subscheme. Moreover, for $y \in Y_{\min}^{ss}$, one has $\phi(y) \in Z^{ss}_{\min}(\PP)$ if and only if $y \in UZ_{\min}^{ss}$ analogously to \cite[Equation (3) of $\S$7]{BDHK2}. Consequently, $\phi$ induces a closed immersion
\[ (Y_{\min}^{ss} \setminus U Z_{\min}^{ss})/U \hookrightarrow \PP_{\min}^{ss} \setminus Z_{\min}^{ss} = \PP_{\min}^{ss} \times_{\PP} \PP^{\lambda-s}.\]
Thus one obtains a geometric quotient $(Y^{ss}_{\min}\setminus UZ_{\min}^{ss})/\hU$ as a locally closed subscheme of the projective scheme $\PP/\!/ \lambda(\GG_m)$; in particular, this $\hU$-quotient is not necessarily projective.

The reductive GIT quotient $q_L: \PP \dashrightarrow \PP/\!/L$ factors via $\PP/\!/ \lambda(\GG_m)$ and we will construct our $H$-quotient as a closed subscheme of $\PP/\!/L$. A reductive Hilbert--Mumford argument shows  $\PP^{L-ss} \subset \PP^{ss}_{\min} \setminus Z^{ss}_{\min}$, and since $(Y_{\min}^{ss} \setminus U Z_{\min}^{ss})/U$ is closed in the latter, the inclusion
\[\PP^{L-ss}  \times_{\PP} (Y_{\min}^{ss} \setminus U Z_{\min}^{ss})/U \hookrightarrow \PP^{L-ss} \]
is closed and  $Y/\!/H:=q_L( \PP^{L-ss}  \times_{\PP} (Y_{\min}^{ss} \setminus U Z_{\min}^{ss})/U)$ is a closed subscheme of $\PP/\!/L$. The composition $Y_{\min}^{ss} \setminus U Z_{\min}^{ss} \ra Y/\!/H$ is thus a projective good $H$-quotient. 

If in addition, we assume that \ref{starR0} holds, then $\dom(q) = Y^{H-s}$ as in Theorem \ref{thm nrGIT HM type inclusions}. To conclude that $Y/\!/H$ is given by taking the projective spectrum of the ring of invariants for a well-adapted perturbation of the linearisation, it suffices to show this quotient coincides with the enveloping quotient by \cite[Corollary 3.1.21]{BDHK_handbook}. To prove this, one first shows that $\dom(q)$ is contained in the locally trivial stable locus with respect to $H$ (see \cite[Definition 3.3.2]{BDHK_handbook}) and then deduces that $Y/\!/H$ coincides with the quotient of the locally trivial stable locus and the enveloping quotient as in \cite[$\S$7, Proof of Theorem 2.16]{BDHK2}.\end{proof}

	\subsubsection{The construction of quotients in general via blow-up sequences}\label{sec nrGIT blowups}
	
	If condition \ref{starU} does not hold for the linearised $H$-action on $Y$, then as explained in \cite[$\S$8]{BDHK2} one can instead perform a sequence of equivariant blow-ups in order to arrange that condition \ref{starU} holds on the blow-up. These blow-ups are similar to the blow up sequences used to construct partial desingularisations of reductive GIT quotients (\textit{cf}.\ $\S$\ref{sec red part desing} and \cite{K2}), but where one now considers unipotent stabiliser groups instead. Away from the exceptional locus, we can identify the semistable locus in the blown up space $\widetilde{Y}$ with a certain open subset of $Y$; hence one obtains a quasi-projective quotient of this open subset, and a compactification given by the quotient of $\widetilde{Y}$. To identify this open subset in $Y$, one has to describe the complement to the exceptional divisor in the $\hU$-stable locus for $\widetilde{Y}$. This is more difficult when all of $Z_{\min}$ is blown up (see \cite{Josh_length2} for further discussion).
	
	There are several cases to consider, of increasing complexity:
	\begin{enumerate}
		\item If the generic dimension of $U$-stabilisers on $Z_{\min}$ is zero, then one does not need to blow up all of $Z_{\min}$ and there is an explicit description of the open subset admitting a quotient given by \cite[Theorem 2.10]{BDHK2}. 
		\item If $U$ is abelian, but the generic dimension of the $U$-stabilisers on $Z_{\min}$ is positive, then either
		\begin{enumerate}
			\item the generic dimension of the $U$-stabilisers on $Z_{\min}$ equals the generic dimension of the $U$-stabilisers on $Y_{\min}$ or
			\item the generic dimension of the $U$-stabilisers on $Z_{\min}$ is bigger than the generic dimension of the $U$-stabilisers on $Y_{\min}$.
		\end{enumerate}
		In the first case, one does not need to blow up all of $Z_{\min}$ and it is easier to describe the open subset that one obtains a quotient of; see $\S$\ref{sec easy blowups abelian} for a situation in which this is particularly straight-forward. In the second case, one must blow-up all of $Z_{\min}$ and use jets to describe the open subset admitting a quotient.
		\item For general $U$ such that the generic dimension of the $U$-stabilisers on $Z_{\min}$ is positive, the results of the blow-up procedure is described by \cite[Theorem 2.20]{BDHK2}. In this case one fixes a series of subgroups of $U$ which are normal in $H$ and whose successive quotients are abelian with $\lambda(\GG_m)$ acting by a single weight. One could then perform quotients in stages doing blow-ups at each stage to arrange for the stabilisers to be constant dimension on $Y_{\min}$ (\textit{cf.}\ \cite[Remarks 2.19 and 8.10]{BDHK2}).
	\end{enumerate}

In this paper we will only consider the case where condition \ref{CZminNE} holds, i.e.\ that there exists a point in $Z_{\min}$ whose $U$-stabiliser dimension is minimal amongst all $U$-stabilisers of points of $Y_{\min}$.  This means we are at worst in case $2$(a), and do not need to blow up all of $Z_{\min}$. Later on in $\S$\ref{subsec no blow ups qnt in stages} we will introduce condition \hyperlink{WUU}{(WUU)}, which guarantees this in the setting of Quotienting-in-Stages. In the following subsection, we will explicitly describe the sequence of blow-ups in the case when $U$ is abelian, as we will only need this setting later in the paper.

	\begin{rmk}
		If the action does not satisfy \ref{starR0}, then one can perform reductive blow-ups so that on the blow-up either \ref{starR0} holds or semistability equals Mumford stability for this reductive group action (i.e. the dimensions of the $\overline{L}$-stabilisers of semistable points are constant).
	\end{rmk}

	\subsubsection{The blow-up procedure in a simple abelian setting}\label{sec easy blowups abelian}

Throughout this section, we suppose that the unipotent radical $U$ is abelian, and consider
\[ d_{\min}:= d_{\min}(Y_{\min}^{ss},U) \quad \text{and} \quad d_{\max}:= d_{\max}(Y_{\min}^{ss},U). \]
If \ref{starURss} fails, then the $U$-stabilisers have non-constant dimension on $Y_{\min}^{ss}$, so $d_{\min} < d_{\max}$. In this case, we perform a sequence of $H$-equivariant blow-ups to arrange for \ref{starURss} to hold on the blow-up. The idea is to start by considering points with maximal dimensional $U$-stabilisers. In fact, there are several possible ways to perform this blow-up sequence: one can blow-up the points in $Y_{\min}$ whose stabilisers have maximal $U$-dimension or blow-up the $U$-sweep of the locus of points in $Z_{\min}$ whose stabilisers have maximal $U$-dimension (the latter is equivalent to blowing-up the locus of points in $Y_{\min}$ with maximal dimensional $\hU$-stabiliser). In \cite{BDHK2}, the former approach is taken, but we will take the latter approach. This approach involves blowing up a smaller locus and it is easier to show this process terminates, as one is always blowing up a locus of codimension at least 2. 

Throughout this section, we will also assume that condition \ref{CZminNE} holds. Let us introduce the centre of each blow-up and construct the first blow-up in each sequence using the notation of Definition \ref{def centres of blowups}.
		
\begin{defn}
For a linearised $H$-action on $Y$ such that \ref{starURss} fails, we define
\[ D(Y):=C(Y_{\min}^{ss},\hU).\]
We let $\pi_{(1)} : \widetilde{Y}_{(1)} \ra Y$ denote the blow-up of $Y$ along the closure of $D(Y)$.
\end{defn}

\begin{rmk}\label{rmk on centre blowup ab uni}
By definition, $D(Y)$ is the closed subscheme of $Y^{ss}_{\min}$ on which the dimension of the stabiliser for $\hU$ are maximal; equivalently this is the $U$-sweep of $C(Z^{ss}_{\min},U)$. As $U < H$ is a normal subgroup, the scheme $D(Y)$ and its closure in $Y$ are $H$-invariant. Therefore, $\pi_{(1)}$ is a $H$-equivariant blow-up and the induced $H$-action on $\widetilde{Y}_{(1)}$ admits an ample linearisation given by pulling back the linearisation on $Y$ and perturbing by a small multiple of the exceptional divisor. The assumption \ref{CZminNE} ensures that $Z_{\min} \not\subseteq D(Y)$. Consequently, the minimal weight space $\widetilde{Z}_{(1),{\min}}$ in the blow-up $\widetilde{Y}_{(1)}$ is the strict transform of the minimal weight space $Z_{\min}$ in $Y$. Assuming as we always do that $UZ_{\min} \subsetneq Y_{\min}$ (as otherwise we proceed as in Remark \ref{rmk Uhat part of proj completion constr} \eqref{rmk when UZ equals Ymin}), the centre of this blow-up has codimension at least 2, so the blown up space is not isomorphic to the space we started with.
\end{rmk}

Let us show that this first blow-up reduces the dimension of maximal stabilisers; we will follow the argument \cite[Proposition 8.8]{BDHK2} even though their blow-up procedure differs slights from ours.

\begin{lemma}\label{lemma uni stab drops blowups}
Assume that \ref{starURss} fails, but \ref{CZminNE} holds. Then the dimensions of $U$-stabilisers in the $\lambda$-minimal weight spaces drops after performing the blow-up $\pi_{(1)} : \widetilde{Y}_{(1)} \ra Y$
\[ d_{\max}(\widetilde{Z}_{(1),{\min}}^{ss},U) < d_{\max}(Z^{ss}_{\min},U) = d_{\max} \]
where $\widetilde{Z}_{(1),{\min}}^{ss}$ denotes the semistable locus in the minimal weight space for the blow-up $\widetilde{Y}_{(1)}$.
\end{lemma}
\begin{proof}
This proof is an expanded version of the argument outlined in \cite[Proposition 8.8]{BDHK2}. 
By Lemma \ref{lem stab cant incr on blowup} below and the fact that $\widetilde{Y}_{(1)} \ncong Y$ (see Remark \ref{rmk on centre blowup ab uni}), it suffices to show that if $\tilde{z}\in \widetilde{Z}_{(1),{\min}}^{ss}$ lies over $z \in C(Z^{ss}_{\min},U)$, then $\tilde{z}$ is not fixed by $U':= \stab_{U}(z)$. We have
\[ Z_{\min}^{U'} \subset C(Z_{\min},U,d_{\max}) \subset Z_{\min}. \]
By embedding $Y$ in a projective space, we can assume without loss of generality that $Y = \PP^m$ and as $Z_{\min}^{U'}$ and $Z_{\min}$ are linear subspace we can assume we have taken coordinates so there are $0 \leq  n \leq m' \leq m$ such that $Z_{\min}^{U'}$ (resp. $Z_{\min}$) is the vanishing locus of the first $m'$ (resp. n) coordinates and $p([x_0: \cdots : x_m]) = [0 : \cdots :0 : x_{n} : \cdots : x_m]$. We claim the $U'$-action on $\PP^m$ with respect to the blocks $n \leq m' \leq m$ has the form
\[ u' \mapsto \left( \begin{array}{ccc} A(u') & B(u') & 0 \\ 0 & I & 0 \\ 0 & 0 & I \end{array} \right). \]
Indeed, the final column has this form, as $U'$ acts trivially on $Z_{\min}^{U'}$. With respect to the two blocks given by $n \leq m$ (given by the minimal $\lambda$-weight), this matrix is block upper triangular, as $\lambda$ grades $U$. Finally, the middle column has this form, as $p(u \cdot x) = x$ for $x \in Z_{\min}$. 

Since the exceptional divisor is given by the projectivised normal bundle, $\tilde{z}$ is represented by some $[\xi] \in \PP(T_zZ_{\min} / T_z(Z_{\min}^{U'}))=\pi_{(1)}^{-1}(z)$. We can take coordinates so $z = [0: \cdots :0 :1]$ and then take local coordinate $x_j/x_m$ around $z$. From the above block form of the action, one sees for any $\xi \in T_zZ_{\min} \setminus T_z(Z_{\min}^{U'})$, 

there is a $u' \in U'$ such that $u' \xi - \xi \notin T_zZ_{\min}$. Indeed, in the above representation of the $U'$-action, we have some $u'$ where $B(u') \neq 0$, as $Z_{\min}^{U'}  \neq Z_{\min}$. Therefore $\tilde{z}$, which is represented by some $[\xi]$, is not fixed by $U'$. 
\end{proof}

We used the following straight-forward lemma in the above proof.

\begin{lemma}\label{lem stab cant incr on blowup}
Let $\pi: \widetilde{X} \ra X$ be a $G$-equivariant blow-up for any linear algebraic group $G$. Then for any subgroup $G' < G$ and any $x \in \widetilde{X}$, we have
\[ \dim \stab_{G'}(x) \leq \dim \stab_{G'}(\pi(x))\]
with equality holding for all points outside the exceptional divisor.
\end{lemma}

We can finally prove that \ref{CZminNE} enables a blow-up procedure that does not blow up all of $Z_{\min}$.

\begin{prop}\label{prop ab uni blowups seq} 
Let $H = U \rtimes L$ be a group with abelian unipotent radical internally graded by $\lambda: \GG_m \ra L$. Let $Y$ be a projective $H$-scheme with a very ample adapted linearisation $\cL$ such that \ref{starR0} and \ref{CZminNE} hold. Then there is a sequence of equivariant blow-ups along $H$-invariant closed subschemes resulting in a projective scheme $\widetilde{Y}$ admitting a very ample adapted linearisation\footnote{This linearisation is on a line bundle obtained by tensoring the pullback of $\cL$ with small multiples of the exceptional divisors for each blow-up.} such that on $\widetilde{Y}$ both \ref{starR0} and \ref{starURss} hold. Moreover, \ref{CZminNE} ensures that not all of $Z_{\min}$ is blown up in this procedure.
\end{prop}
\begin{proof} 
If \ref{starURss} holds already for $Y$, then no blow-ups are needed. Otherwise, let $\widetilde{Y}_{(0)} = Y$ and for $i =0,\dots$, define $\widetilde{Y}_{(i+1)}$ to be the blow-up of $\widetilde{Y}_{(i)}$ along the closure of $D(\widetilde{Y}_{(i)})$. Since the dimensions of the $U$-stabilisers on the semistable locus in the $\lambda$-minimal weight space decreases at each stage, this procedure terminates with a scheme $\widetilde{Y}:=\widetilde{Y}_{(n)}$ which has constant dimensional $U$-stabilisers on the semistable locus in the $\lambda$-minimal weight space equal to $d_{\min}$. We claim that $\dim \stab_{U}$ is constant and equal to $d_{\min}$ on all of $\widetilde{Y}_{\min}^{ss}$. If $\tilde{y} \in \widetilde{Y}_{\min}^{ss}$ flows to $\tilde{z} = \tilde{p}(\tilde{y}) \in \widetilde{Z}^{ss}_{\min}$ under $\lambda$, then
\[ d_{\min} \leq \dim \stab_{U}(\tilde{y}) \leq \dim \stab_{U}(\tilde{z}) = d_{\min}.\]
Thus \ref{starURss} holds for $\widetilde{Y}$, as does \ref{starR0} by Lemma \ref{lem stab cant incr on blowup}, as it already held for $Y$.
\end{proof}

In this case, we define the following stable locus.
	
\begin{defn}\label{def hat stable locus easy abelian case}
Assume that \ref{starR0} and \ref{CZminNE} holds, where $U$ is abelian. Then the \emph{stable locus}\footnote{If additionally \ref{starURss} holds, then this coincides with the stable locus in Definition \ref{def nonred ss}.} for the linearised $H$-action on $Y$ is
\[ Y^{H-s}:=\{ y \in Y^s_{\min} \setminus U Z^s_{\min}: \: \dim \stab_U(p(y)) = d_{\min} \}. \]
\end{defn}

The following result is a simplified version of \cite[Theorems 2.10 and 2.20]{BDHK2} which shows this stable set admits a quasi-projective geometric $H$-quotient (see also \cite[$\S$3.3]{Josh_length2}).
	
\begin{prop}\label{prop easy blowups abelian} 	
Let $H = U \rtimes L$ be a group with abelian unipotent radical internally graded by $\lambda: \GG_m \ra L$ acting on a projective scheme $Y$ with respect to a very ample adapted linearisation $\cL$. Assume that \ref{starR0} and \ref{CZminNE} hold. If $\pi : \widetilde{Y} \ra Y$ is the $H$-equivariant blow-up of Proposition \ref{prop ab uni blowups seq}, then $\pi$ is an isomorphism over the open set $ Y^{H-s}$ and thus via $\pi$, the set $ Y^{H-s}$ admits a quasi-projective geometric $H$-quotient with natural projective completion $\widetilde{Y}/\!/H$.
\end{prop}
\begin{proof} 
We note that $Y^{H-s}$ is disjoint from the centre of each blow-up appearing in Proposition \ref{prop ab uni blowups seq} and so is contained in the open set over which $\pi$ is an isomorphism. In fact, we claim $Y^{H-s}$ is contained in the intersection of the open set over which $\pi$ is an isomorphism with the stable set $ \widetilde{Y}^{H-s}$ in the blow-up. For simplicity, suppose $H = \hU$. Since not of all of $Z_{\min}$ is blown up, the minimal weight space in the blow-up is the strict transform of $Z_{\min}$ and therefore $Y^{\hU-s}$ is contained in the intersection of the open set over which $\pi$ is an isomorphism with $ \widetilde{Y}^{\hU-s} = \widetilde{Y}_{\min} \setminus U \widetilde{Z}_{\min}$. 
\end{proof}

This stable set could potentially be empty; for example, when $H = \hU$, this can happen if we have $UZ_{\min} \subsetneq X_{\min}$, but $UC(Z_{\min},U,d_{\min}) = p^{-1}(C(Z_{\min},U,d_{\min}))$ as in the next example.

\begin{ex}
Let $H = \widehat{\GG_a} = \GG_a  \rtimes \GG_m$ act on $Y:=\PP^2$ with coordinates $[x:y:z]$ via the representation $\GG_a  \rtimes \GG_m \ra \GL_3$
\[ (a,t) \mapsto \left(\begin{array}{ccc} t & a & 0 \\ 0 & t^{-1} & 0 \\ 0 & 0 & t^{-1}\end{array} \right).\]
Then $Z_{\min} = \{ x = 0 \} \cong \PP^1$ and $Y_{\min} \ra Z_{\min}$ is an $\AA^1$-bundle. Furthermore, the only point in $Z_{\min}$ with non-trivial $\GG_a$-stabiliser is $y = [0: 0:1]$ which is fixed by $\GG_a$; thus $d_{\min} = 0$ and $C:=C(Z_{\min},U,0)= Z_{\min} \setminus \{y \}$. Note that $UZ_{\min} \subsetneq X_{\min}$, but $UC = p^{-1}(C)$ and so we have $Y^{H-s} = \emptyset$ using the above definition. Since \ref{starU} fails, we need to do blow-ups, but as \ref{CZminNE} holds (and \ref{starR0} holds trivially) we can apply the above result. In this case, $\widetilde{Y}$ is just the blow-up of $Y$ at $y$ and $\widetilde{Z}_{\min} = \PP^1$ (as it is the blow-up of $Z_{\min}$ at $y$) and one now has $U\widetilde{Z}_{\min} = \widetilde{Y}_{\min}$, so also $\widetilde{Y}^{H-s} = \emptyset$. In this case, one should instead define the stable locus in $\widetilde{Y}$ to be $\widetilde{Y}_{\min}$ with $U$-quotient $\widetilde{Y}_{\min} \ra Z_{\min}$ as in Remark \ref{rmk Uhat part of proj completion constr} \eqref{rmk when UZ equals Ymin}, which under $\pi$ would give a geometric quotient of the open set $p^{-1}(C)=UC$ in $Y$ via the map $p^{-1}(C) \ra C$.
\end{ex}

	\section{Quotients of unstable HKKN strata}\label{sec quot HKKN strata}
	
	Let $G$ be a reductive group acting on a projective scheme $Y$ with respect to an ample linearisation $\cL$. As in $\S$\ref{sec HKKN strat}, associated to a choice of norm $||-||$ on $\Lie G$ there is a HKKN stratification 
	\[ Y = \bigsqcup_{\beta \in \cB} S_\beta \]
	into finitely many $G$-invariant locally closed subschemes $S_\beta \subset Y$. If $Y^{ss} \neq \emptyset$, then this semistable locus is the minimal stratum indexed by $0 \in \cB$ and admits a good quotient: the reductive GIT quotient. If the stable locus is non-empty, this GIT quotient restricts to a geometric quotient on the stable locus. The unstable HKKN strata $S_\beta$, which correspond to non-zero indices $\beta \in \cB$, admit a description
	\[ S_\beta \cong G \times^{P_\beta} Y_\beta^{ss}\]
	for a parabolic subgroup $P_\beta \subset G$ and a locally closed subscheme $Y_\beta^{ss} \subset Y$. Therefore, constructing a $G$-quotient of an open subscheme of $S_\beta$ is equivalent to constructing a $P_\beta$-quotient of an open subset of $Y_\beta^{ss}$. The goal of this section is to explain how to construct such a quotient.
	
\subsection{Categorical quotients from borderline linearisations}
	
The problem of constructing a quotient of the $G$-action on an unstable stratum $S_\beta \subset Y$, or equivalently the $P_\beta$-action on $Y_\beta^{ss}$, was studied in \cite[$\S$2]{HK}. As explained in $\S$\ref{sec intro mod unstable}, it is often preferable to consider the $P_\beta$-action on $Y_\beta^{ss}$ rather than the $G$-action on $S_\beta$, as the parabolic group admits more characters which can be used to twist the linearisation. The characters of $P_\beta$ correspond to the characters of the Levi $L_\beta <P_\beta$, and amongst these there is (after the choices we have made) a distinguished character, $\chi_{-\beta}$ from which we obtain the canonical $P_\beta$-linearisation $\cL_\beta$ by twisting the original linearisation $\cL$ by the character $\chi_{-\beta} : L_\beta \ra \GG_m$. This gives an ample $P_\beta$-linearisation on the closure of $Y_\beta^{ss}$, but does not extend to an ample $G$-linearisation on the closure of $S_\beta$, since the character $\chi_{-\beta}$ may not extend to $G$. 

As we will see in Lemma \ref{lemma parabolic graded}, the unipotent radical of $P_\beta$ is always graded by $\lambda_\beta$, and with respect to this grading the canonical linearisation is borderline, in the sense of Definition \ref{def borderline adapted}. Thus, this canonical linearisation provides a categorical quotient of the unstable stratum, via the following, which we generalised in Proposition \ref{prop canonical linearisation}.
	
	\begin{prop}\cite[Proposition 3.2]{HK} \label{prop cat quot same}
		The projective variety $Z_\beta /\!/_{\!\cL_\beta} L_\beta$ is a categorical quotient for the $L_\beta$-action on $Z_\beta^{ss}$, the $P_\beta$-action on $Y_\beta^{ss}$, and the $G$-action on $S_\beta$.
	\end{prop}
	
	In general this categorical quotient of the $G$-action on $S_\beta$ and of the $P_\beta$-action on $Y_\beta^{ss}$ is far from being an orbit space, as it identifies each $y \in Y_\beta^{ss}$ with $p_\beta(y):= \lim_{t \ra 0} \lambda_\beta(t) \cdot y \in Z_\beta^{ss}$, which in general only lies in the boundary of the orbit of $y$.

\begin{ex}
For an unstable HN stratum $S_\tau$ for sheaves on a polarised projective scheme, the map $p_\tau : Y_\tau^{ss} \ra Z_\tau^{ss}$ corresponds to sending a sheaf to its associated graded for the HN filtration. Hence the above categorical quotient is just the product of moduli spaces of semistable sheaves for each Hilbert polynomial $P_i$ in the HN type $\tau$.
\end{ex}	

\begin{ex} 
Consider configurations of ordered tuples of $n\geq 4$ points on $\PP^1$ up to $\SL_2$-action, with the natural linearisation and the Killing form on $\SL_2$. Such a configuration is semistable as long as no more than $n/2$ of the points coincide anywhere, and \cite{Kirwan_thesis} shows that the HKKN stratification measures how many points do coincide. In this case the above categorical quotient given above is most unsatisfactory: for any $\beta$ it is simply a point; see \cite{Josh_length2} for more details.
\end{ex}

To avoid this type of identification, one would need to first remove $P_\beta Z_\beta^{ss} = U_\beta Z_\beta^{ss}$, ideally by perturbing the $P_\beta$-linearisation $\cL_\beta$ on $Y_\beta^{ss}$ (or strictly speaking its closure) by a sufficiently small rational character of $P_\beta$, so that all points in $Z_\beta^{ss}$ are unstable for this perturbed linearisation. In fact, we will construct a  quotient of the action of the non-reductive group $P_\beta$ on an open subset of semistable points of $Y_\beta^{ss}$, which is disjoint from $Z_\beta^{ss}$, by using a perturbation of the canonical linearisation $\cL_\beta$ and recent results from non-reductive GIT.

\subsection{Applying non-reductive GIT to parabolic actions on unstable strata}
	
In this section, we continue to assume we have a HKKN stratum $S_\beta \cong G \times^{P_\beta} Y_\beta^{ss}$. In order to apply the above results from non-reductive GIT to the $P_\beta$-action on $Y_\beta^{ss}$, we first need to compactify $Y_\beta^{ss}$. Let $X_\beta:=\overline{Y_\beta}$ denote the closure of the locally closed subscheme $Y_\beta \subset Y$. We recall that we gave explicit descriptions \eqref{weight desc Ybeta} of $Y_\beta$ and $Z_\beta$ in $\S$\ref{sec HKKN strat} in terms of torus weights.

	First, we check that the non-reductive group $P_\beta$ has graded unipotent radical.
	
	\begin{lemma}\label{lemma parabolic graded}
		The parabolic subgroup $P_\beta = U_\beta \rtimes L_\beta$ has graded unipotent radical $U_\beta$, where the grading $\GG_m$ is given by the central 1-PS $\lambda_\beta : \GG_m \ra L_\beta$.
	\end{lemma}
	\begin{proof}
		By definition $P_\beta = P(\lambda_\beta)$ as in \eqref{eq parabolic}, it follows that the weights of the conjugation action of $\lambda_\beta(\GG_m)$ on $\Lie P_\beta$ are non-negative and are precisely zero on the Levi factor $L_\beta$. Hence, all the weights of $\lambda_\beta(\GG_m)$ on $\Lie U_\beta$ are strictly positive.
	\end{proof}
	
The associated minimal weight space and attracting open set admit the following descriptions.
	
	\begin{prop}\label{prop nonred GIT for HKKN strata}
		Let $\beta \in \cB$ be a non-zero index of a HKKN stratum $S_\beta \subset Y$. For the $P_\beta$-action on $X_\beta =\overline{Y_\beta}$ with respect to the linearisation $\cL$ and the grading given by $\lambda_\beta : \GG_m \ra L_\beta$ we have
		\[ Z_{\beta,\min}:=Z(X_\beta,\lambda_\beta)_{\min} = Z_\beta \quad \text{and} \quad X_{\beta,\min}:=X_\beta(\lambda_\beta)_{\min} = Y_{\beta} \]
		and the associated retraction is $p_\beta : Y_\beta \ra Z_\beta$.
		Moreover, for $\epsilon \in \QQ_{>0}$ sufficiently small, the twisted linearisation $\cL^{\chi_{-(1+\epsilon)\beta}}$ is adapted for the graded unipotent group $\hU_\beta= U_\beta \rtimes_{\lambda_\beta} \GG_m$.
	\end{prop}
	\begin{proof} 
		Since for a $G$-equivariant closed immersion $i: Y \hookrightarrow \PP^n$ such that $\cL = i^* \cO_{\PP^n}(1)$, the HKKN strata on $Y$ are the fibred products of the HKKN strata on $\PP^n$ with $Y$ and similarly for the fixed loci for $\lambda_\beta$, it suffices to prove that $Z_{\beta,\min} = Z_\beta$ when $Y = \PP^n$ and $\cL= \cO_{\PP^n}(1)$. We can choose a maximal torus $T$ containing $\lambda_\beta(\GG_m)$ and coordinates on $\PP^n$ such that the action of the fixed maximal torus $T < G$ is given by $t \mapsto \diag(\alpha_0(t), \dots , \alpha_n(t))$ for characters $\alpha_i : T \ra \GG_m$. 
		
		Recall that we use an invariant inner product on $\mathfrak{t}$ to identify characters and co-characters of $T$ and the associated norm is denoted $|| - ||$. The weights of the $\lambda_\beta(\GG_m)$-action on $H^0(\PP^n,\cO_{\PP^n}(1))^*$ are given by the projections $\frac{\alpha_i \cdot \beta}{||\beta||^2} \beta$ of the $T$-weights $\alpha_i \in \mathfrak{t}$ to the line spanned by $\beta \in \mathfrak{t}$. Thus to determine the minimum $\lambda_\beta(\GG_m)$-weight $\omega_{\min}$, we want to minimise $\alpha_i \cdot \beta$ over all weights $\alpha_i$ of points in $Y_\beta$. As $\beta \in \cB$, we have $Z_\beta \neq \emptyset$ and so from the descriptions given in \eqref{weight desc Ybeta} of $\overline{Y_\beta}$ and $Z_\beta$, we see that there is a point in $Z_\beta$ with at least one $T$-weight $\alpha_i$ satisfying $\alpha_i \cdot \beta = || \beta ||^2$. The weights $\alpha_j$ of points in $\overline{Y_\beta}$ all satisfy $\alpha_j \cdot \beta \geq || \beta ||^2$. Hence, the minimal  $\lambda_\beta(\GG_m)$-weight corresponds to $\beta \in \mathfrak{t}$ and by definition $Z_{\beta,\min} \subset X_\beta^{\lambda(\GG_m)}$ is the minimal $\lambda_\beta(\GG_m)$-weight space (\textit{cf.}\ Definition \ref{def Zmin}); that is,
		\[ Z_{\beta,\min} = \{ [p_0 : \dots : p_n] : p_i = 0 \text{ if } \alpha_i \cdot \beta \neq ||\beta||^2 \},\]
		which is equal to $Z_\beta$ by  \eqref{weight desc Ybeta}. From the description in  \eqref{weight desc Ybeta} of $Y_\beta$, we see that $Y_\beta$ is the open subscheme of $X_\beta =\overline{Y_\beta}$ consisting of points $y$ such that $\lim_{t \ra 0} \lambda_\beta(t) \cdot y \in Z_\beta$, and so it follows that $X_{\beta,\min} = Y_\beta$ (\textit{cf.}\ Definition \ref{def Zmin}). Note that if we twist $\cL$ by any scalar multiple of the character corresponding to $\beta$, then the subsets $X_{\beta,\min}$  and $Z_{\beta,\min}$ remain unchanged.
		
		For the final statement, recall that the (rational) character $\chi_{-(1+\epsilon)\beta} :P_\beta \ra \GG_m$ merely shifts the $\lambda_\beta(\GG_m)$ weights by $-(1 + \epsilon)||\beta||$. In order for the $\hU_\beta$-action with respect to this twisted linearisation to be adapted, we need the origin to separate the minimal weight from all the other weights. Since the minimum weight for the original linearisation $\cL$ corresponds to $\beta$, we need to shift this weight just beyond the origin so that this weight is negative and all the other weights are positive. This is achieved by taking $\epsilon \in \QQ_{>0}$ sufficiently small.
	\end{proof}
	
	We fix $\epsilon \in \QQ$ sufficiently small as above, so that the perturbation $\cL_\beta^{\text{per}}:=\cL^{\chi_{-(1+\epsilon)\beta}} $ of the canonical linearisation $\cL_\beta$  on $\overline{Y_\beta}$ is adapted for the $\hU_\beta$-action. Now we can prove Theorem \ref{thm nred quot unstable strata}.

	\begin{proof}[Proof of Theorem \ref{thm nred quot unstable strata}]
		The proof follows from Theorem \ref{mainthm nrGIT quotient}. For the second statement, we just need to determine the semistable locus for the action of the reductive group $\overline{L_\beta} := L_\beta /\lambda_\beta(\GG_m)$ on $Z_\beta$ with respect to the perturbed linearisation $\cL_\beta^{\text{per}}$, as by Definition \ref{def nonred ss} we have
		\[ X_\beta^{P_\beta-s}= \{ y \in X_\beta^{\hU_\beta-s} : p_\beta(y) \in Z_\beta^{\overline{L_\beta}-s}(\cL_\beta^{\text{per}})\}. \]
Since \ref{starR0} holds, semistability coincides with stability for the $\overline{L_\beta}$-action on $Z_\beta$ linearised by $\cL_\beta$. By variation of (reductive) GIT quotients, this linearisation $\cL_\beta$ lies inside a VGIT chamber and so the stable locus is unaltered on passing to the small perturbation $\cL_\beta^{\text{per}}$; that is,
\begin{equation}\label{eq ss for per linear}
		Z_\beta^{\overline{L_\beta}-s}(\cL_\beta^{\text{per}}) = Z_\beta^{\overline{L_\beta}-s}(\cL_\beta)= Z_\beta^{\overline{L_\beta}-ss}(\cL_\beta)=Z_\beta^{{L_\beta}-ss}(\cL_\beta)=:Z_\beta^{ss}
		\end{equation}
which completes the proof.		\end{proof}
	
	\begin{rmk}\
		\begin{enumerate}
			\item As $Z_\beta$ is unstable for the $\lambda_\beta(\GG_m)$-action on $\overline{Y_\beta}$ linearised by $\cL^{\text{per}}_\beta$, we see that a semistable point $y \in Y_\beta$ cannot be identified with its limiting flow $p_\beta(y) \in Z_\beta$; hence the quotients constructed by the above theorem avoid the collapsing described in Proposition \ref{prop cat quot same}. Indeed they are geometric quotients and so are orbit spaces.
			\item The above theorem also holds on replacing condition \ref{starU0} with condition \ref{starU} (or \ref{starURss}). If this condition fails, then one can construct a quotient of the $P_\beta$-action on an open subset of $\overline{Y}_\beta$ by performing a sequence of blow-ups as described in $\S$\ref{sec nrGIT blowups}. 
		\end{enumerate}
	\end{rmk}
	
\subsubsection{Application to moduli of sheaves of fixed HN type}\label{sec sheaves nonred git setup}

Let $\tau = (P_1, \dots, P_l)$ be a length $l \geq 2$ HN type for sheaves on a polarised projective scheme $(B,\cO_B(1))$. As explained in $\S$\ref{sec instab sheaves}, there is a close relationship between the HN stratification on the Quot schemes used for the construction of moduli of sheaves and the associated HKKN stratification. In particular,  constructing a moduli space of sheaves of HN type $\tau$ is equivalent to constructing a quotient of the action of a parabolic group $P_\tau$ on (an open subset of) $Y_\tau^{ss}$. We would ideally like to apply Theorem \ref{thm nred quot unstable strata} above to the $P_\tau$-action on $X_\tau:=\overline{Y_\tau}$ linearised by $\cL_\tau^{\text{per}}$. However, for sheaves of fixed HN type the various stabiliser conditions in Definition \ref{def star cond} do not hold in general. Moreover, for length $l > 2$, we have the following issues concerning both the reductive and unipotent stabilisers (\emph{cf.} \cite{Josh_length2}).
\begin{enumerate}[label={(\roman*)}]
\item The central torus $T_\tau :=Z(L_\tau) \cong \GG_m^{l-1}$, which scales the pieces in the HN graded sheaf, acts trivially on $Z_\tau$ but non-trivially on $X_\tau$. Thus the $T_\tau$-weight picture for $X_\tau$ linearised by $\cL_\tau$ looks like a cone with a vertex at the origin. Consequently semistability equals stability for the action of $\overline{L_\tau} = L_\tau/\lambda(\GG_m)$ on $Z_\tau$ if and only if $\tau$ is a coprime HN type (i.e. semistability coincides with stability for all HN subquotients) of length $2$. Moreover, for $l >2$, the $\overline{L_\tau}$-stable locus is empty, which implies that the $P$-stable locus in Definition \ref{def nonred ss} is also empty. 
\item In length $l >2$, the unipotent radical $U_\tau$ is no longer abelian. Moreover, although the unipotent stabiliser groups often admit sheaf theoretic interpretations, the dimension of these unipotent stabilisers can vary. In particular, for $l > 2$, the map $p_\tau : Y_\tau^{ss} \ra Z_\tau^{ss}$ does not preserve the dimensions of unipotent stabiliser groups in general, thus the blow-up procedure in \cite{BDHK2} is almost intractable as we may blow-up all of $Z_{\tau}^{ss}$.
\end{enumerate}
These statements are proved in Propositions \ref{reductive sheaf stabs} and \ref{prop uni stab sheaves} below.

These issues are not specific to sheaves of fixed HN type. For moduli of objects in an abelian linear category (for example, moduli of quiver representations, of Higgs sheaves or other decorated sheaves) one also encounters these problems for HN types of length $l > 2$. This is the main motivation for the Quotienting-in-Stages procedure we develop later in this paper.

Let us explain the issues mentioned above in the case of sheaves of HN type $\tau$.

\begin{prop} \label{reductive sheaf stabs}
Consider the linearised action of $P_\tau = U_\tau \rtimes L_\tau$ on $X_\tau = \overline{Y_\tau}$.
\begin{enumerate}[label=\emph{(\roman*)}]
\item The central torus $T_\tau:=Z(L_\tau) \cong \GG_m^{l-1}$ acts trivially\footnote{However, $T_\tau$ acts non-trivially on $Y_\tau$ if there exists a sheaf $\cE$ such that $\cE \ncong \gr(\cE)$.} on $Z_\tau$.
\item If $l > 2$, then $Z_\tau^{\overline{L_\tau}-s} = \emptyset$ and $Z_\tau^{\overline{L_\tau}-ss} = Z_\tau^{ss}$. Thus semistability does not coincide with stability for the $\overline{L_\tau}$-action on $Z_\tau$.
\item If $l =2$, then $Z_\tau^{\overline{L_\tau}-s} = Z_\tau^{\overline{L_\tau}-ss}$ if and only if $\tau$ is coprime.
\end{enumerate}
\end{prop}
\begin{proof}
Equation \eqref{Ztau} expresses $Z_\tau$ as a product of closed subschemes $R_n(P_i)$ of the Quot schemes $Q_n(P_i)$. As $L_\tau \cong (\prod_{i=1}^l \GL_{P_i(n)} ) \cap \SL_{P(n)}$ and the diagonal $\GG_m$ in each $\GL_{P_i(n)}$ acts trivially on each Quot scheme $Q_n(P_i)$, we see that $T_\tau$ acts trivially on $Z_\tau$. Since $\overline{T_\tau}:=T_\tau/\lambda_\tau(\GG_m) \cong \GG_m^{l-2} < \overline{L}_\tau$ acts trivially on $Z_\tau$, we see for $l > 2$ that  GIT semistability and stability do not coincide for $\overline{L}_\tau$. For $l =2$, we have a surjection $R_{\tau}:=\SL_{P_1(n)} \times \SL_{P_2(n)} \twoheadrightarrow \overline{L_\tau}$ with finite kernel and so (semi)stability with respect to these groups coincides. The final statement then follows, as GIT semistability coincides with stability for the $\SL_{P_i(n)}$-action on $Q^{ss}(P_i)$ if and only if for sheaves with Hilbert polynomial $P_i$ semistability coincides with stability.
\end{proof}

Let us now consider the various stabiliser conditions for the unipotent group $U_\tau$. If $l > 2$, then $U_\tau$ is not abelian and so we need to consider stabilisers for the subgroups appearing in the lower central series $U^{\{l-1\}}_\tau= \{ I \} < U_{\tau}^{\{l-2\}} < \cdots  < U_{\tau}^{\{0\}}= U_\tau$ (see Remark \ref{rmk concrete desc parabolic and lcs} below). The $U_{\tau}^{\{l-2\}}$-stabilisers of quotient sheaves $q_{\cE} : \cO_B(-n)^{\oplus P(n)} \twoheadrightarrow \cE$ are related to certain filtered endomorphism groups for the HN filtration $0 = \cE^{(0)} \subsetneq \cE^{(1)} \subsetneq \cdots \subsetneq \cE^{(l)} = \cE$. Let $\End_{-j}(\cE)$ denote the global sections of the sheaf $\cH om_{-j}(\cE,\cE)$, whose sections over $U \subset B$ are 
\[ \cH om_{-j}(\cE,\cE)(U):=\{ \varphi \in \cH om(\cE,\cE)(U)  : \varphi(\cE^{(i)}(U)) \subset \cE^{(i-j)}(U) \text{ for all }\: i\}.\]
	
\begin{prop}\label{prop uni stab sheaves}
For $q_{\cE} : \cO_B(-n)^{\oplus P(n)} \twoheadrightarrow \cE$ in $Y_\tau^{ss}$, the following statements hold.
		\begin{enumerate}[label=\emph{(\roman*)}]
\item \label{prop stab 1b} We have $\Lie \stab_{U_\tau}(q_{\cE}) \cong \End_{-1}(\cE)$ and, for $0 \leq i \leq l-1$, we have
			\[ \Lie \stab_{U_\tau^{\{i\}}}(q_{\cE}) \cong \End_{-(i+1)}(\cE).\]
\item \label{prop stab 2} If $q_{\cE} \in Z_{\tau}^{ss}$, then $\Lie \stab_{U_\tau^{\{i\}}}(q_{\cE}) = \oplus_{j>k+i} \Hom(\cE_j,\cE_k)$ for $0 \leq i \leq l-1$.

		\end{enumerate}
In particular, the $P_\tau$-action on $\overline{Y}_\tau$ satisfies condition \ref{starURss} with respect to the lower central series if and only if for $ 1 \leq j < l$, the function $\dim \End_{-j}(-)$ is constant on the set of all sheaves over $X$ with HN type $\tau$. 
\end{prop}
\begin{proof}
Statement (i) for $i = 0$ follows as the stabiliser in $U_\tau$ of $q_\cE$ consists of all automorphisms of $\cE$ of the form $\text{Id}_\cE +\alpha$, where $\alpha \in \End_{-1}(\cE)$ for the HN filtration and the version for higher $i$ follow similarly. Statement (ii) follows as for $q_{\cE} \in Z_{\tau}^{ss}$, we have $\cE = \gr(\cE)$. 

\end{proof}

Thus, as explained in \cite{Josh_length2}, even in HN length $l = 2$ conditions \ref{starU0}, \ref{starU0Rss}, \ref{starU} and \ref{starURss} fail in general. However, in length $l=2$ the above Proposition \ref{prop uni stab sheaves} simplifies as follows:
\[ \Lie \stab_{U_\tau}(q_{\cE}) \cong \Hom(\cE_2,\cE_1) \cong \Lie \stab_{U_\tau}(p_\tau(q_{\cE})) \]  
where $p_\tau(q_{\cE}) \in Z_\tau^{ss}$ is the quotient sheaf for the associated graded $\gr(\cE) = \cE_1 \oplus \cE_2$. This means that $p_\tau$ preserves the dimension of $U_\tau$ stabilisers, which is enough to ensure that we do not blow up the whole minimal weight space when performing the blow-up process to obtain a quotient. In this way, \cite{Josh_length2} constructs moduli spaces of sheaves of fixed HN type for length $l=2$.

However, in general this property of $p_\tau$ fails along with  \ref{starU0}, \ref{starU0Rss}, \ref{starU} and \ref{starURss} for $l>2$. Let us describe why, for $l=3$. In this case $\End_{-2}(\cE) = \Hom(\cE_3,\cE_1)$ and $\End_{-1}(\cE)$ is given by a long exact sequence
\[ 0 \ra \End_{-2}(\cE) \ra \End_{-1}(\cE) \ra  \Hom(\cE_{3},\cE_{2}) \oplus \Hom(\cE_2,\cE_1) \stackrel{\delta(\cE)}{\ra} \Ext^1_{-2}(\cE,\cE) \ra \cdots \] 
where $\delta(\cE) \neq 0$ in general. For the associated graded, we have $\delta(\gr(\cE)) = 0$ and 
\[  \End_{-1}(\gr(\cE)) = \Hom(\cE_3,\cE_1) \oplus \Hom(\cE_3,\cE_2) \oplus \Hom(\cE_2,\cE_1). \]
Hence, for length $l = 3$, we have
\[ \dim \Lie \stab_{U_\tau}(q_{\cE}) = \dim \Lie \stab_{U_\tau}(p_\tau(q_{\cE})) \iff \delta(\cE) = 0\]
and so $p_\tau$ will not in general preserve $U_\tau$-stabiliser dimension, as $\delta(\cE)$ will be non-zero in general. As the length increases, we need to consider $l-2$ long exact sequences to relate the dimension of $\End_{-1}(\cE)$ with $\End_{-1}(\gr(\cE)) = \oplus_{i > j} \Hom(\cE_i,\cE_j)$. If the coboundary map in this long exact sequences is zero, then these groups have the same dimension.

\section{Quotienting-in-Stages for parabolics: Definitions and Results}\label{sec quotient stages par} \label{subsec no blow ups qnt in stages}

In this section, we introduce a new method for constructing quotients by parabolic group actions in stages, where we use a sequence of different 1-PSs which are central in the Levi subgroup to grade the unipotent radical of a maximal parabolic at each stage. Our primary motivation is to construct quotients of unstable HKKN strata $S_\beta \cong G \times^{P_\beta} Y_\beta^{ss}$ representing moduli of objects of a fixed HN type in a linear abelian category. Our construction will be suited to the case when the centre of $L_\beta$ acts trivially on $Z_\beta$ and so one cannot directly apply Theorem \ref{thm nred quot unstable strata}. We will focus on the case of a parabolic subgroup $P$ of the special linear group $\SL_N$.

\subsection{Parabolic subgroups of the special linear group} \label{subsec parabolic subgps of SL}

For concreteness, we will assume that our parabolic subgroup $P < \SL_N$ is block upper triangular, since we can arrange this by conjugating by an element in $\SL_N$. In particular, we assume that $P$ is a parabolic subgroup associated to a diagonal (rational) 1-PS with decreasing weights as follows.

\begin{defn}
For a (rational) 1-PS $\lambda : \GG_m \ra \SL_N$ of the form 
	\begin{equation}\label{eq rep conj class}
	\lambda(t) = \diag (t^{r_1}, \dots , t^{r_N}) \:\text{ with } \:  r_1 \geq r_2 \geq \cdots \geq r_N \: \text{ such that } \sum_{i=1}^N r_i = 0,
	\end{equation}
we let $P = P(\lambda)$ be the associated parabolic subgroup and we write $P = U \rtimes L$ where $U$ is the unipotent radical and $L$ is the Levi factor. Let $l=l(\lambda)\geq 2$ denote the number of distinct weights of $\lambda$, and write these weights as $\beta_1 > \dots > \beta_l$ with multiplicities $\underline{m} = \underline{m}(\lambda) =(m_1,\dots, m_l)$. We refer to $l$ and $\underline{m}$ as the \emph{length} and \emph{multiplicities} of $\lambda$ (or of $P$). If $l =2$, then $P$ is a \emph{maximal parabolic subgroup} and $U$ is abelian, and if $l> 2$, then $U$ is non-abelian. 
\end{defn}
	
Note that the multiplicities sum to $N$ and $\sum_{i=1}^l m_i \beta_i =0$ as $\lambda$ is a 1-PS of $\SL_N$. Moreover, $l(\lambda)$ and $\underline{m}(\lambda)$ only depend on the conjugacy class of $\lambda$. The parabolic group $P$ is explicitly determined from $l$ and $\underline{m}$ as follows.
	
\begin{rmk}\label{rmk concrete desc parabolic and lcs}
For $\lambda$ as in \eqref{eq rep conj class}  with length $l = l(\lambda)$ and multiplicities $(m_1,\dots, m_l)$, we have
\[ P= P(\lambda) = \left\{\left(\begin{array}{ccccc} A_{11} & A_{12} & \cdots & \cdots & A_{1l}\\ 0 & A_{22} & A_{23} & \cdots & A_{2l} \\ \vdots & \ddots & \ddots &  & \vdots \\  \vdots & & 0 & \ddots & A_{l-1 \: l} \\ 0 & \cdots & \cdots & 0 & A_{ll} \end{array} \right) \in \SL_{N} : A_{ij} \in \Mat_{m_i \times m_j} \right\}.\]
The Levi subgroup $L = \{ A \in P_\beta : A_{ij} = 0\text{ for } 1 \leq i<j\leq l \} < P$ consists of block diagonal matrices and is the intersection of the block diagonal embedding $\prod_{i=1}^l \GL_{m_i} \hookrightarrow \GL_N$ with $\SL_N$. The centre of the Levi is a rank $l-1$ torus, denoted
\[ T: =Z(L) =  \{ \diag(t_1 I_{m_1}, \dots , t_l I_{m_l}) : t_i \in \GG_m , \prod_{i=1}^l t_i^{m_i} =1 \} \cong \GG_m^{l-1}.\]
The unipotent radical $U =\{ A \in P_\beta : A_{ii} = I_{m_i} \text{ for } 1 \leq i \leq l \} < P$ has lower central series given by filtering by off-diagonals:
\[U^{\{l-1\}}= \{ I \} < U^{\{l-2\}} < \cdots < U^{\{k\}}: = \{ A \in U : A_{ij} = 0 \text{ if }  j - i \leq k \} <\cdots < U^{\{0\}}= U. \]

\end{rmk}

\begin{defn}\label{def H and R beta}
For $P=P(\lambda)$ as above, we let $R := \prod_{i=1}^l \SL_{m_i}$ be the semisimple part of $L$. We define \[  H := U \rtimes R < P = U \rtimes L, \hspace{4pt}  \text{   and  }\hspace{10pt} \hH := H\times T \cong H\times \Gm^{l-1}.\]
\end{defn}
	
Note that there is a surjection $R \ra L/T$ with finite kernel. Similarly, there is a surjection $H \ra P/T$ with finite kernel. Consequently, from the perspective of GIT, quotienting by $\hH$ and $P$ are equivalent, since quotienting by finite groups poses no problem. Our construction will take a quotient of the $\hH$-action rather than of the $P$-action, in $l-1$ stages.

\subsection{Important notation for the parabolic group action}\label{sec important notn}
Suppose our parabolic group $P < \SL_N$ acts on an irreducible projective scheme $X$ with respect to an ample linearisation $\cL$. Before describing the Quotienting-in-Stages procedure for such an action, we give some definitions of subgroups of $P$ and loci in $X$ that will play a prominent role in the construction. After this in Definition \ref{defn quotienting-in-stages stable locus} we define the locus of which we will construct a quotient: the Quotienting-in-Stages stable locus, denoted $X^{P-qs}$.

\subsubsection{Subgroups of P}

Instead of using the lower central series of Remark \ref{rmk concrete desc parabolic and lcs}, which filters $U$ using the off-diagonals, we will consider a normal series for $U$ such that the successive quotients are unipotent radicals of maximal parabolic groups; that is, parabolics of length $l=2$. In fact, there are several ways in which we could do this, depending on a choice of an ordering of the subgroups of $U$ which are unipotent radicals of maximal parabolic subgroups containing $P$. 

\begin{defn}\label{def unipotent square brackets ordering}
For $1 \leq k \leq l-1$, define the subgroup 
\[U^{[k]}:=\{ A \in U : A_{ij} = 0 \text{ if }  j \leq  k  \text{ or } i > k\} < U.\] We let we $U^{(i)}:= \prod_{j \leq i} U^{[j]}$ to be the subgroup of $U$ generated by the first $i$ of these subgroups. Then we have the \emph{row filtration} of $U$, a normal series filtering $U$ by its rows from top to bottom
\[ \{I\} <U^{(1)} < \cdots < U^{(l-1)} =U\] 
whose successive quotients $U_{i}:=U^{(i)}/U^{(i-1)} \cong U^{[\sigma(i)]}/ U^{[i]} \cap U^{[i-1]}$ are abelian.
\end{defn}

\begin{rmk}
 Though we will not make use of this idea in the present paper, one could instead take any permutation $\sigma \in S_{l-1}$, which determines a different ordering $U^{[\sigma(1)]},\dots, U^{[\sigma(l-1)]}$ of these groups, and obtain a filtration by $U_{\sigma}^{(i)}:= \prod_{j \leq i} U^{[\sigma(j)]}$ with abelian subquotients. For example, the permutation $\sigma: i \mapsto l-i$ corresponds to filtering $U$ by its columns, rather than rows. 
 \end{rmk}

\begin{defn}\label{defn row filtr groups}
 Let us introduce some further subgroups and notation as follows. These groups will feature heavily in the Quotienting-in-Stages procedure outlined in $\S$\ref{sec outline quotient in stages} below.
\begin{enumerate}
\item For $1 \leq i < l=l(\lambda)$, define 1-PSs $\lambda^{(i)}$ and $ \lambda^{[i]} : \GG_m \ra Z(L)$ by
\[ \quad \quad \lambda^{(i)}(t) :=  \diag(t^{\beta_1} I_{m_1}, \dots , t^{\beta_i}I_{m_i}, t^{\beta_{>i}}I_{m_{>i}}) \quad \text{and} \quad \lambda^{[i]}(t) :=  \diag(t^{\beta_{\leq i}} I_{m_{\leq i}}, t^{\beta_{>i}}I_{m_{>i}})\]
where $m_{>i} = \sum_{j > i} m_j$ and $m_{\leq i} = \sum_{j \leq i} m_j$, and the exponents are defined by averaging: \[\beta_{>i} := \frac{\sum_{j>i} \beta_j m_j}{m_{>i}} \quad \quad \beta_{\leq i} := \frac{\sum_{j\leq i} \beta_j m_j}{m_{\leq i}}.\]
\item The unipotent group $U^{[i]}< U$ is the unipotent radical of the maximal parabolic subgroup $P^{[i]}:=P(\lambda^{[i]}) >P $ associated to the 1-PS $\lambda^{[i]}$. Consequently $\lambda^{[i]}$ grades $U^{[i]}$ and $U^{[i]}$ is normal in $P$, since it is normal in the larger group $P^{[i]}$.
\item The group $U^{(i)}< U$ is the unipotent radical of the parabolic subgroup $ P(\lambda^{(i)}) > P$ associated to the 1-PS $\lambda^{(i)}$ and thus $\lambda^{(i)}$ grades $U^{(i)}$.
\item We write $U^{[i-1,i]}:= U^{[i-1]}\cap U^{[i]} \leq U$. 

\item We obtain an increasing filtration $P^{(i)}:= U^{(i)}\rtimes L^{(i)}$ of $P$, where for $i < l-1$, \[L^{(i)}:= \{ A \in L : A_{jj} = I_{m_j} \text{ if }  j > i \}\]  and $L^{(l-1)}:=L$, giving an increasing filtration of $L$. We denote the successive quotients by $L_i := L^{(i)}/L^{(i-1)}$ and $U_i := U^{(i)}/U^{(i-1)}$ and $P_i := P^{(i)}/P^{(i-1)}$.
\item By restricting these filtrations to the semisimple part $R< L$, we obtain groups denoted $R^{(i)} < L^{(i)}$ with successive quotients $R_i$. Similarly, restricting the filtration to $H$ yields \[H^{(i)} = U^{(i)} \rtimes R^{(i)} < P^{(i)}\] with successive quotients $H_i = U_i \rtimes R_i$.
\item Let $T^{(i)}:= \prod_{j \leq i} \lambda^{[j]}(\GG_m) =\prod_{j \leq i} \lambda^{(j)}(\GG_m) \cong \Gm^{i}$ be the $i$ dimensional subtorus of $T= Z(L)$ generated by the first $i$ 1-PS's. We write $\hH^{(i)} := H^{(i)} \rtimes T^{(i)}$.
\item The 1-PS $\lambda^{[i]} : \GG_m \ra P$ descends to a 1-PS $\lambda_{i} : \GG_m \ra P /P^{(i-1)}$ of length two and the unipotent radical of the associated parabolic $P(\lambda_{i})$ is isomorphic to $U_{i}$. Furthermore, $\lambda_{i}$ grades $U_{i}$. We write $\hH_i :=   \hH^{(i)}/\hH^{(i-1)}\cong H_i \rtimes \lambda_i(\GG_m)$.
\end{enumerate}
\end{defn}

The basic idea of quotienting-in-stages is to quotient inductively by the groups $\hH_i$, since the stabilisers in these groups have nicer properties than those of the full $\hH$ at once. Moreover, at each stage of our quotienting process, we will use a different 1-PS in the central torus $T = Z(L)$ to grade $U^{[i]}$; consequently, we will only need to consider stabilisers for the semisimple part $R$ of $L$.

\begin{rmk} These filtrations have length $l-1$, but one could filter $L$ and $R$ (and thus also $P$ and $H$) in $l$ steps, using all $l$ of the block rows. However, since the above filtrations are only required to quotient by the non-reductive unipotent radical $U$, which naturally has a length $l-1$ row-filtration, we have made the filtrations of $L$ and $R$ the same length by combining the bottom two rows in the final step. Thus  $R_{i} = \SL_{m_i}$ for $1 \leq i \leq l-2$ and $R_{l-1} = \SL_{m_{l-1}} \times \SL_{m_{l}}$.
\end{rmk}

\subsubsection{Loci in $X$}

We continue to use the various subgroups of $P$ from Definition \ref{defn row filtr groups} above.

\begin{defn} \label{defn important loci}
We define the following subsets of $X$, by analogy with Definition \ref{def Zmin}.
	\begin{enumerate}
		\item The \emph{$[i]$-minimal weight space} is $Z^{[i]}_{\min} = Z(X,\lambda^{[i]})_{\min}$ and the \emph{$[i]$-attracting open set} is $X^{[i]}_{\min} := X(\lambda^{[i]})_{\min}$, which are related by the \emph{$[i]$-retraction} $p_{[i]} : X^{[i]}_{\min}  \ra Z^{[i]}_{\min} $ sending $x \mapsto \lim_{t \ra 0} \lambda^{[i]}(t) \cdot x$. The corresponding (semi)stable sets are $Z^{[i],(s)s}_{\min} := (Z^{[i]}_{\min})^{R^{(i)}-(s)s}$ and $X^{[i],(s)s}_{\min}:=p_{[i]}^{-1}(Z^{[i],(s)s}_{\min} )$. 

		\item The \emph{$(i)$-minimal weight space} is $Z^{(i)}_{\min} = Z(X,\lambda^{(i)})_{\min}$ and the \emph{$(i)$-attracting open set} is $X^{(i)}_{\min} := X(\lambda^{(i)})_{\min}$. We define the \emph{$(i)$-retraction} $p_{(i)} : X^{(i)}_{\min}  \ra Z^{(i)}_{\min} $ sending $x \mapsto \lim_{t \ra 0} \lambda^{(i)}(t) \cdot x$, and the corresponding (semi)stable sets $Z^{(i),(s)s}_{\min} := (Z^{(i)}_{\min})^{R^{(i)}-(s)s}$ and $X^{(i),(s)s}_{\min}:=p_{(i)}^{-1}(Z^{(i),(s)s}_{\min} )$ in an analogous way. 
		
		\item As a special case of (2), for $\lambda = \lambda^{(l-1)}$, the \emph{$\lambda$-minimal weight space} is $Z_{\min} = Z(X,\lambda)_{\min}$ and the \emph{$\lambda$-attracting open set} is $X_{\min} := X(\lambda)_{\min}$, which are related by the \emph{$\lambda$-retraction} $p : X_{\min}  \ra Z_{\min} $ sending $x \mapsto \lim_{t \ra 0} \lambda(t) \cdot x$. The corresponding (semi)stable sets are $Z^{(s)s}_{\min} := (Z_{\min})^{R-(s)s}$ and $X^{(s)s}_{\min}:=p^{-1}(Z^{(s)s}_{\min})$.
	\end{enumerate}
\end{defn}

\brm As a rough conceptual shorthand, it may help the reader to bear in mind the case of objects in an abelian category where $\lambda$ corresponds to a length $l$ filtration: then $Z^{[i]}_{\min}$ is those objects for which the filtration is split into two pieces \emph{at the $i$th point}, and $Z^{(i)}_{\min}$ is those objects for which the filtration is split \emph{at each of the first $i$ steps}. The maps $p_{[i]}$ and $p_{(i)}$ then correspond to taking a kind of \lq partial associated graded'; that is, an associated graded with respect to certain coarsenings of the original filtration defined by $\lambda$ (see $\S$\ref{sec desc important loci sheaves} below). \erm

\begin{rmk}\label{rmk square and round minimal weight spaces}\
	\begin{enumerate}
 \item Note that here $Z^{(s)s}_{\min}$ is defined to be the (semi)stable set for the semisimple group $R < L$ rather than the group $\overline{L} = L /\lambda(\GG_m)$ as in Definition \ref{def Zmin}. If $T$ acts trivially on $Z_{\min}$ (which will later follow from Assumption \hyperlink{QiS}{(QiS)}; see Remark \ref{rmk centre of L acts trivially on Zmin}), then $Z^{ss}_{\min} :=Z^{R-ss}_{\min}= Z^{\overline{L}-ss}_{\min}$, so for the \emph{semistable}s locus it makes no difference which group we use; however, for $l>2$, we have $Z^{\overline{L}-s}_{\min} = \emptyset$, while $Z^{s}_{\min} :=Z^{R-s}_{\min}$ may be non-empty.

		\item We can also extend the definition of $\lambda^{[i]}$ in Definition \ref{defn row filtr groups} to $i = l$. Then $\lambda^{[l]}$ is trivial and $p_{[l]} : X = X^{[l]}_{\min} \ra X = Z^{[l]}_{\min}$ is the identity map.
		\item If we are interested in an unstable HKKN stratum $S_\beta \cong G \times^{P_\beta} Y_\beta^{ss}$ associated to the reductive group $G = \SL_N$ acting on a projective variety $Y$, where the parabolic subgroup $P_\beta$ is determined by a 1-PS $\lambda_\beta$ of $G$, then we can apply this procedure to the $P_\beta$-action on $X := \overline{Y_\beta}$. Since $\lambda= \lambda_\beta$, we have that $p : X_{\min}  \ra Z_{\min}$ is $p_\beta : Y_\beta \ra Z_\beta$ and $Z^{ss}_{\min} := Z_\beta^{R-ss} = Z_\beta^{ss}:=Z_\beta^{L_\beta-ss}$ if $Z(L_\beta)$ acts trivially on $Z_\beta$.
	\end{enumerate}
\end{rmk}

\subsubsection{The Quotienting-in-Stages stable set} Now we have introduced all this notation, we can define a subset of $X$ with an explicit Hilbert--Mumford type description. We will eventually prove that this set admits a geometric $P$-quotient under certain assumptions.

\begin{defn} \label{defn quotienting-in-stages stable locus} 
The \emph{Quotienting-in-Stages stable locus} for the linearised $P$-action on $X$ is
\[X^{P-qs} :=\left\{ 
\begin{array}{ll}
 &  p_{[i]}(x) \notin U^{[i-1]}Z^{[i-1]}_{\min} \text{ for all } 2\leq  i\leq l \\
x\in X^{s}_{\min} : & \dim \stab_{U^{[i]}}(p_{[i]}(x)) = d^{[i]}_{\min}  \text{ for all } 1\leq  i\leq l-1\\
& \dim \stab_{U^{[i-1,i]}}(p_{[i]}(x)) = d^{[i-1,i]}_{\min}  \text{ for all } 2\leq  i\leq l-1
\end{array}\right\},\]
where we recall $U^{[i-1,i]}:= U^{[i-1]} \cap U^{[i]}$ and, using the notation of Definition \ref{def centres of blowups}, we set
\[ d^{[i]}_{\min}:= d_{\min}(X^{[i],ss}_{\min},U^{[i]}) \quad \text{ and } \quad d^{[i-1,i]}_{\min}:= d_{\min}(X^{[i],ss}_{\min},U^{[i-1,i]}),\]
which are the generic stabiliser dimensions of $U^{[i]}$ and $U^{[i-1,i]}$ on $X^{[i],ss}_{\min}$ respectively.
\end{defn}

\begin{rmk}
We note that $U^{[i]}Z^{[i]}_{\min} = U^{(i)}Z^{[i]}_{\min}$, as $U^{[i]} \subset U^{(i)}$ and the complement lies in the Levi $L^{[i]}$ of $P(\lambda^{[i]})$ and, as $\lambda^{[i]}$ is central in $L^{[i]}$, this Levi preserves $Z^{[i]}_{\min}$.
\end{rmk}

A natural question is whether, or under what conditions, the locus $X^{P-qs}$ is non-empty. In Remark \ref{rmk QiS stable set is open} below, we will see that, under the so-called \emph{Quotienting-in-Stages Assumption} \hyperlink{QiS}{(QiS)} and the \emph{Upstairs Unipotent Stabiliser Assumption} \ref{starUstagesibullet} introduced in \S\ref{def star conditions stages}, the locus $X^{P-qs}$ is a non-empty open subset of $X_{\min}$. In fact, under the \emph{Upstairs Unipotent Stabiliser Assumption}, the stabiliser conditions in this definition hold automatically, and so
\[ X^{P-qs} = \{x\in X^{s}_{\min} \: : \:p_{[i]}(x) \notin U^{[i-1]}Z^{[i-1]}_{\min} \text{ for all } 2\leq  i\leq l \}.\]

\subsubsection{Example: unstable sheaves}\label{sec desc important loci sheaves}

Let $\tau = (P_1, \dots, P_l)$ be a length $l \geq 2$ HN type for sheaves on a polarised projective scheme $(B,\cO_B(1))$. As in $\S$\ref{sec sheaves nonred git setup}, we consider the action of $P_\tau = P(\lambda_\tau) = U_\tau \rtimes L_\tau < \SL_{P(n)}$ on the closure $X:=\overline{Y_\tau}$ of $Y_\tau$ in $Q_n$. Recall that the canonical linearisation $\cL_\tau$ is borderline for the grading of $U_\tau$ given by the 1-PS $\lambda_\tau$. The Quotienting-in-Stages procedure will determine a small perturbation of this linearisation.

The (rational) 1-PS $\lambda_\tau$ determines a length $l$ filtration of $V =\CC^{P(n)}$
\begin{equation}\label{filt V}
	0 = V^{(0)} \subset V^{(1)}=\CC^{P_1(n)} \subset \cdots \subset V^{(i)}=\CC^{P_{\leq i}(n)} \subset \cdots \subset V^{(l)}=\CC^{P(n)}
\end{equation}
where $P_{\leq i}:= \sum_{j \leq i} P_j$ and $P_{> i}:= \sum_{j > i} P_j$. For a quotient sheaf $q_{\cE} : \cO_B(-n)^{\oplus P(n)} \twoheadrightarrow \cE$ in $Q_n$, the filtration \eqref{filt V} induces a length $l$ filtration of $\cE$
\begin{equation}\label{filt E}
	0 = \cE^{(0)} \subset  \cE^{(1)} \subset \cdots  \subset \cE^{(i)}  = q_{\cE}(V^{(i)}\otimes \cO_B(-n))\subset \cdots \subset \cE^{(l)}=\cE
	\end{equation}
and we write $\cE_i:=\cE^{(i)}/\cE^{(i-1)}$ for the successive subquotients; this filtration is the HN filtration of $\cE$ if $q_{\cE} \in Y_\tau^{ss}$. The limit of $\lambda(t) \cdot q_{\cE}$ as $t \ra 0$ corresponds to the associated graded sheaf $\gr(\cE^{(\bullet)})=\oplus_i \cE_i$ of the filtration \eqref{filt E}; see \cite[Lemma 4.4.3]{HL} and \cite[$\S$5.1]{HK}. Thus the morphism $p_\tau : Y_\tau^{ss} \ra Z_\tau^{ss}$ takes a sheaf to its HN-graded sheaf. In particular, we have the following statements concerning the minimal weight space for $\lambda_\tau$ (see also Propositions \ref{prop nonred GIT for HKKN strata} and \ref{reductive sheaf stabs}).

\begin{lemma}\label{lemma Z for sheaves}
For the $P_\tau=P(\lambda_\tau)$-action on $X:=\overline{Y_\tau}$, the following statements hold.
\begin{enumerate}
\item The retraction $p : X_{\min} \ra Z_{\min}$ coincides with the retraction $p_\tau : Y_\tau \ra Z_\tau$.
\item The semisimple part of $R_\tau = \Pi_{i=1}^l \SL_{P_i(n)}$ of $L_\tau$ has stable locus 
\[ Z_\tau^{s}:= (Z_\tau)^{R_\tau-s} \cong Q^{s}_n(P_1) \times \cdots \times Q^{s}_n(P_l)\]
where $Q^{s}_n(P_i)$ is the open subscheme of $Q_n(P_i):=\quot(\cO_B(-n)^{\oplus P_i(n)},P_i)$ of quotient sheaves $q_i : \cO_B(-n)^{\oplus P_i(n)} \surj \cE_i$ such that $\cE_i$ is stable and $H^0(q_i(n))$ is an isomorphism. 
\item We have $X^{ss}_{\min} = Y_\tau^{ss}$ and $X_{\min}^s = Y_\tau^{s}:= p_\tau^{-1}(Z_\tau^s)$ is the locus of quotient sheaves $q_{\cE} \in Y_\tau^{ss}$ such that each subquotient $\cE_i$ in the HN filtration of $\cE$ is stable.
\end{enumerate}
\end{lemma}

Next consider the 1-PS $\lambda_{\tau}^{[i]}(t)$ with only two weights; this induces length 2 filtrations 
\begin{equation}\label{filt E two}
 V^{[\bullet, i]}:=(0 \subset V^{(i)} \subset V) \quad \text{ and } \quad \cE^{[\bullet, i]}:=(0 \subset \cE^{(i)} \subset \cE),
 \end{equation}
which are coarsenings of \eqref{filt V}  and \eqref{filt E}. The limit of $\lambda^{[i]}(t) \cdot q_{\cE}$ as $t \ra 0$ corresponds to the associated graded sheaf $\cE^{(i)} \oplus \cE/\cE^{(i)}$.

\begin{lemma}\label{lem Z[i] for sheaves}
For $m \gg 0$ depending on $\tau$, we have for all $1 \leq i \leq l-1$, an
a $L_\tau$-equivariant isomorphism  
\[Z_{\min}^{[i]} \cong  X \times_{Q_n(P)} (Q_n(P_{\leq i}) \times Q_n(P_{>i})). \]
Moreover a quotient sheaf $q_{\cE} \in Y_\tau^{ss}$ lies in $U_\tau^{[i]}Z_{\min}^{[i]} $ if and only if $\cE \cong \cE^{(i)} \oplus \cE/\cE^{(i)}$, i.e.\ the inclusion of the $i$th sheaf in the HN filtration $\cE^{(i)} \subset \cE$ is split.

\end{lemma}
\begin{proof}
As outlined above, the components of the $\lambda^{[i]}_\tau(\GG_m)$-fixed locus are products of two Quot schemes, corresponding to the two sheaves in the associated graded of the filtration induced by this 1-PS, and indexed by pairs of Hilbert polynomials $(P',P'')$ such that $P = P' + P''$. One such component is given by the pair $(P_{\leq i}, P_{>i})$ and we claim this is the minimal weight space. Since $Y_{\tau}^{ss}$ is open in $X$, the minimal weight for $\lambda_{\tau}^{[i]}$ acting on $X$ is the same as the minimal weight for the action on $Y_{\tau}^{ss}$. For any $q_\cE \in Y_{\tau}^{ss}$, a quick calculation shows the Hilbert--Mumford weight of $\lambda^{[i]}$ with respect to the twisted linearisation $\cL_\tau$ is equal to zero; this computation follows from \cite[Lemma 4.4.4]{HL} and is similar to \cite[$\S$4.3]{Hoskins}.  By taking $m$ sufficiently large, we can ensure that the only pair of polynomials $(P',P'')$ giving this minimum weight is $(P_{\leq i}, P_{>i})$. 

The final statement follows, as any sheaf quotient sheaf in $U_\tau^{[i]}Z_{\tau,\min}^{[i]} =P_\tau Z_{\tau,\min}^{[i]}$ is isomorphic to a sheaf in $Z_{\min}^{[i]}$ and this is precisely the locus where $\cE^{(i)} \subset \cE$ is split.
\end{proof}

Similarly to Lemmas \ref{lemma Z for sheaves} and \ref{lem Z[i] for sheaves} one can prove the following result.

\begin{lemma} \label{lemma Z(i) for sheaves} The minimal weight space for the 1-PS $\lambda^{(i)}_\tau$ admits the following description.
	\bnu \item For $m \gg 0$ depending on $\tau$, we have for all $1 \leq i \leq l-1$, an
a $L_\tau$-equivariant isomorphism  
\[Z_{\min}^{(i)} \cong  X \times_{Q_n(P)} (Q_n(P_{1}) \times \dots \times Q_n(P_i) \times Q_n(P_{>i})). \]
Moreover a quotient sheaf $q_{\cE} \in Y_\tau^{ss}$ lies in $U_\tau^{(i)}Z_{\min}^{(i)} $ if and only if $\cE \cong \cE_1 \oplus \dots \cE_i \oplus \cE/\cE^{(i)}$.
\item  The (semi)stable locus for the group  $R^{(i)}$ are given by
\[Z^{(i),(s)s}_{\min} := (Z_{\min}^{(i)})^{R_\tau^{(i)}-(s)s} \cong Q^{(s)s}_n(P_1) \times \cdots \times Q^{(s)s}_n(P_i) \times Q_n(P_{\geq i}).\]
\enu 
\end{lemma}

Let us now determine the Quotienting-in-Stages stable locus in the case of sheaves. Initially we ignore the conditions on unipotent stabilisers and focus on the first part of the definition.

\begin{prop}\label{prop first part of QiS stable for sheaves}
For $q_{\cE} \in X^{s}_{\min} = Y_\tau^s$ and $2\leq  i\leq l$, we have $p_{[i]}(q_\cE) \in  U_{\tau}^{[i-1]}Z^{[i-1]}_{\min}$ if and only if $\cE^{(i-1)} \subset \cE^{(i)}$ is split. 
\end{prop}
\begin{proof}
For $i <l$, we have $p_{[i]}(q_{\cE}) = q_{\cE_{\leq i}}\oplus q_{\cE_{> i}}$ with corresponding quotient sheaf $\cE^{(i)} \oplus \cE/\cE^{(i)}$. The 1-PS $\lambda^{[i-1]}_\tau(t)$ induces the following filtration of $\cE^{(i)} \oplus \cE/\cE^{(i)}$
\[ 0 \subset \cE^{(i-1)} \subset \cE^{(i)} \oplus \cE/\cE^{(i)}, \]
which is split (i.e. $p_{[i]}(x) \in U_\tau^{[i-1]}Z^{[i-1]}_{\min}$) if and only if $\cE^{(i-1)} \subset \cE^{(i)}$ is split. For $i = l$, we have $p_{[l]} = \mathrm{Id}$ and so $p_{[l]}(x) \in U_{\tau}^{[l-1]}Z^{[l-1]}_{\min}$ if and only if $\cE^{(l-1)} \subset \cE^{(l)}=\cE$ is split.
\end{proof}

The next step is to describe the unipotent stabilisers for $U_\tau^{[i]}$ and $U_\tau^{[i-1,i]}$.

\begin{prop}\label{prop QiS unipotent stab for sheaves}
For any $q_\cE \in Y^{ss}_\tau$, the following statements hold.
\begin{enumerate}
\item For $1 \leq i \leq l-1$, we have $\Lie \stab_{U_{\tau}^{[i]}}(q_\cE) \cong \Hom(\cE/\cE^{(i)},\cE^{(i)})\cong \Lie \stab_{U_{\tau}^{[i]}}(p_{[i]}(q_\cE))$.
\item For $2 \leq i \leq l-1$, $\Lie \stab_{U_{\tau}^{[i-1,i]}}(q_\cE) \cong \Hom(\cE/\cE^{(i)},\cE^{(i-1)}) \cong \Lie \stab_{U_{\tau}^{[i-1,i]}}(p_{[i]}(q_\cE))$.
\end{enumerate}
\end{prop}
\begin{proof}
Let us prove the first statement as the second follows similarly. Similarly to Proposition \ref{prop uni stab sheaves} (i), we have $\Lie \stab_{U_{\tau}^{[i]}}(q_\cE) \cong \End_{-1}(\cE^{[\bullet, i]})$ for the filtration $\cE^{[\bullet, i]}:=(0 \subset \cE^{(i)} \subset \cE)$. Since this is a length 2 filtration, this filtered endomorphism group is the homomorphism group from the quotient sheaf to the subsheaf, which gives the first isomorphism. The second isomorphism follows as $p_{[i]}(q_\cE)$ has corresponding quotient sheaf $\cF:=\cE^{(i)}\oplus \cE/\cE^{(i)}$ and the induced length two filtration on $\cF$ is $0 \subset \cF^{(i)} = \cE^{(i)} \subset \cF$ with $\cF/\cF^{(i)}  \cong \cE/\cE^{(i)}$.
\end{proof}

We can characterise the quotienting-in-stages stable locus for sheaves by introducing a notion of $\tau$-stability for sheaves of HN type $\tau$ generalising the length 2 notion in \cite{Josh_length2}.

\begin{defn}\label{def tau stable} 
We say a sheaf of $\cE$ of HN type $\tau$ is \emph{$\tau$-stable} if the following conditions hold:
\begin{enumerate}
\item Each successive quotient $\cE_i$ in the HN filtration is stable,
\item Each inclusion $\cE^{(i-1)} \subset \cE^{(i)}$ in the HN filtration is non-split,
\item For $1\leq i \leq l-1$ and $2\leq j \leq l-1$ we have \[ \dim\Hom(\cE/\cE^{(i)},\cE^{(i)}) = d_{\min}^{[i]} \quad \text{and} \quad  \dim\Hom(\cE/\cE^{(j)},\cE^{(j-1)}) = d^{[j-1,j]}_{\min}.\]
\end{enumerate} 
\end{defn}

\begin{cor} \label{cor what is X P-qs for sheaves}
	For, $q_\cE\in Y_\tau^{ss}$, we have $q_\cE \in X^{P-qs}$ if and only if $\cE$ is $\tau$-stable.
\end{cor}

\subsection{The Quotienting-in-Stages procedure}\label{sec outline quotient in stages} Let us return to the more general case described in Definition \ref{def unipotent square brackets ordering}, in which $P<\SL_N$ is a parabolic subgroup of length $l$ acting on a projective scheme $X$ with respect to an ample linearisation $\cL$. Our approach is to construct a quotient of this action in $l-1$ stages, where we use $\lambda^{[i]} < T$ to grade the unipotent radical $U^{[i]}$ at the $i$th stage. This gives an inductive procedure, where the base case and inductive step look like quotienting by a group whose unipotent radical is the unipotent radical of a maximal parabolic group; in particular, this unipotent radical is abelian and so we can apply the results in $\S$\ref{sec easy blowups abelian}.

\subsubsection{The base case and inductive step: quotients by length 2 parabolics}\label{sec quotient parabolic of length 2}
The base case and inductive step in the Quotienting-in-Stages procedure both involve quotienting by a unipotent radical of a maximal (i.e.\ length $2$) parabolic group, which is studied in \cite{Josh_length2}. If $l =2$, we have $P = P(\lambda)$ for a 1-PS $\lambda$ with two distinct weights $\beta_1 > \beta_2$ with multiplicities $(m_1,m_2)$. The central torus $T \cong \GG_m$ of $L$ is generated by $\lambda$, which grades the (abelian) unipotent radical $U$. Thus we can directly apply the results of non-reductive GIT: in length $2$, the surjection $R \ra \overline{L}$ has finite kernel and so (semi)stability for these groups coincides.

In length $l=2$, suppose that semistability coincides with stability for $\overline{L}$ (or equivalently $R$); thus $Z^{s}_{\min} = Z^{ss}_{\min}$. Additionally suppose that \ref{starURss} holds (note that as $U$ is abelian, this requires $\dim \stab_U(-)$ to be constant on $X^{ss}_{\min}$). Then by Theorem \ref{mainthm nrGIT quotient}, we obtain a projective geometric $P$-quotient 
\[ q : X^{P-s} = X^{s}_{\min} \setminus U Z^{s}_{\min} \ra X /\!/ P. \] 

If \ref{starURss} fails in length $l = 2$, let $d_{\min}$ be the minimal dimension of $\dim \stab_U(-)$ on $X^{ss}_{\min}$ and we perform unipotent blow-ups. If \ref{CZminNE} holds, we will not blow-up all of the minimal weight space and by Proposition \ref{prop easy blowups abelian} (and \cite{Josh_length2}), the stable locus of Definition \ref{def hat stable locus easy abelian case},
\[ X^{\hH-\hs} =  \{ x \in X^s_{\min} \setminus U Z^s_{\min} : \dim \stab_{U}(p(x)) = d_{\min} \}, \]
admits a geometric and quasi-projective $\hH_\beta$-quotient. This quotient is projective if additionally
	\[ \dim \stab_{U}(x) = d_{\min} \text{ for all } x \in Z_{\min}^{ss}.\]

\subsubsection{Quotienting-in-Stages for $l>2$}

Recall the groups introduced in Definition \ref{defn row filtr groups}. We will construct a quotient in stages so that after stage $i$, we will have quotiented by the group $\hH^{(i)} = H^{(i)} \rtimes T^{(i)}$, and to go from stage $i-1$ to stage $i$, we take a quotient by the group $\hH_i = H_i \rtimes \lambda_i(\GG_m)$, where the graded unipotent part of this action $U_i \rtimes \lambda_i(\GG_m) $ is induced by an action of $U^{[i]} \rtimes \lambda^{[i]}(\GG_m)$ on $X$.

\begin{const}\label{const qnt in stages} (\emph{Quotienting-in-Stages Construction})
	
\noindent \textit{Base stage:} Let $X_1:= X$ and suppose that the action of $\hH_1 = \hH^{(1)}$ satisfies all the assumptions of Theorem \ref{mainthm nrGIT quotient}, so we can take a NRGIT quotient giving a rational map
\[q_1 : X_1 \dashrightarrow X_2\] 
to a projective scheme $X_2$, which restricts to a geometric quotient on its domain of definition.

\noindent \textit{Inductive step:} Assume we have constructed successive quotients $q_j : X_j \ra X_{j+1}$ of the $\hH_{j}$-action on $X_{j}$ for all $j < i$, whose composition 
\[ q_{(i-1)} := q_i \circ \cdots \circ q_1 : X \dashrightarrow X_{i} \]
is a $\hH^{(i-1)}$-quotient of an open subset of $X$. On $X_i$, there is an induced linearised action of $\hH_{i}$ coming from the linearised $\hH^{(i)}$-action on $X_1$ and assuming we can apply Theorem \ref{mainthm nrGIT quotient}, we can take a NRGIT quotient of this $\hH_{i}$-action \[q_{i} : X_{i} \dashrightarrow X_{i+1}\] and obtain a quotient of an open subset of $X$ by defining 
\[q_{(i)} = q_i \circ \cdots \circ q_1 : X \dashrightarrow X_{i+1}\] to 
be the composition. Then $\dom(q_{(i)}) \subset X$ is open subset admitting a projective geometric $\hH^{(i)}$-quotient given by the map $q_{(i)}$. 

\noindent \textit{Final step (l-1):} On the $\hH_{l-2}$-quotient $X_{l-1}$, there is an induced linearised action of $\hH_{l-1}$ coming from the $\hH^{(l-1)}$-action on $X$. Again assuming we can apply Theorem \ref{mainthm nrGIT quotient}, we can take a $\hH_{l-1}$-quotient $q_{l-1} : X_{l-1} \dashrightarrow X_{l}$, and let \[ q:= q_{(l-1)} = q_{l-1} \circ \cdots \circ q_1 : X_1 \dashrightarrow X_{l}\] be the composition. Thus $q|_{\dom(q)} : \dom(q) \ra X_{l}$ is a projective geometric $\hH^{(l-1)}$-quotient. 
\end{const} 

If at each stage we are in the best case covered by Theorem \ref{mainthm nrGIT quotient}, then the domain of definition of each $q_i$ admits an explicit Hilbert--Mumford description as the $\hH_i$-stable locus in $X_i$ given by Definition \ref{def nonred ss}. This procedure consists of $l-1$ quotients:
\[ \xymatrix{ X_{1} \ar@{-->}[r] & X_{2} \ar@{-->}[r] & X_{3} \cdots \cdots \ar@{-->}[r] & X_{l-1} \ar@{-->}[r] & X_{l} \\
		X_{1}^{\hH_{1}-s} \ar@{_{(}->}[u] \ar@{->>}[ru]^{q_1} & X_{2}^{\hH_{2}-s} \ar@{_{(}->}[u] \ar@{->>}[ru]^{q_2} & \cdots \cdots \cdots  & X_{l-1}^{\hH_{l-1}-s} \ar@{_{(}->}[u] \ar@{->>}[ru]^{q_{l-1}}
	} \]
such that the composition $q : X \dashrightarrow X_l$ gives a projective geometric $P$-quotient of its domain of definition. 

In $\S$\ref{sec proof}, we will explicitly determine the domain of definition of $q$ under certain assumptions.  A large part of $\S$\ref{sec proof} is devoted to understanding what happens when the conditions of Theorem \ref{mainthm nrGIT quotient} fail at some stage. The basic idea is to perform a sequence of (reductive or non-reductive) blow-ups on $X$ similarly to $\S$\ref{sec red part desing} and $\S$\ref{sec nrGIT blowups}. This process results in a blow-up of $X$ for which the above construction can be carried out as written, and it remains to relate the resultant stable locus on this blow-up with the original $X$.

Before we describe this sequence of quotients more carefully, we pause to remark on the choices we are making at each stage.
	
\begin{rmk}\label{rmk choices filtering}  
At each stage, one has to twist the induced linearisation by a (rational) character in order for it to be well-adapted. Since the reductive groups $R_i$ are products of special linear groups and have no non-trivial characters, the twisting character can only be non-trivial on the grading 1-PS $\lambda_{i}(\GG_m)$ and thus the space of twisting characters is 1-dimensional. In this 1-dimensional space of characters, only a small interval will give characters that make the linearisation well-adapted. Thus, there is in effect no choice of character\footnote{Technically, at the next stage, the possible choices of character needed to shift will depend on the previous choices, but this is purely bookkeeping and does not effect the quotient.}, since we are always constrained to lie in a certain NRVGIT chamber in the sense of \cite{Berczi2018a}.  In fact, we will soon see in \S\ref{subsec a single character works for the whole qnt} that, in the case we consider, there is a single character of $\hH$ that induces well-adapted characters of each $\hH_i$ at each stage. Hence, we may as well choose this character from the beginning and fix the linearisation once and for all, removing the need to twist at any subsequent point.
\end{rmk}

\subsubsection{The Downstairs Stabiliser Assumption (D)} \label{subsec downstairs suffice to obtain a quotient}
Our next task is to precisely state the assumptions we need at each stage in order to apply Theorem \ref{mainthm nrGIT quotient} and execute Construction \ref{const qnt in stages}. First we need some more notation.

\begin{defn}
	\label{def stable loci stage i}
	For the linearised $\hH_i$-action on $X_i$, we introduce the following notation:
	\begin{enumerate}
		\item The \emph{stage $i$ minimal weight space} is $Z_{i,\min} := Z(X_i, \lambda_i)_{\min}$,
		\item The \emph{stage $i$ attracting open set} is $X_{i,\min} :=X_i(\lambda_i)_{\min}$,
		\item The \emph{stage $i$ retraction} is $p_i : X_{i,\min} \ra Z_{i,\min}$ is given by $p_i(x):= \lim_{t \ra 0} \lambda_i(t) \cdot x$,
		\item The \emph{stage $i$ (semi)stable minimal weight space} is $Z_{i,\min}^{(s)s}:=Z_{i,\min}^{R_i-(s)s}$,
		\item The \emph{stage $i$ (semi)stable locus} is $X_{i,\min}^{(s)s} := p_i^{-1}(Z_{i,\min}^{(s)s})$.
	\end{enumerate}
\end{defn}

With this notation introduced, we can state the conditions we will assume when applying Construction \ref{const qnt in stages}.

\begin{ass} \hypertarget{D}{(\emph{The Downstairs Stabiliser Assumption})} \label{stab stages assumption}
\begin{enumerate}
\item We say the \emph{downstairs unipotent stabiliser assumption} (DU) holds if 
		\begin{equation*}\label{hUbullet,Rss} 
		\forall \: 1 < i < l \: : \: \:
		\dim \stab_{U_i}(-) \text{ is constant on }  X_{i,\min}^{ss} \tag*{(DU)},
	\end{equation*}
	\item We say the \emph{downstairs reductive stabiliser assumption} (DR) holds if
	\begin{equation*}\label{Rbulletzero}
				\forall \: 1 < i < l \: : \: \: \dim \stab_{R_i}(-) = 0 \text{ on }  Z_{i,\min}^{ss}. \tag*{(DR)}
	\end{equation*}
	\item We say the \emph{$R$-stable locus and semistable locus coincide} if $X^{ss}_{\min} = X^{s}_{\min}$, or equivalently:
	\[ \dim \stab_{R}(-) = 0 \text{ on }  Z_{\min}^{ss}. \] 
	
		\end{enumerate}
We say the \emph{Downstairs Stabiliser Assumption} \hyperlink{D}{(D)} holds if all the above assumptions hold.
\end{ass}

If \hyperlink{D}{(D)} holds, Theorem \ref{mainthm nrGIT quotient} tells us that  $X_i^{\hH_i-s}=X_{i,\min}^{s} \setminus U_i Z_{i,\min}^{s}$ admits a projective geometric $\hH_i$-quotient $ q_i : X_i^{\hH_i-s} \ra X_{i+1}$.

\bdf \label{def qis map}  \label{def qis s map}
If \hyperlink{D}{(D)} holds, we define the \emph{Quotienting-in-Stages map} \[q:=q_{(l -1)} := q_{l -1} \circ \cdots \circ q_1 : X=X_1 \dashrightarrow X_{l}\]
to be the composition of the rational maps obtained in the above construction. The following restriction of $q$ will also play an important role in our argument: we further define \[q^{s}= q\mid_{X^s_{\min}} X^{s}_{\min} \dashrightarrow X_{l}\] to be the restriction of $q$ to $X^s_{\min}$. In the same way, we define each $q^s_{(i)}$ to be the restriction of $q_{(i)}$ to $X^s_{\min}$, and $q_i^s$ to be the restriction of $q_i$ to the open locus $q_{(i-1)}(X^s_{\min}) \subset X_i$.
\edf

One of our main aims for the remainder of this paper will be to show that, under certain assumptions, the domain of definition of $q$ is equal to the more explicit quotienting-in-stages stable locus $X^{P-qs}$ defined at Definition \ref{defn quotienting-in-stages stable locus}.

\brm We will later show in Corollary \ref{cor domain q in Xsmin} that $\dom(q) \subset X^{s}_{\min}$,  The latter will follow from the fact that $T$ acts trivially on $Z_{\min}$, which is a consequence of Assumption \hyperlink{QiS}{(QiS)} below.  \erm

The results of this subsection are summarised in the following corollary.

\begin{cor} \label{cor exists open set with qnt}
	If \hyperlink{D}{(D)} holds, then Construction \ref{const qnt in stages} yields an open subset $\dom(q) \subset X$ that has a projective geometric $P$-quotient.
	
\end{cor}
\begin{proof}
Since each $q_i$ is a geometric quotient of its domain of definition by Theorem \ref{mainthm nrGIT quotient}, the morphism $q : X \dashrightarrow X_{l}$ gives a projective geometric $P$-quotient of $\dom(q) \subset X$. \end{proof}

\subsubsection{Constructing the quotient from a single character} \label{subsec a single character works for the whole qnt}

Recall the discussion on twisting the linearisation at each stage in Remark \ref{rmk choices filtering}. We will now prove that, under the hypothesis \hyperlink{D}{(D)}, we can find a single character $P$ such that the Quotienting-in-Stages morphism $q$ is given by taking invariant sections of this twisted linearisation. 

The character groups of $P$ and its Levi $L$ are isomorphic; let us explicitly write down the latter. Recall that the 1-PS $\lambda$ has $l$ distinct weights $\beta_1 > \dots > \beta_m$ with multiplicities $m_i$ and $L \cong (\prod_{i=1}^l \GL_{m_i}) \cap \SL_N$. Any (rational) character $\chi : L \ra \GG_m$ is given by $(r_1, \dots , r_l) \in \QQ^l$ satisfying $\sum_{i=1}^l  m_i r_i = 0$; namely, we have
\[ \chi(g_1, \dots ,g_l) = \prod_{i=1}^l \det(g_i)^{r_i}.\]
Let us write $\chi_{\lambda}$ for the character corresponding to the rational tuple $(\beta_1, \dots , \beta_l)$ of weights of $\lambda$. Let $\cL_0$ be obtained from $\cL$ by twisting by a suitable rational multiple of $\chi_{\lambda}$ to obtain a borderline linearisation in the sense of Definition \ref{def borderline adapted}. For any character $\chi$, let $\cL_\chi$ be the linearisation obtained by twisting $\cL_0$ by $\chi$. 

\begin{prop}\label{prop character for Q-i-S}
If \hyperlink{D}{(D)} holds, then there is a small rational character $\chi : P \ra \GG_m$ such that $\cL_\chi$ induces well-adapted linearisations at each stage of Construction \ref{const qnt in stages}. Hence, we have:

\bnu \item The ring of $P$-invariant sections is finitely generated:
\[ R(X,\cL_\chi)^P \subset R(X,\cL_\chi) := \bigoplus_{r \geq 0} H^0(X,\cL_\chi^{\otimes r})\] 

\item The Quotienting-in-Stages map $q : X = X_1 \dashrightarrow X_l$  coincides with the map
\[ q_{H,\chi} :  X \dashrightarrow X/\!/_{\cL_\chi} H :=\proj(R(X,\cL_\chi)^P).\]
In particular, this latter map is a good $P$-quotient of its domain of definition
\[ X^{P-ss}(\cL_\chi) := \{ x \in X : \exists \sigma \in H^0(X,\cL_\chi^{\otimes r})^P \text{ for } r > 0 \text{ with } \sigma(x) \neq 0 \}. \]
\item Furthermore, $\dom(q) = X^{P-ss}(\cL_{\chi}) \subset X^{P-ss}(\cL_0) = p^{-1}(Z^{\overline{L}-ss}_{\min})$ and we obtain a map
\[X_l \cong X/\!/_{\cL_\chi} H \ra X/\!/_{\cL_0} H \cong  Z_{\min}/\!/_{\cL_0} L.\]
\enu 
\end{prop}
\begin{proof}
Since \hyperlink{D}{(D)} holds, we do not need to perform any blow-ups in the Quotienting-in-Stages construction and thus, at each stage, we just need to perturb the linearisation by twisting by a character of the quotienting groups $\hH_i$ to make it well-adapted, so that we can apply Theorem \ref{mainthm nrGIT quotient} at each stage. It thus suffices to show the choices of rational character used to perturb the borderline linearisation $\cL_0$ at each stage can be lifted to a single choice of rational character $\chi : P \ra \GG_m$ such that 
\[ X_l \cong \proj(R(X,\cL_\chi)^P). \]
Once we have shown this, we already know $q : X \ra X_l$ is a surjective geometric (hence categorical) $P$-quotient, and the remaining parts then follow from Proposition \ref{prop ss locus of perturbation of canonical linearisation}.

To prove the existence of the character $\chi$, we note that the quotienting group $\hH_i$ at each stage in the Quotienting-in-Stages procedure has 1-dimensional character group, arising form the 1-dimensional torus $\lambda_i(\GG_m)$ which descends from $\lambda^{[i]}(\GG_m) < P$. Let $\chi_i$ be the generator of the character group of $\hH_i$ which is dual to the 1-PS $\lambda_i$. Then at each stage we take a twist of the induced linearisation by small rational multiple $\epsilon_i \chi_i$, with the choice of $\epsilon_i$ depending on the choices of $\epsilon_j$ required for well-adaptedness for $j < i$. Recall from Definition \ref{defn row filtr groups} (1) that the 1-PS $\lambda$ for which $P = P(\lambda)$ corresponded to weights $\beta_1 > \dots > \beta_l$ with multiplicities $m_i$ and the 1-PS $\lambda^{[i]}$ has two weights $\beta_{\leq i} > \beta_{>i}$ with multiplicities $m_{\leq i}$ and $m_{<i}$.  We need to find a rational character $\chi_\epsilon$, corresponding to a vector $(r_1(\epsilon),\dots, r_l(\epsilon)) \in \QQ^l$ with $\sum_{i=1}^l m_i r_i(\epsilon) = 0$, such that for $1 \leq i \leq l-1$ we have
\[ \epsilon_i =: <\chi_i,\lambda_i>_{\mathfrak{t_i}} =  <\chi_\epsilon,\lambda^{[i]}>_{\mathfrak{t}} := \beta_{\leq i} \left(\sum_{j \leq i} m_jr_j(\epsilon) \right) + \beta_{>i}\left(\sum_{j > i} m_j r_j(\epsilon)\right) \]
where $\mathfrak{t} = \Lie (T)$ and $\mathfrak{t_i} = \Lie(T_i)$ for $T_i:=\lambda_i(\GG_m)$. If we set $\epsilon_0=\epsilon_l = 0$, this system of equations admits the solution
\[ r_i(\epsilon):= \frac{1}{m_i} \left( \frac{\epsilon_i}{\beta_{\leq i} - \beta_{>i}} - \frac{\epsilon_{i-1}}{\beta_{\leq i-1} - \beta_{>i-1}}\right) \text{ for } 1 \leq i \leq l,\]
where we note that the denominators appearing in this formula are non-zero.
\end{proof}

\begin{cor}\label{cor domain q in Xsmin}
If \hyperlink{D}{(D)} holds and $T$ acts trivially on $Z_{\min}$, then $\dom(q) \subset X_{\min}^s$.
\end{cor}
\begin{proof}
Since $T= Z(L) \cong \GG_m^{l-1}$ acts trivially on $Z_{\min}$, we have $Z^{\overline{L}-ss}_{\min} = Z^{R-ss}_{\min}=:Z^{ss}_{\min}$ and this in turn is equal to $Z^{s}_{\min}$ by  the third part of \hyperlink{D}{(D)}. Then the result follows from the last statement in Proposition \ref{prop character for Q-i-S} above as $X_{\min}^s:=p^{-1}(Z^{s}_{\min})$.
\end{proof}

In the next section, we will see that the Quotienting-in-Stages Assumption \hyperlink{QiS}{(QiS)} implies that $Z(L)$ acts trivially on $Z_{\min}$, allowing us to apply this Corollary.

\subsection{The Quotienting-in-Stages and Upstairs Stabiliser Assumptions}
\label{subsec qis and u assumptions}
In this section, we introduce the assumptions appearing in Theorem \ref{mainthm2}, which we arrange into two groups, called the \emph{Quotienting-in-Stages Assumption} \ref{quotienting in stages assumptions} \hyperlink{QiS}{(QiS)} and the \emph{Upstairs Stabiliser Assumption} \ref{Upstairs Stabiliser Assumptions} \hyperlink{U}{(U)} (or its weakened form \hyperlink{WUU}{(WUU)}). 

Assumption \hyperlink{QiS}{(QiS)} will be assumed to hold throughout, and is used to obtain an explicit description of the domain of definition of the Quotienting-in-Stages map. In contrast to this, the \emph{Upstairs Stabiliser Assumption} \hyperlink{U}{(U)} specifies ideal conditions on the stabilisers in $X$ under which we can use non-reductive GIT to obtain a projective geometric quotient. If \hyperlink{U}{(U)} fails, we can perform a sequence of blow-ups and, under the so-called \emph{Weak Upstairs Unipotent Stabiliser Assumption} \hyperlink{WUU}{(WUU)}, we obtain an explicit geometric and quasi-projective quotient of $X^{P-qs}$ (see Theorem \ref{prop blowup process if upstairs unip fails}, which completes the proof of Theorem \ref{mainthm2}).

\subsubsection{The Quotienting-in-Stages Assumption}

Here we use the notation from Definition \ref{defn important loci}.  

\begin{ass} \hypertarget{QiS}{(\emph{The Quotienting-in-Stages Assumption})}  \label{quotienting in stages assumptions} We say that:
	\begin{enumerate}
		\item\label{filt} \emph{The minimal weight spaces are filtered} if $Z^{(j)}_{\min} \subset Z^{(i)}_{\min} \subset  Z^{[i]}_{\min}$ for all $1 \leq i < j \leq l-1$.
		\item\label{split}  \emph{The preimages of $[i]$-minimal weight spaces are filtered} if for all $0 < i < j < l-1$  \[p_{[j+1]}^{-1}(Z^{[i]}_{\min} \cap Z^{[j+1]}_{\min} ) \subset p_{[j]}^{-1}(Z^{[i]}_{\min} \cap Z^{[j]}_{\min});\]  that is, for all $x \in X$, if $p_{[j+1]}(x) \in Z^{[i]}_{\min}\cap Z_{\min}^{[j+1]}$, then $p_{[j]}(x) \in Z^{[i]}_{\min}\cap Z_{\min}^{[j]}$.

	 		\item\label{non-deg} \emph{Non-degeneracy} holds: there exists $x \in X^{s}_{\min}$ with $p_{[j]}(x) \notin U^{[j-1]}Z^{[j-1]}_{\min}$ for all $ 2 \leq j \leq l$. 
	 		\end{enumerate}
When all these conditions hold, we say the \emph{Quotienting-in-Stages Assumption} \hyperlink{QiS}{(QiS)} holds.
\end{ass}

\begin{rmk} \label{rmk on qis assumptions} The first two assumptions are modelled on the behaviour of unstable HKKN strata for GIT problems associated to moduli of objects in an abelian category.  The final assumption, non-degeneracy, is equivalent to the non-emptiness of the Quotienting-in-Stages stable set when \hyperlink{QiS}{(QiS)} holds. If non-degeneracy fails, we should interpreted this as the Hilbert-Mumford criterion giving an empty semistable locus. For moduli of objects in an abelian category, non-degeneracy is a property of the relevant Harder-Narasimhan type which requires the existence of an object for which all extensions in the Harder--Narasimhan filtration are non-split. In the sheaf case, we interpret these assumptions in $\S$\ref{sec QiS Ass for sheaves}: we prove that the first two assumptions hold for unstable sheaves and give a sheaf-theoretic interpretation of non-degeneracy.
\end{rmk}

\begin{ex} \label{ex curves QiS doesn't always hold}
	It should be noted that \hyperlink{QiS}{(QiS)} is not always satisfied, even for the problem of performing quotients of unstable HKKN strata. Indeed, if we consider quintic plane curves $C\subset \PP^2$ up to the action of $\SL_2$, with the Killing form and the natural linearisation, we see that there is an HKKN stratum for which \[Z_\beta = \langle Y^5, X^4Y \rangle, \] and the associated Kempf one-parameter subgroup $\lambda_{\beta}$ has (up to scaling) weight vector $(2,1,-3)$  with respect to the natural coordinates. Thus we have $\lambda^{[1]}$ with weight vector $(2,-1,-1)$ up to scaling, and $\lambda^{[2]}$ has weight vector a scalar multiple of $(1,1,-2)$. Since the length of $\lambda_\beta$ is $3$, we have $\lambda^{(2)} = \lambda_{\beta}$. It is then easy to see that Assumption \eqref{filt} of \hyperlink{QiS}{(QiS)} is not satisfied, as
	\[Z^{[2]}_{\min} = \langle X^4Z \rangle \quad \text{and} \quad Z^{[1]}_{\min} = \langle Y^5 \rangle. \] 
\end{ex}

Before we proceed, let us note some straight-forward consequences of \hyperlink{QiS}{(QiS)}.

\begin{rmk}\label{rmk centre of L acts trivially on Zmin}
One immediate consequence of Assumption \eqref{filt}  in \hyperlink{QiS}{(QiS)} is that $T$ acts trivially on $Z_{\min}$. Indeed, this torus is generated by $\lambda^{[i]}(\GG_m)$ for $1 \leq i \leq l-1$, and we have $Z_{\min} \subset Z^{[i]}_{\min}$ for all such $i$, where each $Z^{[i]}_{\min}$ is fixed pointwise by $\lambda^{[i]}(\GG_m)$.
\end{rmk}

\begin{lemma} \label{lemma consequences of filt}
Suppose as in Assumption  \eqref{filt} of \hyperlink{QiS}{(QiS)} that $Z_{\min} \subset Z_{\min}^{[i]}$ for $1 \leq i \leq l-1$. Then we have:
	\begin{enumerate}[label=\emph{(\roman*)}]
		\item $X_{\min} \subset X_{\min}^{[i]}$ and, for $x \in X_{\min}$, we have $(p \circ p_{[i]})(x) = p(x)$,
		\item For $x \in X_{\min}$, we have $(p_{[i]} \circ p_{[j]})(x) = (p_{[j]} \circ p_{[i]})(x)$ for all $j \leq i$,
		\item $Z_{\min}^{(s)s} \subset Z_{\min}^{R^{(i)}-(s)s} \subset Z_{\min}^{[i],(s)s}$,
		\item $X_{\min}^{(s)s} \subset  X_{\min}^{[i],(s)s}$,
		\item If $x \in  X_{\min}^{(s)s}$, then $p_{[i]}(x) \in  X_{\min}^{(s)s}$.
	\end{enumerate}
\end{lemma}
\begin{proof}
	Note that the images of the 1-PSs $\lambda^{[i]}$ and $\lambda$ lie in the cental torus $T=Z(L)$. Using a sufficiently large power of our ample linearisation on $X$ we can embed $X$ into $\PP^N$ equivariantly and choose coordinates on $\PP^N$ so that the $T$-action is diagonalised. Then we can evaluate the weights of any 1-PSs of $T$ by pairing the $T$-weights with this 1-PS. By definition, $Z^{[i]}_{\min}$ is the minimal weight space for $\lambda^{[i]}$ (in $X$) and $X^{[i]}_{\min}$ is the locus of points in $X$ with at least one coordinate in $\PP^N$ having this minimal weight. Then $p_{[i]} : X^{[i]}_{\min} \ra Z^{[i]}_{\min}$ is the projection onto this minimal weight space. Assumption \eqref{filt} means all the $T$-weights that are minimal for $\lambda$ are minimal for $\lambda^{[i]}$ for any $ i$. Hence we can choose coordinates on $\PP^N$ such that $p$ is the projection on the first $n$ coordinates and $p_{[i]}$ is the projection onto the first $m$ coordinates with $n \leq m$. In particular, this shows (i). Similarly (ii) follows from Assumption \eqref{filt}.
	
For (iii), as $R^{(i)} < R$, we have $Z_{\min}^{(s)s} \subset Z_{\min}^{R^{(i)}-(s)s}$ and then the claim follows from (i) together with the functoriality of the semistable set with respect to equivariant closed immersions. For (iv), if $x \in X^{(s)s}_{\min}$, then $p(x) \in Z_{\min}^{(s)s} \subset Z_{\min}^{[i],(s)s}$ by (iii); therefore, by applying Lemma \ref{equiv pullback of semistable loci} to $p$, we conclude that $p_{[i]}(x)\in Z_{\min}^{[i],(s)s}$ as required.
	
Finally for (v), we have that $p(p_{[i]}(x)) =  p(x) \in Z^{(s)s}_{\min}$ using (i), and so $p_{[i]}(x) \in  X_{\min}^{(s)s}$.
\end{proof}

To complete the proof, we state the following result from about the behaviour of semistable loci in reductive GIT under equivariant pullbacks.

\begin{lemma}\label{equiv pullback of semistable loci} \cite[Chapter 1 \S5]{Mumford}
	Let $G$ be a reductive group acting linearly on polarised schemes $(X,\cL_X)$ and $(Y,\cL_Y)$ and let $f : X \ra Y$ is a $G$-equivariant morphism such that $f^*\cL_Y \cong \cL_X$. Then $f^{-1}(Y^{G-ss}) \subset X^{G-ss}$, and if $f$ is quasi-affine we also have $f^{-1}(Y^{G-s}) \subset X^{G-s}$.
\end{lemma}

\begin{rmk}\label{rmk consequences of filt}
 By an almost identical proof to Lemma \ref{lemma consequences of filt}, one can show from the Assumption \eqref{filt} of \hyperlink{QiS}{(QiS)} that $Z^{(i)}_{\min} \subset Z^{(j)}_{\min} \subset  Z^{[j]}_{\min}$ for $1 \leq j \leq i \leq l-1$ that the following statements hold.
	\begin{enumerate}[label=\emph{(\roman*)}]
		\item $X_{\min}^{(i)} \subset X_{\min}^{(j)}$  and for $x \in X_{\min}^{(i)}$, we have $p_{(i)} \circ p_{(j)}(x) = p_{(i)}(x)$,
		\item $X_{\min}^{(i)} \subset X_{\min}^{[j]}$ and for $x \in X_{\min}^{(i)}$, we have $p_{(i)} \circ p_{[j]}(x) = p_{(i)}(x)$ and $p_{[i]} \circ p_{[j]}(x) = p_{[j]} \circ p_{[i]}(x)$,
		\item $Z_{\min}^{(i),(s)s} \subset Z_{\min}^{(j),(s)s}$ and $Z_{\min}^{(i),(s)s} \subset Z_{\min}^{[j],(s)s}$,
		\item $X_{\min}^{(i),(s)s} \subset X_{\min}^{(j),(s)s}$ and $X_{\min}^{(i),(s)s} \subset X_{\min}^{[j],(s)s}$,
		\item If $x \in  X_{\min}^{(i),(s)s}$, then $p_{(j)}(x)$ and $p_{[j]}(x) \in  X_{\min}^{(i),(s)s}$.
\end{enumerate}
\end{rmk}

\begin{rmk}
One can generalise Lemma \ref{lemma consequences of filt} to higher weight spaces and associated attracting subschemes, which are higher Bia{\l}ynicki-Birula strata, which could be useful for constructing a quotient when non-degeneracy fails.
\end{rmk}

\subsubsection{An interpretation of (QiS) for sheaves of fixed HN type}\label{sec QiS Ass for sheaves}

Let us return to our running example of sheaves of fixed HN type $\tau = (P_1, \dots, P_l)$ on $(B,\cO_B(1))$. As in $\S$\ref{sec desc important loci sheaves}, we consider the action of $P_\tau$ on $X:=\overline{Y_\tau} \hookrightarrow Q_n$ linearised by $\cL_\tau$.

\begin{prop} \label{prop QiS 1 and 2 hold for sheaves}
Parts \eqref{filt} and \eqref{split} of \hyperlink{QiS}{(QiS)} hold for the linearised $P_\tau$-action on $X:=\overline{Y_\tau}$; that is the following statements hold.
\begin{enumerate}
\item For all  $0 < i  \leq j < l$, we have $Z^{(j)}_{\min} \subset Z^{(i)}_{\min} \subset Z^{[i]}_{\min}$.
\item For all $0 < i < j < l-1$, we have \[p_{[j+1]}^{-1}(Z^{[i]}_{\min} \cap Z^{[j+1]}_{\min} ) \subset p_{[j]}^{-1}(Z^{[i]}_{\min} \cap Z^{[j]}_{\min}).\]   
\end{enumerate}	
\end{prop}
\begin{proof}
The first claim follows immediately from the description of $Z_{\min} = Z_\tau$ and $Z^{[i]}_{\min}$ and $Z^{(i)}_{\min}$ given respectively in \eqref{Ztau} and Lemmas \ref{lem Z[i] for sheaves} and \ref{lemma Z(i) for sheaves} above. 

For the second, let $0 < i < j < l-1$. If $q_{\cE} \in X^{[j+1]}_{\min}$, then $p_{[j+1]}(q_{\cE}) = q_{\cE_{\leq j}} \oplus q_{\cE_{>j+1}}$ has corresponding quotient sheaf $\cF:=\cE^{(j+1)} \oplus \cE/\cE^{(j+1)}$ and similarly $p_{[j]}(q_{\cE})$ has corresponding quotient sheaf $\cG:=\cE^{(j)} \oplus \cE/\cE^{(j)}$. As at \eqref{filt E two}, the 1-PS $\lambda_{\tau}^{[i]}(t)$ induces filtrations
\[ \cF^{[\bullet, i]}= (0 \subset \cF^{(i)}=\cE^{(i)} \subset \cF ) \quad \text{and} \quad \cG^{[\bullet, i]}= (0 \subset \cG^{(i)}=\cE^{(i)} \subset \cG).\]
If $p_{[j+1]}(q_{\cE}) \in Z^{[i]}_{\min}$, then the inclusion $\cE^{(i)}\subset \cE^{(j+1)}$ is split (by $\phi: \cE^{(j+1)} \twoheadrightarrow \cE^{(i)}$), but then $\cF^{(i)} \subset \cE^{(j)}$ is also split (by $\phi_{|_{\cE^{(j)}}}$); that is, $p_{[j]}(q_{\cE}) \in Z^{[i]}_{\min}$.
\end{proof} 

Let us finally relate the non-degeneracy assumption \eqref{non-deg} to \lq non-degenerate' HN types. 

\begin{defn} \label{def HN type non-degen}
A HN type $\tau$ is \emph{non-degenerate} if there exists a sheaf $\cE$ of HN type $\tau$ such that the inclusion $\cE^{(i-1)} \hookrightarrow \cE^{(i)}$ in the HN filtration is non-split for each $i$ and the quotients $\cE_i := \cE^{(i-1)}/\cE^{(i)}$ are stable for all $i$. 
\end{defn}

The next result is an immediate corollary of Proposition \ref{prop first part of QiS stable for sheaves}.	
	
\begin{cor} \label{cor tau non-degen interp}
For a HN type $\tau$, Assumption \eqref{non-deg} of \hyperlink{QiS}{(QiS)} holds if and only if $\tau$ is non-degenerate. 
\end{cor}

\subsubsection{The Upstairs Stabiliser Assumption (U)}
\label{def star conditions stages} 
We can now give the stabiliser conditions \lq upstairs' on $X$ that will suffice to allow us to carry out Construction \ref{const qnt in stages}. 

\begin{ass}[\hypertarget{U}{\emph{The Upstairs Stabiliser Assumption}}]\label{Upstairs Stabiliser Assumptions}
Recall that $U^{[i-1,i]} := U^{[i-1]}\cap U^{[i]}$.
\begin{enumerate}
		\item We say that the \emph{upstairs unipotent stabiliser assumption} (UU) holds if for all $i$, 
\begin{equation}\label{starUstagesibullet}
		\dim \stab_{U^{[i]}}(-) \: \text{ and } \: \dim \stab_{U^{[i-1,i]}}(-) \: \text{ are constant on } X_{\min}^{[i],ss}. \tag*{(UU)}
		\end{equation}			
\item We say that the \emph{upstairs reductive stabiliser assumption} (UR) holds if for all $ 1 < i < l$,
			\begin{equation}\label{starR0stagesibullet}
			\dim \stab_{R^{(i)}}(z) = 0 \: \text{ for all } z \in Z_{\min}^{(i),ss}. \tag*{(UR)}
	\end{equation} 
 \end{enumerate}
 We say the \emph{Upstairs Stabiliser Assumption} \hyperlink{U}{(U)} holds if both the upstairs unipotent and reductive stabiliser assumptions hold.
\end{ass}

\brm The fact that we have to fix the dimension of stabilisers for both $U^{[i]}$ and $U^{[i-1,i]}$ is eventually explained by Proposition \ref{prop stab relns U}; in fact, if $X$ is connected, it is equivalent to fixing their difference by semi-continuity of the dimension of stabiliser groups. \erm 

We will eventually consider the case when \hyperlink{U}{(U)} fails and perform blow-up sequence for reductive and unipotent stabilisers. For the unipotent stabilisers we will use \S\ref{sec easy blowups abelian} based on \cite{BDHK2}. As observed in \S\ref{sec easy blowups abelian}, as well as \cite{BDHK} and \cite{Josh_length2}, we can obtain an explicit description of the open set we get a quotient of, provided the minimal weight space is not entirely blown up at any stage of the process. To ensure this, we introduce the following weaker assumption than \hyperlink{U}{(U)}.

\begin{ass}[\hypertarget{WUU}{\emph{The Weak Unipotent Upstairs Stabiliser Assumption}}]\label{WUU}
If $X^{P-qs} \neq \emptyset$ and there exists $z \in Z_{\min}$ having minimal dimensional  stabilisers for the groups $U^{[i]}$ and $U^{[i-1,i]}$
i.e.\
	\begin{equation}\label{cond not all of Zmin blownup} \tag*{[$Z_{\min} \not\subseteq C$]}
		\dim \stab_{U^{[i]}}(z) = d^{[i]}_{\min} \quad \text{and} \quad \dim \stab_{U^{[i-1,i]}}(z) = d^{[i-1,i]}_{\min} \quad \text{ for all } i,\end{equation}
we say that the \emph{weak upstairs unipotent stabiliser assumption} (WUU) holds.
\end{ass}

\brm \label{rmk following WUU can take the z in Zmins}
Note that if \hyperlink{WUU}{(WUU)} holds, then the set of $z \in Z_{\min}$ having minimal dimensional  stabilisers for $U^{[i]}$ and $U^{[i-1,i]}$ is non-empty and open. If $X$ is irreducible, so is $Z_{\min}$ and $Z^s_{\min} \subset Z_{\min}$ is a dense open, which intersects the open locus where the stabilisers for $U^{[i]}$ and $U^{[i-1,i]}$ are minimal, because we assumed $X^{P-qs} \neq \emptyset$.
\erm

\subsubsection{The Upstairs Stabiliser Assumption for Sheaves}\label{sec U for sheaves}

Let us return to our running example of sheaves of fixed HN type $\tau$. By Lemma \ref{lemma Z(i) for sheaves}, the Upstairs Reductive Assumption \ref{starR0stagesibullet} holds precisely when $\tau$ is a coprime HN type (semistability coincides with stability for each HN-subquotient). Unfortunately, while the quotienting-in-stages process tackles the problem with the central torus, the situation for the unipotent stabilisers is not so good. Often both \ref{starUstagesibullet} and \hyperlink{WUU}{(WUU)} fail to hold for sheaves. However, they are more readily verified for Higgs sheaves (and via the spectral correspondence Higgs sheaves correspond to certain sheaves on a higher dimensional base): one can show \hyperlink{WUU}{(WUU)} holds if one works with a Higgs HN type whose underlying sheaf HN type is non-empty \cite{HHJ}.

Unipotent stabilisers of sheaves (and more general moduli problems in a linear abelian category) are better behaved with respect to the quotienting-in-stages framework as follows.

\blm \label{lem sheaves p preserves stabs}
For all $q_{\cE} \in X_{\min}^{ss} = Y_\tau^{ss}$ and all $i$, the dimensions of the $U_\tau^{[i]}$-stabilisers (resp. $U_\tau^{[i-1,i]}$-stabilisers) of $q_{\cE}$ and $p_{[i]}(q_{\cE})$ coincide.
\elm 
\bpf
This follows from Proposition \ref{prop QiS unipotent stab for sheaves}.
\epf 

\begin{rmk} \label{rmk QiS stable set is open}
If \hyperlink{QiS}{(QiS)} and \hyperlink{UU}{(UU)} hold, we claim that the Quotienting-in-Stages stable locus $X^{P-qs}$ is open and non-empty. Indeed, as $\dim \stab_{U^{[i]}}(-)$  is constant on $X_{\min}^{[i],ss}$, then $U^{[i]}Z^{[i],(s)s}_{\min} \subset X_{\min}^{[i],(s)s}$ is a closed subset by \cite[Lemma 5.4]{BDHK2}. Note that $X^{(s)s}_{\min} \subset X_{\min}^{[i],(s)s}$ by Lemma \ref{lemma consequences of filt} (v). Hence, $X^{P-qs}$ is open in $X^{s}_{\min}$, which itself is open in $X_{\min}$. Thus provided $X_{\min}$ is open in $X$, which for example will be the case if $X$ is irreducible\footnote{See $\S$\ref{sec remove irred assumption} for how to remove this irreducibility assumption.}, the Quotienting-in-Stages stable locus $X^{P-qs}$ is open in $X$. Finally, non-degeneracy \eqref{non-deg} of \hyperlink{QiS}{(QiS)} implies this set is non-empty. 
\end{rmk}

\subsection{Roadmap of the proof of Theorem \ref{mainthm2}}\label{sec roadmap}

In the next section, we will prove our main result, Theorem \ref{mainthm2}, on this quotienting-in-stages procedure. Let us first outline the proof.

In Corollary \ref{cor exists open set with qnt}, we saw that if \hyperlink{D}{(D)} holds then we obtain a projective geometric $P$-quotient of some open subset $\dom(q) \subset X$. In \S\ref{sec describing domain q with downstairs stab}, we first show under assumption \hyperlink{D}{(D)} that the domain of definition of this quotient map $q$ has an explicit description (see Theorem \ref{thm with downstairs assumptions}). 
This proof is a form of induction on the length $l$, making use of the auxiliary conditions $D(i)$, $E(i)$, $F(i)$, which are defined below in Definition \ref{defn DEF}. The basic idea is that, as a consequence of \hyperlink{QiS}{(QiS)}, the stage $i$ minimal weight space $Z_{i,\min}$ downstairs in $X_i$ is the image under $q_{(i-1)}$ of the $[i]$-minimal weight space upstairs (see Lemma \ref{lem can take limits r}). After this, the most subtle point in the proof of Theorem \ref{thm with downstairs assumptions} lies in \lq pushing forward' reductive GIT stability of points in minimal weight spaces along the quotient maps $q_{(i)}$ (see Proposition \ref{prop claim 4}). 

The next step is to replace \hyperlink{D}{(D)} with \hyperlink{U}{(U)} and in this case show that the domain of $q$ coincides with the Quotienting-in-Stages stable locus $X^{P-qs}$ introduced in Definition \ref{defn quotienting-in-stages stable locus}. For this, we prove some comparison results for upstairs versus downstairs unipotent stabiliser dimensions in \S\ref{subsec two stab comparison results}, which allow us to deduce Proposition \ref{prop no blow ups get proj qnt} that gives the part of Theorem \ref{mainthm2} under the assumptions \hyperlink{QiS}{(QiS)} and \hyperlink{U}{(U)}: namely $X^{P-qs}= \dom (q)$ is a non-empty open, which admits a projective geometric $P$-quotient.

In \S\ref{sec upstairs red fails}, we turn our attention to what happens when \hyperlink{QiS}{(QiS)} and  \ref{starUstagesibullet} hold, but \ref{starR0stagesibullet} fails. In this case we much perform partial desingularisations as in \cite{K2} for the reductive groups $R^{(i)}$ until we achieve \ref{starR0stagesibullet} on the blown up space $\widehat{X}$, whilst also preserving \hyperlink{QiS}{(QiS)} and \ref{starUstagesibullet} (see Proposition \ref{prop blowup process upstairs red fails}). Then, in Proposition \ref{prop qnt in stages no blow ups} we simply apply Proposition \ref{prop no blow ups get proj qnt} to $\widehat{X}$ to deal with this case.  Finally in $\S$\ref{sec quot stages with blowups}  we deal with the case when \ref{starUstagesibullet} fails by performing a blow-up sequence concerning unipotent stabilisers as described in $\S$\ref{sec easy blowups abelian}, which is based on an earlier arXiv version of \cite{BDHK2}. To achieve our result, we must avoid the bad situation where the whole minimal weight space is blown up at any stage, and this is where we use the assumption \hyperlink{WUU}{(WUU)} to complete the proof of Theorem \ref{mainthm2}.

\section{Quotienting-in-stages for parabolics: Proof of the Main Theorem}
\label{sec proof}
In this section, we assume that we have a linearised action of a parabolic group $P<\SL_n$ on an irreducible projective scheme $X$ as in $\S$\ref{sec important notn}. We now proceed with the proof of Theorem \ref{mainthm2}, as outlined in $\S$\ref{sec roadmap} above. Throughout this section, we  assume that the Quotienting-in-Stages Assumption \hyperlink{QiS}{(QiS)} holds.

\subsection{Describing the stable locus under the Downstairs Stabiliser Assumption (D)}\label{sec describing domain q with downstairs stab}

As a first step in our proof of Theorem \ref{mainthm2}, we will obtain the following result.

\begin{thm} \label{thm with downstairs assumptions}
Suppose that \hyperlink{QiS}{(QiS)} and \hyperlink{D}{(D)} hold. Then we have 
\[ \dom(q^s) = \dom(q)= \left\{x\in X^{s}_{\min}  \:\middle\vert\:p_{[j]}(x) \notin U^{[j-1]}Z^{[j-1]}_{\min} \text{ for all } 2\leq  j\leq l\right\}\]
and this locus has a projective geometric $P$-quotient.
\end{thm}

Recall that by Corollary \ref{cor domain q in Xsmin} we know that when \hyperlink{QiS}{(QiS)} and \hyperlink{D}{(D)} hold we have $\dom (q) \subset X^s_{\min}$, so that $q = q^s$ and to prove the above theorem it suffices to prove $\dom (q^s)$ is the explicit set described above.

Our proof of Theorem \ref{thm with downstairs assumptions} is inductive, making use of the following conditions. 

\begin{defn} \label{defn DEF} (Description of $D(i), E(i)$ and $F(i)$) Let $ 1 \leq  i \leq l-1$.
	\begin{enumerate}
		\item We say that the \emph{$i$th Description holds} (denoted $D(i)$) if
		\[\dom (q^s_{(i)}) \cap Z^{[i+1]}_{\min} \neq \emptyset.\] 
		We write $D(\leq k)$ if the above holds for all $i\leq k$.
		\item We say that the \emph{$i$th Explicit description holds} (denoted $E(i)$) if
		\[ \dom(q^s_{(i)})= \cE(i):= \left\{x\in X^{s}_{\min} \setminus U^{[i]}Z^{[i],s}_{\min} \:\middle\vert\: p_{[j]}(x) \in \dom (q_{(j-1)}) \text{ for all } 2\leq j\leq i\right\} \] We say $E(\leq k)$ holds if $E(i)$ holds for all $i\leq k$.
		\item We say that the \emph{$i$th Final description holds} (denoted $F(i)$) if
		\[ \dom(q^s_{(i)}) = \cF(i):= \left\{x\in X^{s}_{\min} \setminus U^{[i]}Z^{[i],s}_{\min} \:\middle\vert\:p_{[j]}(x) \notin U^{[j-1]}Z^{[j-1]}_{\min} \text{ for all } 2\leq  j\leq i\right\}. \]We say $F(\leq k)$ holds if $F(i)$ holds for all $i\leq k$.
	\end{enumerate}
\end{defn}

 In particular, for $i =l -1$, we get
\[ \cF(l-1) =\left\{x\in X^{s}_{\min}  \:\middle\vert\:p_{[j]}(x) \notin U^{[j-1]}Z^{[j-1]}_{\min} \text{ for all } 2\leq  j\leq l\right\}  \] 
since $\lambda^{[l]}$ is trivial, thus $p_{[l]} =\mathrm{Id}_X$ (see Remark \ref{rmk square and round minimal weight spaces}). Thus in Theorem \ref{thm with downstairs assumptions}, we will show that $\dom(q^s) = \cF(l-1)$. In fact, under the upstairs stabiliser assumptions \hyperlink{U}{(U)}, one has $\cF(l-1)= X^{P-qs}$.

\brm \label{rmk another defn of F E}
In addition to the above (set-theoretic) descriptions, we can characterise the sets $\cF(i)$ and $\cE(i)$ scheme-theoretically using the equalities:
	\[  \cE(i):=\left( X^{s}_{\min} \setminus U^{[i]}Z^{[i],s}_{\min}\right) \cap \left(\bigcap_{2 \leq j \leq i} p_{[j]}^{-1}(\dom(q_{(j-1)}) \cap Z^{[j]}_{\min} ) \right)\]
\[ \cF(i):=X^{s}_{\min} \setminus \left( U^{[i]}Z^{[i],s}_{\min} \cup \bigcup_{2 \leq j \leq i} p_{[j]}^{-1} (U^{[j-1]}Z^{[j-1]}_{\min} \cap Z^{[j]}_{\min}) \right). \]
\erm

We now indicate the strategy of our proof.  We will inductively show that $F(i)$ holds,  by proving the following implications
\begin{equation}\label{s chain of implications}
	\begin{array}{l}
	D(\leq i-1) \\ E(\leq i-1) \\ F(\leq i-1) 
	\end{array} \implies \begin{array}{l}
D(\leq i-1) \\ E(\leq i)\\ F(\leq i-1) 
\end{array}  \implies \begin{array}{l}
D(\leq i-1) \\ E(\leq i)\\ F(\leq i) 
\end{array} \implies \begin{array}{l}
D(\leq i) \\ E(\leq i)\\ F(\leq i) 
\end{array} 
\end{equation}

For $i =1$, under the assumption $[\hU_1;R_1\text{-ss}]$ on $X_1 = X$, we have by definition of $q^s_{(1)}$ that
\[ \dom(q^s_{(1)}) = X^{s}_{\min} \setminus  U^{[1]}Z^{[1],s}_{\min} \]
and, in particular, $E(1)$ and $F(1)$ hold, which form the base step of our inductive proof. 

Before we begin the induction, let us note some immediate consequences of condition $D(i)$.

\subsubsection{Consequences of $D(i)$}

\begin{lemma} \label{lem can take limits r}
If $D(i)$ holds, then the following hold for the map $q_{(i)} : X \dashrightarrow X_{i+1}$.
\begin{enumerate}[label=\emph{(\arabic*)}]
\item For $x \in q_{(i)}^{-1}(Z_{i+1,\min})$, we have $p_{[i+1]}(x) \in \dom (q_{(i)})$,
\item For $x \in \dom (q_{(i)})$, we have $q_{(i)}(x) \in X_{i+1,\min}$ if and only if $x \in X^{[i+1]}_{\min}$ and also $p_{[i+1]}(x) \in \dom (q_{(i)})$,

\item We have $q_{(i)}^{-1}(Z_{i+1,\min}) = (U^{[i]} Z^{[i+1]}_{\min}) \cap \dom (q_{(i)})$,

\item In particular, for all $x\in \dom (q_{(i+1)})$, we have $p_{[i+1]}(x) \in \dom(q_{(i)})$.
\end{enumerate}
\end{lemma}
\begin{proof} 
As a conequence of $D(i)$, we have \begin{equation}\dom (q_{(i)}) \cap Z^{[i+1]}_{\min} \neq \emptyset. \label{eqn d(i)} \end{equation}

	The map $q_{(i)}$ is a $\hH^{(i)}$-quotient which is constructed by taking invariants and so looks like a projection $\PP(V^*) \dashrightarrow \PP(W^*)$ where $V = H^0(X, \cL^{\otimes N}) \supset W:=V^{\hH^{(i)}}$ and $N$ is sufficiently large. Then Equation \eqref{eqn d(i)} implies that the maximal weights for $\lambda^{[i+1]}(\GG_m)$ acting on both $V$ and $W$ coincide. Hence, we can choose a basis of $V$ of the form 
	\[V= <\alpha_1, \dots , \alpha_n, \alpha_1',\dots ,\alpha'_{n'},\beta_1, \dots \beta_m, \beta'_1,\dots ,\beta'_{m'} >\]
	where $W = < \alpha_1, \dots , \alpha_n, \alpha_1',\dots ,\alpha'_{n'}>$ and the maximal weight space for $\lambda^{[i+1]}(\GG_m)$ in $V$ is spanned by $\alpha_1, \dots , \alpha_n, \beta_1, \dots \beta_m$. With respect to this basis, the projection $\PP(V^*) \dashrightarrow \PP(W^*)$ corresponding to $q_{(i)}$ is the projection onto the first $n+n'$ coordinates (corresponding to the elements $\alpha_j$ and $\alpha'_j$). Moreover, the projection $p_{[i+1]} : X \dashrightarrow Z^{[i+1]}_{\min}$ to the minimal weight space is the projection onto the coordinates corresponding to $\alpha_j$ and $\beta_j$. Thus a point lies in $X^{[i+1]}_{\min}$ precisely when one of these coordinates is non-zero. 
	
	For (1), we use the above descriptions of $q_{(i)}$ and $p_{[i+1]}$. Let $x \in \dom (q_{(i)})$ with $q_{(i)}(x) \in Z_{i+1,\min}$; then all the $\alpha'_j$ coordinates of $x$ are zero. However, as $x \in \dom(q_{(i)})$, there must be one $\alpha_j$ coordinate of $x$ that is non-zero, which also means that $p_{[i+1]}(x)$ has a non-zero $\alpha_j$ coordinate and thus is in the domain of definition of $q_{(i)}$.
	
	For (2), we also use these descriptions of $q_{(i)}$ and $p_{[i+1]}$. Since $x \in \dom(q_{(i)})$, one of the $\alpha_j$ or $\alpha'_j$ coordinates of $x$ is non-zero. Then $p_{i+1}(q_{(i)}(x))$ is the projection onto the $\alpha_j$ coordinates of $x$ and so $q_{(i)}(x) \in X_{i+1,\min}$ if and only if one of the  $\alpha_j$ coordinates of $x$ is non-zero. Similarly $x \in X^{[i+1]}_{\min}$ and $p_{[i+1]}(x) \in \dom(q_{(i)})$ is equivalent to one of the $\alpha_j$ coordinates of $x$ being non-zero.
	
	For (3), let us first show the inclusion of the right side in the left. If $x \in \dom(q_{(i)}) \cap U^{[i]}Z^{[i+1]}_{\min}$, then we claim that $q_{(i)}(x) \in Z_{i+1,\min}$. Since $q_{(i)}$ is $U^{(i)}$-invariant, we can assume $x \in Z^{[i+1]}_{\min}$; then its image $z:=q_{(i)}(x) \in X_{i+1}$ is fixed by the 1-PS $\lambda_{i+1}$ induced by the $\lambda^{[i+1]}$-action on $X$. Furthermore, as the induced linearisation $\cL_{i+1} \ra X_{i+1}$ satisfies $q_{(i)}^*\cL_{i+1} \cong \cL^{\otimes N}$, the $\lambda_{i+1}$-weight of $\cL_{i+1}$ over $z$ equals the $\lambda^{[i+1]}$-weight of $\cL^{\otimes N}$ over $x$, as both of these points are fixed by these 1-PSs and $z=q_{(i)}(x)$. By Equation \eqref{eqn d(i)}, the minimal weight for $\lambda_{i+1}$ acting on $\cL_{i+1} \ra X_{i+1}$ is the minimal weight as for $\lambda^{[i+1]}$ acting on $\cL^{\otimes N} \ra X$. Thus the $\lambda^{[i+1]}$-fixed point $x$ has minimal $\lambda^{[i+1]}$-weight if and only if the $\lambda_i$-fixed point $z$ has minimal $\lambda_{i+1}$-weight. Therefore, $z \in Z_{i+1,\min}$ as required. For the inclusion of the left side in the right, suppose that $x \in \dom(q_{(i)})$ and $z=q_{(i)}(x) \in Z_{i+1,\min}$. By (1), we know that $\overline{x}:=p_{[i+1]}(x) \in \dom(q_{(i)})$ and also $q_{(i)}(\overline{x}) = p_{i+1}(q_{(i)}(x))= z$. Arguing as above, the $\lambda^{[i+1]}$-weight of $\overline{x}$ agrees with the $\lambda_{i+1}$-weight of $z$, which is minimal and thus $\overline{x} \in Z^{[i+1]}_{\min}$. Since $q_{(i)}(\overline{x}) = z = q_{(i)}(x)$ and $q_{(i)}$ is a geometric $\hH^{(i)}$-quotient, $x$ lies in the $\hH^{(i)}$-orbit of $\overline{x} \in Z^{[i+1]}_{\min}$. Since $\lambda^{[i+1]}$ is central in the Levi $L$, the action of the Levi preserves $Z^{[i+1]}_{\min}$, thus $x \in \hH^{(i)}Z^{[i+1]}_{\min}=U^{[i]}Z^{[i+1]}_{\min}$ as required.	
	
For (4), since $x \in \dom(q_{(i+1)}) = q_{(i)}^{-1}(\dom(q_{i+1}))$, we have $q_{(i)}(x) \in \dom(q_{i+1}) \subset X_{i+1,\min}$ and so we conclude that $p_{[i+1]}(x) \in \dom(q_{(i)})$ using (2).
\end{proof}

Lemma \ref{lem can take limits r} tells us that, provided $D(i)$ holds, the image of $Z^{[i+1]}_{\min}$ under the quotient map $q_{(i)}$ will be the new minimal weight space $Z_{i+1\min}$. It is then natural to suppose that the same might be true for the relevant stable loci. This is a more subtle matter, but we shall confirm this is the case in Proposition \ref{prop claim 4} below; first we need a technical preliminary lemma.

\begin{lemma} \label{lem an induction for claim 4} 
Suppose that \hyperlink{Qis}{(QiS)} and \hyperlink{D}{(D)} hold. Fix $1 < i < l$ and suppose that $D(\leq i-1)$ holds. For $1<j\leq i$ and $x \in \dom (q_{(j-1)}) \cap X^{(i),s}_{\min}$ with $p_{[j]}(x) \in \dom (q_{(j-1)})$, the image $q_{(j-1)}(p_{[j]}(x)) \in Z_{j,\min}$ is stable for the action of $R^{(j,i)}:=R^{(i)}/R^{(j-1)}$. 
\end{lemma} 
\begin{proof} 
We prove the claim by induction on $j$. For $j=2$, suppose $x \in \dom (q_{(1)}) \cap X^{(i),s}_{\min}$ with $\overline{x}:=p_{[2]}(x) \in \dom(q_{(1)})$; then we need to show that $q_{(1)}(\overline{x}) \in Z_{2,\min}$ is stable for $R^{(2,i)}$. Consider the action of the larger non-reductive group
\[ \hH':=U^{(1)}\rtimes (T^{(1)} \times R^{(i)}) > \hH^{(1)} =U^{(1)}\rtimes (T^{(1)} \times R^{(1)}) \]
on $X=X_1$, where the unipotent radical $U^{(1)}=U^{[1]}$ is graded by $\lambda^{[1]}(\GG_m) = T^{(1)}$ with minimal weight space $Z^{[1]}_{\min}$. Recall that the $\hH^{(1)}$-quotient of $q_{(1)} : X_1 \dashrightarrow X_2$ is constructed by first taking a $U^{(1)}$-quotient of $X^{[1]}_{\min}$ and then taking a residual reductive GIT quotient of $T^{(1)} \times R^{(1)}$ acting on a projective completion $Y:= \overline{X^{[1],ss}_{\min}/U^{[1]}}$  described in Remark \ref{rmk Uhat part of proj completion constr}. Since \ref{hUbullet,Rss} holds and the $R^{(i)}$-semistable locus $Z^{[1],R^{(i)}-ss}_{\min}$ is contained in the $R^{(1)}$-semistable locus $Z^{[1],ss}_{\min}$, we can take a NRGIT quotient of the $\hH'$-action on $X=X_1$; this is the reductiveGIT  quotient of $X_2$ by $R^{(2,i)}$ or equivalently the reductive GIT quotient of $Y$ by $T^{(1)} \times R^{(i)}$. Hence, we have the following commutative diagram of quotients 

\[ \xymatrixcolsep{8pc} \xymatrix{ 
X \ar@{-->}[rd]_{q_{(1)} = q_{\hH^{(1)}}} \ar@{-->}[r]^{\quad q_{U^{(1)}}} \ar@{-->}@/^2pc/[rr]^{q_{\hH'}} &  Y \ar@{-->}[d]^{q_{T^{(1)} \times R^{(1)}}} \ar@{-->}[r]^{q_{T^{(1)} \times R^{(i)}}\quad\quad\quad \quad} & Y/\!/(T^{(1)} \times R^{(i)}) \\
& X_2 = Y/\!/(T^{(1)}\times R^{(1)}) \ar@{-->}[ru]_{q_{R^{(2,i)}}} 
}
\]
where the subscript denotes the group by which we have quotiented. To prove that $q_{(1)}(\overline{x})$ is $R^{(2,i)}$-stable, we want to apply Lemma \ref{lemma on GIT stability images red} below to the action of the action of the product of $G_1 := T^{(1)} \times R^{(1)}$ and $G_2:= R^{(2,i)}$ on $Y$. For this, we first need to verify that $q_{U^{(1)}}(\overline{x}) \in Y$ is stable for $G_1 \times G_2$; by Theorem \ref{thm nrGIT HM type inclusions} \emph{\ref{nrGIT HM inclusion 1}}, it suffices to check that $p_{[1]}(\overline{x}) \in Z^{[1],R^{(i)}-s}_{\min}$, as we already know that $\overline{x} =p_{[2]}(x) \in  \dom(q_{(1)}) \subset X^{[1]}_{\min} \setminus U^{(1)} Z^{[1]}_{\min}$. By Assumption \eqref{filt} of \hyperlink{Qis}{(QiS)}, we have that $Z^{(i)}_{\min} \subset Z_{\min}^{[1]}$ and in fact this is a closed $R^{(i)}$-invariant subscheme; therefore, $Z^{R^{(i)}-s}_{\min} = Z^{(i)}_{\min} \cap Z_{\min}^{[1], R^{(i)}-s}$. Since $x \in X^{(i),s}_{\min}$ we have also $\overline{x} =p_{[2]}(x) \in X^{(i),s}_{\min}$ by Remark \ref{rmk consequences of filt} (v); thus $p_{(i)}(\overline{x}) \in Z^{s}_{\min} \subset Z_{\min}^{R^{(i)}-s}$. Since $p_{(i)} \circ p_{[1]} = p_{(i)}$ on $X^{(i)}_{\min}$ by Remark \ref{rmk consequences of filt} (ii), we deduce that $p_{[1]}(\overline{x}) \in  Z^{[1],R^{(i)}-s}_{\min}$, which completes the proof of the base case $j = 2$.

For the inductive step, fix $k < i$ and suppose the claimed statement holds for all $j \leq k$. To prove the statement for $k+1$, we take $x \in \dom (q_{(k)}) \cap X^{(i),s}_{\min}$ such that $\overline{x}:=p_{[k+1]}(x) \in \dom (q_{(k)})$ and want to show that $q_{(k)}(\overline{x}) \in Z_{k+1,\min}$ is stable for the action of $R^{(k+1,i)}$. We have $q_{(k)} = q_k\circ q_{(k-1)}$, where $q_k : X_k \dashrightarrow X_{k+1}$ is the $k$th stage quotient. Let $y := q_{(k-1)}(\overline{x})$; then $y \in \dom(q_k) \subset X_{k,\min} \setminus U_k Z_{k,\min}$, as $ \overline{x}\in \dom (q_{(k)}) = q_{(k-1)}^{-1}(\dom(q_k))$. We want to apply the inductive hypothesis to $\overline{x} \in \dom(q_{(k-1)})$; for this, note that $x \in X^{(i),s}_{\min}$ implies also $\overline{x} \in X^{(i),s}_{\min}$ by Remark \ref{rmk consequences of filt} (v), and also $p_{[k]}(\overline{x}) \in \dom(q_{(k-1)})$, as $y := q_{(k-1)}(\overline{x}) \in X_{k,\min}$ and so we can apply Lemma \ref{lem can take limits r} (2) since $D(k-1)$ holds. By the inductive hypothesis applied to $\overline{x} \in \dom(q_{(k-1)}) \cap X^{(i),s}_{\min}$, we know that $q_{(k-1)}(p_{[k]}(\overline{x})) \in Z_{k,\min}$ is $R^{(k,i)}$-stable. Therefore
\[ p_{k}(y) = p_k(q_{(k-1)}(\overline{x})) = q_{(k-1)}(p_{[k]}(\overline{x}))\]
is $R^{(k,i)}$-stable. It remains to show that $q_{k}(y) = q_{(k)}(\overline{x}) \in Z_{k+1,\min}$ is $R^{(k+1,i)}$-stable. For this, let $G_1:=\lambda_k(\GG_m) \times R_k$ and $G_2 := R^{(k+1,i)}$; then, as in the base case, we consider the following composition of quotients 
\[ \xymatrixcolsep{4pc} \xymatrix{ 
X_k \ar@{-->}[r]^{q_{U_k} \quad \quad } & Y_k:=\overline{X^{ss}_{k,\min}/U_k} \ar@{-->}[r]^{\quad \quad q_{G_1}} & X_{k+1} \ar@{-->}[r]^{q_{G_2} \quad \quad} & X_{k+1} /\!/ G_2 }\] 
where the subscript on each quotient denotes the quotienting group and $q_k$ is the composition of the first two of these rational maps. The resulting quotient on the right is the NRGIT quotient of $X_k$ by $U_k \rtimes (G_1 \times G_2)$ (using $\lambda_k$ to grade $U_k$) and also the reductive GIT quotient of $Y_k$ by $G_1 \times G_2$. To deduce that $q_k(y) = q_{G_1}(q_{U_k}(y))$ is stable with respect to $G_2=R^{(k+1,i)}$, we want to apply Lemma \ref{lemma on GIT stability images red} to the action of $G_1 \times G_2$ on $q_{U_k}(y) \in Y_k$; for this, it remains to check that $q_{U_k}(y)$ is stable for the product $G_1 \times G_2$. Since we have shown $y \in X_{k,\min} \setminus U_k Z_{k,\min}$ and $p_{k}(y)$ is $R^{(k,i)}$-stable, it follows that $q_{U_k}(y)$ is $(G_1 \times G_2)$-stable by Theorem \ref{thm nrGIT HM type inclusions} \emph{\ref{nrGIT HM inclusion 1}} applied to the NRGIT quotient of $X_k$ by $U_k \rtimes (T_k \times R^{(k,i)})$.
\end{proof}

To complete the proof of this lemma, we need to prove the following lemma on performing reductive GIT quotients in stages.

\begin{lemma}\label{lemma on GIT stability images red}
Let $G = G_1\times G_2$ be a product of reductive groups acting linearly on a projective scheme $X$ with reductive GIT quotient $q_G: X \dashrightarrow X \git G$. If $x \in X$ is $G$-stable, then $q_{G_1}(x)$ is $G_2$-stable, where $q_{G_1} : X\dashrightarrow X\git G_1$ is the $G_1$-quotient. 
\end{lemma}
\begin{proof}
Let $y:=q_{G_1}(x) \in Y:=X\git G_1$ and write $q_{G_2} : Y \dashrightarrow Y \git G_2$; then $y \in Y^{G_2-ss} = \dom(q_{G_2})$ as $x \in X^{G-ss} = \dom(q_G)$ and $q_G = q_{G_2} \circ q_{G_1}$. To prove that $y$ is $G_2$-stable, it suffices to show that its orbit is closed in the semistable locus and its stabiliser is zero dimensional.

Suppose that there is an orbit $G_2 \cdot y'$ of $y' = q_{G_1}(x') \in Y^{G_2-ss}$ contained in the closure of $G_2 \cdot y$. Then $q_{G_2}(y) = q_{G_2}(y')$ and also $q_G(x) = q_{G}(x')$. Since $x$ is $G$-stable, $x$ and $x'$ lie in the same $G$-orbit and thus also $y$ and $y'$ lie in the same $G_2$-orbit. Hence $G_2 \cdot y$ is closed in $Y^{G_2-ss}$.

It remains to show that $\stab_{G_2}(y)$ is zero dimensional. If $\dim\stab_{G_2} (y) >0$, then for each $s \in \stab_{G_2}(y)$ there is $g_s \in G_1$ such that $g_ssx=x$, as the $G$-stable point $x$ is also $G_1$-stable  (for example, this follows by the Hilbert-Mumford criterion) and so $q_{G_1}^{-1}(y) = G_1 \cdot x$. The projection of the set \[S :=\{g_s\cdot s  \mid s \in \stab_{G_2}(y)\} \subseteq \stab_{G_1\times G_2}(x) \leq  G_1\times G_2\] to $G_2$ is $\stab_{G_2}(y)$; thus $S$ also has strictly positive dimension, which contradicts $x \in X^{G-s}$.
\end{proof}

We can now show, when $D(\leq i-1)$ holds, that the image of stable $[i]$-minimal weight space $Z^{[i],s}_{\min}$ under the quotient map $q_{(i-1)}$ is the stage $i$ stable minimal weight space $Z_{i,\min}^s$.

\begin{prop}\label{prop claim 4}
Assume that \hyperlink{QiS}{\hyperlink{QiS}{(QiS)}} and \hyperlink{D}{(D)} hold. If $D(\leq i-1)$ holds, then for all $x\in \dom (q_{(i-1)}) \cap Z^{[i],s}_{\min}\cap X^{(i),s}_{\min}$ we have $q_{(i-1)}(x) \in Z_{i,\min}^s$.
\end{prop}
\begin{proof}
Our goal is to use Lemma \ref{lem an induction for claim 4} with $j = i$. This applies as for $x\in \dom (q_{(i-1)}) \cap Z^{[i],s}_{\min}\cap X^{(i),s}_{\min}$, we have $p_{[i]}(x) = x \in \dom (q_{(i-1)})$. We conclude that $q_{(i-1)}(p_{[i]}(x))$ is stable for the action of $R_i$. Since $p_{[i]}(x) = x$, this completes the proof that $q_{(i-1)}(x) \in Z_{i,\min}^{s}$. 
\end{proof}

\subsubsection{The induction} We now begin our induction. Throughout this section, we assume that \hyperlink{QiS}{(QiS)} and \hyperlink{D}{(D)} hold.

\begin{lemma}\label{lem a^s holds r}
If $F(\leq i)$ holds, then $D(i)$ holds. 
\end{lemma} 
\bpf 
By the non-degeneracy assumption \eqref{non-deg} of \hyperlink{QiS}{(QiS)}, there exists $x \in X^{s}_{\min}$ with $p_{[j]}(x) \notin U^{[j-1]}Z^{[j-1]}_{\min}$ for all $ 1 < j < l$. Since $X^{s}_{\min} \subset X^{[i+1],s}_{\min}$ by Lemma \ref{lemma consequences of filt} (iv), we have $y:=p_{[i+1]}(x) \in Z^{[i+1]}_{\min}$. Therefore, it suffices to show that also $y \in \cF(i)$, since we have $\cF(i) = \dom(q^s_{(i)})$ as $F(i)$ holds. By assumption, we have $y = p_{[i+1]}(x) \notin U^{[i]}Z^{[i]}_{\min}$. As $x \in X^{s}_{\min}$, we have by $y=p_{[i+1]}(x) \in X^{s}_{\min}$ by Lemma \ref{lemma consequences of filt} (v). Thus it remains to check that $p_{[j]}(y) \notin U^{[j-1]}Z^{[j-1]}_{\min}$ for all $j \leq i$. However, if $p_{[j]}(y) \in U^{[j-1]}Z^{[j-1]}_{\min}$, then as $p_{[j]} \circ p_{[i+1]} = p_{[i+1]}\circ p_{[j]}$ on $X_{\min}$ by Remark \ref{rmk consequences of filt} (ii), we would obtain by Assumption \eqref{split} of \hyperlink{Qis}{(QiS)} that $p_{[j]}(x) \in U^{[j-1]}Z^{[j-1]}_{\min}$, which contradicts our assumption on $x$. Therefore, $p_{[j]}(y) \notin U^{[j-1]}Z^{[j-1]}_{\min}$ and we have shown $y \in \cF(i)$ as required.
\epf

We next prove the middle implication in \eqref{s chain of implications}.

	\begin{prop}\label{prop middle implication r}
 If $E(\leq i)$ and $F(\leq i-1)$ hold, then $F(i)$ also holds; that is,
	\begin{equation}\label{inclusion for F(i) r}
		\dom(q^s_{(i)})  = \cF(i):=\left\{x\in X^{s}_{\min} \setminus U^{[i]}Z^{[i],s}_{\min} \:\middle\vert\:p_{[j]}(x) \notin U^{[j-1]}Z^{[j-1]}_{\min} \text{ for all } 2\leq  j\leq i\right\}. 
	\end{equation}
\end{prop}
\begin{proof}
	Suppose that $x \in \dom(q^s_{(i)})$. Since $E(i)$ holds, we have $x \in \cE(i)$ and in particular, $x \in X^{s}_{\min}$ and $x \notin U^{[i]}Z^{[i]}_{\min}$. Thus to show $x \in \cF(i)$, it remains to check that $p_{[j]}(x) \notin U^{[j-1]}Z^{[j-1]}_{\min}$ for all $j \leq i$. We have $\dom(q^s_{(i)}) \subset \dom(q^s_{(i-1)}) = \cF(i-1)$ as $F(i-1)$ holds; thus $p_{[j]}(x) \notin U^{[j-1]}Z^{[j-1]}_{\min}$ for all $j \leq i-1$. If $p_{[i]}(x) \in U^{[i-1]}Z^{[i-1]}_{\min}$, then $p_{[i]}(x) \notin \cF(i-1) = \dom(q^s_{(i-1)})$ and this contradicts the fact that $x \in \cE(i)$. Therefore, we conclude $x \in \cF(i)$ as required.
	
	Suppose conversely that $x \in \cF(i)$. We claim that $x \notin U^{[i-1]}Z^{[i-1]}_{\min}$. If $x \in U^{[i-1]}Z^{[i-1]}_{\min}$, then there exists some $u \in U^{[i-1]}$ with $ux \in Z^{[i-1]}_{\min}$ and, as $Z^{[i-1]}_{\min}$ is closed and preserved by the $\lambda^{[i]}(\GG_m)$-action, also $p_{[i]}(ux) \in Z^{[i-1]}_{\min}$. The conjugation action of $\lambda^{[i]}(t)$ on $U^{[i-1]}$ has non-negative weights and so $u_0:=\lim_{t \ra 0} \lambda^{[i]}(t) u \lambda^{[i]}(t)^{-1}$ exists in $U^{[i-1]}$. Therefore, we deduce that $p_{[i]}(ux) = u_0 p_{[i]}(x) \in Z^{[i-1]}_{\min}$, which contradicts the assumption that $p_{[i]}(x) \notin U^{[i-1]}Z^{[i-1]}_{\min}$. Thus $x \notin U^{[i-1]}Z^{[i-1]}_{\min}$ as claimed. 
	
	Since by assumption $x \in \cF(i)$, we have $x \in \cF(i-1)= \dom(q^s_{(i-1)}) = \cE(i-1)$, where the equalities follow because $F(i-1)$ and $E(i-1)$ hold. Since $E(i)$ also holds, to show $x \in \dom(q^s_{(i)})$ it suffices to show $x \in \cE(i)$. For this, we already know that $x \in \cE(i-1)$ and  $x \notin U^{[i]}Z^{[i],s}_{\min}$ as $x \in \cF(i)$; thus, it remains to check that $\overline{x}:=p_{[i]}(x) \in \dom(q_{(i-1)}) = \cF(i-1)$. Since $x \in X^s_{\min} 
	$, we have $\overline{x}=p_{[i]}(x) \in X^{s}_{\min}$ by Lemma \ref{lemma consequences of filt} (vi). By assumption, $\overline{x}=p_{[i]}(x) \notin U^{[i-1]}Z_{\min}^{[i-1]}$. If there exists an index $j \leq i-1$, with $p_{[j]}(\overline{x}) \in U^{[j-1]}Z^{[j-1]}_{\min}$, then as $p_{[j]}(\overline{x}) = p_{[i]}(p_{[j]}(x))$ by Lemma \ref{lemma consequences of filt} (iii), we would deduce $p_{[j]}(x) \in U^{[j-1]}Z^{[j-1]}_{\min}$ from Assumption \eqref{split} of \hyperlink{Qis}{(QiS)}, which contradicts the fact that $x \in \dom(q^s_{(i-1)}) = \cF(i-1)$. Therefore, $\overline{x} \in \cF(i-1) = \dom(q^s_{(i-1)})$ as required.
\end{proof}

Before we prove the final implication, we prove the following short lemma relating the $\cE(i)$'s.

\begin{lemma}\label{lemma Es include}
Suppose that $D(i-2)$ holds. Then for $x \in \cE(i)$, we have $x \notin U^{[i-1]}Z^{[i-1]}_{\min}$. In particular, we have $\cE(i) \subset \cE(i-1)$.
\end{lemma}
\begin{proof}
For a contradiction, suppose $x \in U^{[i-1]}Z^{[i-1]}_{\min}$; then we can conclude $p_{[i]}(x) \in U^{[i-1]}Z^{[i-1]}_{\min}$ exactly as in the second paragraph in the proof of Proposition \ref{prop middle implication r}. Since $x \in \cE(i)$, we also know $p_{[i]}(x) \in \dom(q_{(i-1)}) \subset \dom(q_{(i-2)})$. By Lemma \ref{lem can take limits r} $(3)$, we conclude that $q_{(i-2)}(p_{[i]}(x)) \in U_{i-1}Z_{i-1,\min}$. Hence $q_{(i-2)}(p_{[i]}(x)) \notin \dom(q_{i-1})$, which contradicts $p_{[i]}(x) \in \dom(q_{(i-1)})$. 
\end{proof}

Finally we prove the left implication in \eqref{s chain of implications}.

\begin{prop}\label{prop left implication r}
 If $D(\leq i-1)$, $E(\leq i-1)$ and $F(\leq i-1)$ hold, then $E(i)$ holds; that is,
	\[ \dom(q^s_{(i)})= \cE(i) := \left\{x\in X^{s}_{\min} \setminus U^{[i]}Z^{[i],s}_{\min} \:\middle\vert\: p_{[j]}(x) \in \dom (q_{(j-1)}) \text{ for all } 2\leq j\leq i\right\} \]
\end{prop}
\begin{proof}
	Suppose $x \in \dom(q^s_{(i)})$; then $x \in X^{s}_{\min}$ by definition of $q_{(i)}^s$. For all $j \leq i-1$ we have $ \dom(q^s_{(i)}) \subset  \dom(q^s_{(j)}) = \cE(j)$, where the equality is given by $E(\leq i-1)$. In particular, $x \notin U^{[j]}Z^{[j]}_{\min}$ and $p_{[j]}(x) \in \dom(q_{(j-1)})$ for all $j \leq i-1$. Therefore, to show that $x \in \cE(i)$, it remains to check that $x \notin U^{[i]}Z^{[i]}_{\min}$ and $p_{[i]}(x) \in \dom(q_{(i-1)})$. If $x \in U^{[i]}Z^{[i]}_{\min}$, then as $D(i-1)$ holds, we have $q_{(i-1)}(x) \in U_i Z_{i,\min}$ by Lemma \ref{lem can take limits r} (3). Then $q_{(i-1)}(x) \notin \dom(q_i) \subset X^{ss}_{i,\min} \setminus U_i Z_{i,\min}^{ss}$, where this containment follows from Theorem \ref{thm nrGIT HM type inclusions} \emph{\ref{HM type inclusions}}. However, $\dom(q_{(i)}) = q_{(i-1)}^{-1}(\dom (q_i))$ and so this would contradict the fact that $x \in \dom(q^s_{(i)}) \subseteq \dom(q_{(i)})$. If $p_{[i]}(x) \notin \dom(q_{(i-1)})$, then $q_{(i-1)}(x) \notin X_{i,\min}$ by Lemma \ref{lem can take limits r} (2), which we can apply as $D(i-1)$ holds. This again contradicts the fact that $x \in \dom(q_{(i)})= q_{(i-1)}^{-1}(\dom (q_i))$, as $\dom(q_i) \subset X_{i,\min}$. Therefore, we conclude that $x \in \cE(i)$ as required.
	
Suppose conversely that $x \in \cE(i)$. Then by Lemma \ref{lemma Es include}, also $x \in \cE(i-1) =\dom(q^s_{(i-1)})$, where this equality is due to $E(i-1)$. Since $x \in \dom (q^s_{(i-1)})$, to prove that $x \in \dom(q^s_{(i)})$, it suffices to show that the point $y := q_{(i-1)}(x) \in X_i$ is stable for the $\hH_i$-action; that is, we need to show that $y \in X^s_{i,\min} \setminus U_i Z_{i,\min}^s$. By assumption, we have $p_{[i]}(x) \in \dom(q_{(i-1)})$ and so $y = q_{(i-1)}(x) \in X_{i,\min}$ by Lemma \ref{lem can take limits r} (2), which we can apply as $D(i-1)$ holds. Furthermore, $y \notin U_iZ_{i,\min}$, as $x \notin U^{[i-1]}Z^{[i]}_{\min}$ and so this follows by Lemma \ref{lem can take limits r} (2), which we can apply $D(i-1)$ holds. Therefore, it remains to check that $y \in X^s_{i,\min}$ or equivalently $p_i(y) \in Z_{i,\min}^s$; this last claim follows from Proposition \ref{prop claim 4}, because $p_i(y) = q_{(i-1)}(p_{[i]}(x))$ and $p_{[i]}(x) \in \dom(q_{(i-1)})$ and also $p_{[i]}(x) \in Z^{[i],s}_{\min} \cap X^{s}_{\min}$ by Lemma \ref{lemma consequences of filt} (v), as $x \in X^{s}_{\min} \subset X^{[i],s}_{\min}$.
\end{proof}

Now we can describe the domain of $q$, under assumptions on the downstairs stabilisers.

\begin{proof}[Proof of Theorem \ref{thm with downstairs assumptions}]
	Since \hyperlink{D}{(D)} holds, we have by Theorem \ref{mainthm nrGIT quotient} and the definition of $q_{(1)}^s$ that
		\[ \dom(q^s_{(1)}) = X^{s}_{\min} \setminus  U^{[1]}Z^{[1],s}_{\min}, \]
		which shows $E(1)$ and $F(1)$ hold. We prove that the domain of  $q^s $ is $\cF(l-1)$ inductively using Lemma \ref{lem a^s holds r}, together with Propositions \ref{prop middle implication r} and \ref{prop left implication r}. Since \hyperlink{D}{(D)} holds, $\dom(q) \subset X^{s}_{\min}$ by Corollary \ref{cor domain q in Xsmin} (see Remark \ref{rmk centre of L acts trivially on Zmin}) and consequently $\dom(q^s) =  \dom(q)$. Finally this set admits a projective geometric quotient by Corollary \ref{cor exists open set with qnt}.
\end{proof}

\subsection{From Downstairs to Upstairs Stabiliser Assumptions} \label{subsec from downstairs to upstairs}

In this subsection, we will prove Theorem \ref{mainthm2} in the case where \ref{starUstagesibullet} holds (see Proposition \ref{prop qnt in stages no blow ups}) by moving from the \emph{Downstairs Stabiliser Assumptions} \hyperlink{D}{(D)} to the \emph{Upstairs Stabiliser Assumptions} \hyperlink{U}{(U)}.

Our proof is split into two cases. In $\S$\ref{sec upstairs red holds}, when \hyperlink{U}{(U)} holds, we inductively show \hyperlink{D}{(D)} holds so we can apply Theorem \ref{thm with downstairs assumptions}. In $\S$\ref{sec upstairs red fails}, where only \hyperlink{UU}{(UU)} holds, we perform a (reductive) blow-up procedure on $X$ to obtain a sequence of blow ups $\widehat{X} \ra X$ such that \hyperlink{U}{(U)} holds on $\widehat{X}$ and then apply $\S$\ref{sec upstairs red holds}. We deal with the case when \hyperlink{UU}{(UU)} fails in $\S$\ref{sec quot stages with blowups}  below.

\subsubsection{Comparing upstairs and downstairs stabilisers} \label{subsec two stab comparison results}
	
	In this subsection, we will prove results relating the assumptions \hyperlink{D}{(D)} and \hyperlink{U}{(U)}. 
	
		Suppose we have performed Construction \ref{const qnt in stages} up to stage $j$, so we have the rational map \[q_{(j-1)}: X \dashrightarrow X_j\] which is equivariant with respect to the quotient map $\hH^{(j)} \twoheadrightarrow \hH_j$ and whose restriction to its domain of definition gives a geometric quotient by \begin{equation}\label{decomp of Hj-1}
		\hH^{(j-1)} = U^{(j-1)} \rtimes (R^{(j-1)} \times T^{(j-1)}) \quad \text{where} \quad  T^{(j-1)}=\lambda^{[1]}(\Gm)\times\dots\times\lambda^{[j-1]}(\Gm)
	\end{equation} 
	
Now we begin the work of relating the upstairs and downstairs stabilisers.
	
	\begin{lemma}\label{lemma stab exact seq} Suppose that \hyperlink{QiS}{(QiS)} holds, and that \hyperlink{D}{(D)} holds up to the $(j-1)$th stage, so that we can obtain a geometric quotient $q_{(j-1)}$ using Theorem \ref{thm with downstairs assumptions}. Then for $x \in X^{[j]}_{\min}\cap \dom (q_{(j-1)})$, we have an exact sequence
		\begin{equation}
			1 \rightarrow \stab_{U^{[j-1,j]}}(x) \rightarrow \stab_{U^{[j]}}(x)\rightarrow \stab_{U_{j}}( q_{(j-1)}(x)). \label{U stab ses}
		\end{equation}
		If moreover, $x \in Z^{[j]}_{\min}$, then this sequence is also right exact. 
	\end{lemma}
	\begin{proof}
Since $q_{(j-1)}$ is equivariant with respect to the quotient map $U^{[j]} \twoheadrightarrow U_j$, we obtain a natural morphism
		\[ \stab_{U^{[j]}}(x)\rightarrow \stab_{U_{j}}( q_{(j-1)}(x)) \]
		whose kernel is $\stab_{U^{[j]}}(x) \cap U^{[j-1]} = \stab_{U^{[j-1,j]}}(x)$, which proves the first claim.
		
		Now suppose additionally $x \in Z^{[j]}_{\min}$ and write $y := q_{(j-1)}(x)$. We claim this sequence is also right exact. Let $\bar{u} \in \stab_{U_{j}}(y)$; choose a lift $u \in U^{(j)}$ of $\bar{u}$ under the quotient map $U^{(j)} \twoheadrightarrow U_j$. By equivariance of $q_{(j-1)}$, we have $q_{(j-1)}(ux) = \bar{u} y = y = q_{(j-1)}(x)$, and as $q_{(j-1)}$ is a geometric $\hH^{(j-1)}$-quotient, there exists an element $h \in \hH^{(j-1)}$ such that $ux = hx$. Let us write $h = vrt$ with respect to the decomposition \eqref{decomp of Hj-1}, where $v \in U^{(j-1)}$, $r \in R^{(j-1)}$ and $t \in T^{(j-1)}$. Since $\lambda^{[j-1]}$ is a central 1-PS of the Levi $L$ and $R^{(j-1)} \times T^{(j-1)} < L$, the minimal weight space $Z^{[j]}_{\min}$ is preserved by the action of $R^{(j-1)} \times T^{(j-1)} $. Thus $v^{-1}ux = rtx \in Z^{[j]}_{\min}$. We claim there is an element $\tilde{h} \in U^{(j-1)}$ such that $\tilde{u} := \tilde{h} v^{-1}u \in U^{[j]}$ and $\tilde{u}x \in Z^{[j]}_{\min}$. Indeed if we write the block form of $v^{-1}u \in U^{(j)}$ with respect to the two-term weight filtration of $\lambda^{[j]}$: if
		\[ v^{-1}u = \left( \begin{array}{cc} 
		A & B \\ 0 & \id \end{array} \right), \quad \text{then} \quad  \tilde{h}:= \left( \begin{array}{cc}
		A^{-1} & 0 \\ 0 & \id \end{array} \right),\]
has the desired properties: $\tilde{h}$ fixes $Z^{[j]}_{\min}$ setwise since it lies in the Levi $R^{[j]}$ of the parabolic $P^{[j]} = P(\lambda^{[j]})$. Since $\tilde{u} \in U^{[j]}$ and $\tilde{u}x \in Z^{[j]}_{\min}$, it follows that $\tilde{u} \in \stab_{U^{[j]}}(x)$ by Lemma \ref{lem second small lemma}. Finally, the image of $\tilde{u} = \tilde{h} v^{-1}u$ under the quotient $U^{[j]} \twoheadrightarrow U_j$ is $\bar{u}$, as $\tilde{h} v^{-1} \in U^{(j-1)}$, which completes the surjectivity proof. 	
\end{proof}
	
\begin{prop}\label{prop stab relns U} 
In the situation of Lemma \ref{lemma stab exact seq}, suppose additionally that \ref{starUstagesibullet} holds. Take $x \in X^{[j],ss}_{\min}\cap \dom (q_{(j-1)})$ such that $p_{[j]}(x) \in \dom (q_{(j-1)})$ and $q_{(i-1)}(x) \in Y_{j,\min}$. Then 
		\[\dim \stab_{U_{j}}( q_{(j-1)}(x))=\dim \stab_{U^{[j]}}(x) - \dim \stab_{U^{[j-1,j]}}(x) .\]
	\end{prop}
	\begin{proof}

		Let us write $\overbar{x} := p_{[j]}(x)$ and $y:=q_{(j-1)}(x)$. Then $q_{(i-1)}(\overbar{x}) = \overbar{y}$ where $\overbar{y} = p_j(y)$. By applying Lemma \ref{lemma stab exact seq} to $\overbar{x} \in Z^{[j]}_{\min} \cap \dom (q_{(j-1)})$, we have 
		\[\dim \stab_{U_{j}}( \overbar{y})=\dim \stab_{U^{[j]}}(\overbar{x}) - \dim \stab_{U^{[j-1,j]}}(\overbar{x}), \]
		and by applying this lemma also to $x$, we have  
		\begin{equation}\label{ineq stab}
			\dim \stab_{U_{j}}( y) \geq \dim \stab_{U^{[j]}}(x) - \dim \stab_{U^{[j-1,j]}}(x).
		\end{equation}
		
		By the semicontinuity of dimensions of stabilisers, we have 

		\[\dim \stab_{U_{j}}(y) \leq \dim \stab_{U_{j}}( \overbar{y}).\]
		We claim that equality holds here and in  \eqref{ineq stab}: indeed, we have 
		\begin{equation*}\begin{split}
				\dim \stab_{U_{j}}( \overbar{y}) \geq \dim \stab_{U_{j}}(y) &\geq \dim \stab_{U^{[j]}}(x) - \dim \stab_{U^{[j-1,j]}}(x) \\ &=   \dim \stab_{U^{[j]}}(\overbar{x}) - \dim \stab_{U^{[j-1,j]}}(\overbar{x})\\ &= \dim \stab_{U_{j}}( \overbar{y}),
		\end{split} \end{equation*}
		where the second to last equality follows from \ref{starUstagesibullet} as $x \in X^{[j],ss}_{\min}$. 
	\end{proof}
	
\subsubsection{The case when the Upstairs Stabiliser Assumption holds}\label{sec upstairs red holds}
In this section, we assume the $P$-action on $X$ satisfies \hyperlink{QiS}{(QiS)} and also \hyperlink{U}{(U)}.

\begin{prop}\label{prop no blow ups get proj qnt} 
If \hyperlink{U}{(U)} holds for the linearised $P$-action on $X$, then the Quotienting-in-Stages procedure yields a projective geometric $P$-quotient of the locus	
\[\dom (q^s) = \dom (q) = X^{P-qs}. \]  	
\end{prop}
\begin{proof}
We will inductively show that if \hyperlink{U}{(U)} holds, then \hyperlink{D}{(D)} and thus, the result follows from Theorem \ref{thm with downstairs assumptions}, as when  \hyperlink{U}{(U)} holds we have
\[ X^{P-qs}  = \left\{x\in X^{s}_{\min}  \:\middle\vert\:p_{[j]}(x) \notin U^{[j-1]}Z^{[j-1]}_{\min} \text{ for all } 2\leq  j\leq l\right\}. \]

For the base case, if $l =2$, then there is only one stage in the Quotienting-in-Stages procedure and the Upstairs and Downstairs Stabiliser Assumptions (that is, \hyperlink{U}{(U)} and \hyperlink{D}{(D)}) coincide. For the inductive step, we fix $i$ and assume that the conditions in \hyperlink{D}{(D)} hold for the groups $U_j$ and $R_j$ for all $j\leq i$, then 
\[q_{(i)} : X \dashrightarrow X_{i+1}\]
is a geometric $\hH^{(i)}$-quotient of its domain of definition, which is $\dom(q_{(i)}^s)=\cF(i)$ following the inductive proof of Theorem \ref{thm with downstairs assumptions} up to this stage. We need to show that we can deduce the stabiliser conditions in \hyperlink{D}{(D)} for the groups $U_{i+1}$ and $R_{i+1}$ acting on $X_{i+1}$ by using the Assumption \hyperlink{U}{(U)} upstairs on $X$. For any $y\in X^s_{i+1,\min}$ we can write $y = q_{(i)}(x)$ for some $x \in \cF(i) \subset X^{s}_{\min} \subset X^{[i],s}_{\min}$. By Lemma \ref{lem a^s holds r} we can apply Lemma \ref{lem can take limits r} (2) to deduce that $p_{[i+1]}(x) \in \dom (q_{(i)})$. Hence we can apply Proposition \ref{prop stab relns U} to the map $q_{(i)} : X \dashrightarrow X_{i+1}$ to obtain an expression for $\dim\stab_{U_{i+1}}(y)$ in terms of the dimensions of the stabilisers of $x$ under $U^{[i+1]}$ and $U^{[i,i+1]}$. Since \hyperlink{U}{(U)} holds, the latter dimensions are constant on $X^{[i],s}_{\min}$ and thus we see that $\dim\stab_{U_{i+1}}$ is fixed on $X^s_{i+1,\min}$; that is the part of \hyperlink{D}{(D)} concerning the $U_{i+1}$-stabilisers on $X_{i+1}$ holds.

To prove the part of \hyperlink{D}{(D)} concerning the stabilisers for the reductive group $R_{i+1}$ acting on $X_{i+1}$, we need to show $Z_{i+1,\min}^{s} = Z_{i+1,\min}^{ss}$. By Lemma \ref{lem can take limits r} (3), any $z\in Z_{i+1,\min}^{ss}$ is the image under $q_{(i)}$ of some $w \in Z^{[i+1]}_{\min}$. Moreover, as $z$ is $R_{i+1}$-semistable, the point $w$ is $R^{(i+1)}$-semistable by pulling back invariants along $q_{(i)}$ using Lemma \ref{equiv pullback of semistable loci}. To prove that $z\in Z_{i+1,\min}^{s}$, we will show  first that it suffices to show $w \in X^{(i+1),ss}_{\min}$, and then show that $w \in X^{(i+1),ss}_{\min}$.

We claim that if $w \in X^{(i+1),ss}_{\min}$, then $z\in Z_{i+1,\min}^{s}$. In this case, we have $w \in X^{(i+1),ss}_{\min} =  X^{(i+1),s}_{\min} \subset X^{[i],s}_{\min}$, where the equality follows from Assumption \ref{starR0stagesibullet} and the inclusion follows from Remark \ref{rmk consequences of filt}. Since $w \in Z^{[i+1]}_{\min}$, we see that $w \in \dom(q_{(i)}) \cap Z^{[i+1],s}_{\min} \cap X^{(i+1),s}_{\min}$ and deduce that $z = q_{(i)}(w) \in Z_{i+1,\min}^{s}$ from Proposition \ref{prop claim 4}.

It remains to show that $w\in X^{(i+1),ss}_{\min}$. For this, we will show $w$ lies in the domain of definition for a NRGIT quotient $\phi$ for a larger group 
 \[ \hH^{\prime} : = U^{(i)}\rtimes (R^{(i+1)}\times T^{(i+1)}) = \hH^{(i)}\rtimes (R_{i+1} \times \lambda^{[i+1]}(\GG_m))\] 
using arguments similar to Lemma \ref{lem an induction for claim 4}. Then we will show that  that $\dom (\phi) \subset X^{(i),ss}_{\min}$ using arguments similar to Corollary \ref{cor domain q in Xsmin}. Recall that by our inductive assumption, $q_{(i)} : X \dashrightarrow X_{i+1}$ is a geometric $\hH^{(i)}$-quotient of its domain of definition and $X_{i+1}$ is the projective spectrum of the $\hH^{(i)}$-invariants on $X$, as each stage in the quotient satisfies the relevant downstairs assumptions of \hyperlink{D}{(D)} up to this stage. We consider the following composition of quotients 
\[\phi : X \dashrightarrow X_{i+1} \dashrightarrow Y:=X_{i+1} \git (R_{i+1}\times \lambda^{[i]}(\GG_m)),\] 
where the first quotient is $q_{(i)}$ and the second is a reductive GIT quotient. Then $\phi $ is also a NRGIT quotient by $ \hH^{\prime}$ and thus $Y$ is the projective spectrum of the $ \hH^{\prime}$-invariants on $X$. In fact, we will take these quotients with respect to a perturbation $\cL_i \ra X$ of the linearisation. Recall that we obtained a character $\chi$ to twist the borderline linearisation by in Proposition \ref{prop character for Q-i-S} from certain choices $\epsilon_j$. We let $\chi_i$ be the character obtained from this construction by setting $\epsilon_j = 0$ for $j \geq i + 1$ and let $\cL_i$ be the twist of the borderline linearisation by this character. By construction, the weight of $\lambda_{i+1}$ on $Z_{i+1,\min} \subset X_{i+1}$ is zero with respect to $\cL_i$, as the central torus $T$ acts on $Z_{\min}$ with weight zero with respect to the borderline linearisation. Therefore,
\[ Z_{i+1,\min}^{ss} := (Z_{i+1,\min})^{R_{i+1}-ss} = (Z_{i+1,\min})^{(R_{i+1}\times \lambda_{i+1}(\GG_m))-ss}(\cL_i) \subset X_{i+1}^{(R_{i+1}\times \lambda_{i+1}(\GG_m))-ss}(\cL_i),\] 
where we also use $\cL_i$ to denote the induced linearisation on $X_{i+1}$ and the last inclusion follows as $Z_{i+1,\min} \hookrightarrow X_{i+1}$ is an equivariant closed embedding. Thus, since $q_{(i)}(w) = z \in Z_{i+1,\min}^{ss}$ by assumption, we deduce $w \in \dom (\phi)$. To conclude that $w\in X^{(i+1),ss}_{\min}$, we will show that $\dom(\phi) \subset X^{(i+1),ss}_{\min}$. By Assumption \eqref{filt} of  \hyperlink{QiS}{(QiS)}, the torus $T^{(i+1)}$ acts trivially on $Z^{(i+1)}_{\min}$ and this torus must act with weight zero on $Z^{(i+1)}_{\min}$ for the borderline linearisation $\cL_0$, as $Z_{\min} \subset Z^{(i+1)}_{\min}$. Hence,
 \[X_{i+1}^{(R_{i+1}\times \lambda_{i+1}(\GG_m))-ss}(\cL_i) \subset X_{i+1}^{(R_{i+1}\times \lambda_{i+1}(\GG_m))-ss}(\cL_0) = p_{(i+1)}^{-1}(Z^{(i+1),ss}_{\min}) =X^{(i+1),ss}_{\min}\] where the first inclusion is by the same VGIT argument as Proposition \ref{prop character for Q-i-S}. Thus we see in the same way as Corollary \ref{cor domain q in Xsmin} that $\dom (\phi) \subset X^{(i),ss}_{\min}$, which concludes our proof that $z \in Z_{i+1,\min}^s$.

By induction we deduce \hyperlink{D}{(D)} from the \hyperlink{U}{(U)} and the result then follows from Theorem \ref{thm with downstairs assumptions}.
\end{proof}

\subsubsection{The case when only Upstairs Unipotent Stabiliser Assumption holds}
\label{sec upstairs red fails}
In this subsection, we assume that the \emph{Upstairs Unipotent Stabiliser Assumption} \ref{starUstagesibullet} holds, but the \emph{Upstairs Reductive Stabiliser Assumption} \ref{starR0stagesibullet} does not. Our approach will be to reduce to the case where both \ref{starUstagesibullet} and \ref{starR0stagesibullet} hold by performing an equivariant blow-up sequence which involves considering dimensions of stabilisers for the action of the reductive group $R$ on $X$ analogous to the reductive partial desingularisation procedure \cite{K2} outlined in $\S$\ref{sec red part desing}.

Since we will need to consider different minimal weight spaces on successive blow-ups of $X$, we introduce the following notation.

\begin{defn}
For any projective scheme $Y$ with an ample linearised $P$-action, we write $Z_{\min}^{[i]}(Y):=Z_{\min}(Y,\lambda^{[i]})$ for the $\lambda^{[i]}$-minimal weight space in $Y$ and let $Z_{\min}^{[i],(s)s}(Y)$ denote the $R^{(i)}$-(semi)stable locus. We write $p_{[i],Y} : Y_{\min}^{[i]} \ra Z_{\min}^{[i]}(Y)$ for the retraction onto the minimal weight space (or sometimes just $p_{[i]}$ if $Y$ is fixed). Similarly, we write $Z_{\min}^{(i),(ss)}(Y)$ for the semistable $\lambda^(i)$-minimal weight space and $Y_{\min}^{(i),(ss)}$ for the corresponding attracting open. 
\end{defn}

Our goal is to perform sequences of blow-ups for each reductive group $R^{(i)}$ inductively to arrange that $Z^{(i),s}_{\min} = Z^{(i),ss}_{\min}$ on the blow-up. Using the notation of Definition \ref{def centres of blowups}, we introduce the following schemes to define the centre of each blow-up.

\begin{defn}
For a projective scheme $Y$ with an ample linearised $P$-action, we define 
\[ C^{(i)}(Y):=C(Z^{(i),ss}_{\min}(Y),R^{(i)}) \quad \text{and} \quad B^{(i)}(Y):=p_{(i),Y}^{-1}(C^{(i)}(Y)) \quad \quad \text{for } \: 1 \leq i \leq l-1. \]
\end{defn}

By definition, $C^{(i)}(Y)$ is the closed subscheme of $Z^{(i),ss}_{\min}(Y)$ on which the dimension of the $R^{(i)}$-stabiliser groups are maximal; provided there is a strictly semistable point, this is a proper subscheme disjoint from the stable locus $Z^{(i),s}_{\min}(Y)$. In the reductive partial desingularisation procedure for the $R^{(i)}$-action on $Z^{(i),ss}_{\min}(Y)$, this scheme (or its closure in $Z_{\min}^{(i)}(Y)$, strictly speaking) is the centre of the first blow-up. Hence also $B^{(i)}(Y)$ is a closed subscheme of $Y_{\min}^{(i),ss}$. 

\begin{lemma}
The schemes $B^{(i)}(Y)$, and thus also their closures, are $P$-invariant.
\end{lemma}
\begin{proof}
Let us just prove the case for $i = 1$ as the other cases are similar. It suffices to show invariance under the group $\hH = H \rtimes T$, as this group surjects onto $P$. Since the central torus $T$ commutes with $R^{(1)}$, we just need to show invariance under $H$. We argue as in Proposition \ref{prop middle implication r}: if $h\in H$ and $y\in Y^{(1)}_{\min}$ then $p_{(1)}(hy) =h_0p_{(1)}(y)$ where $h_0:=\lim_{t \ra 0} \lambda^{(1)}(t)\cdot h \cdot \lambda^{(1)}(t)^{-1}$. We can write this element as $h_0 = h_1r_1$ where $r_1 \in R^{(1)}$ and $h_1$ has zero first row except for an identity matrix in the first block; thus $h_1$ commutes with $R^{(1)}$. Therefore,  
\[\stab_{R^{(1)}}(p_{(1)}(hy)) = \stab_{R^{(1)}}(h_1r_1p_{(1)}(y)) = r_1\stab_{R^{(1)}}(p_{(1)}(y))r_1^{-1}\] 
from which the $H$-invariance follows.
\end{proof}

We can now define the first blow-up in this procedure, assuming that $Z^{(1),s}_{\min} \neq Z^{(1),ss}_{\min}$.

\begin{defn}
Let $\widehat{X}_{(1)}$ denote the blow-up of $X$ along the closure of $B^{(1)}(X)$, with exceptional divisor denoted by $E_{(1)}$. As the centre of this blow-up is $P$-invariant, there is an induced $P$-action on $\widehat{X}_{(1)}$ which we linearise by the pulling back the linearisation on $X$ and perturbing by a small multiple of the exceptional divisor $E^{1}$ as in \cite{K2}. 
\end{defn}

\begin{prop}\label{prop properties first blowup}
Assume that $X$ satisfies \hyperlink{QiS}{(QiS)} and \ref{starUstagesibullet}. If $Z^{(1),s}_{\min} \neq Z^{(1),ss}_{\min}$, then the blow-up $\pi_1: \widehat{X}_{(1)} \ra X$ has the following properties.
\begin{enumerate}[label=\emph{(\alph*)}]
\item $X^{s}_{\min}$ is disjoint from the centre of the blow-up $\overline{B^{(1)}(X)}$,
\item The strict transform of $Z_{\min}^{[i]}$ (resp. $Z_{\min}^{(i)}$) is $Z_{\min}^{[i]}(\widehat{X}_{(1)})$ (resp. $Z_{\min}^{(i)}(\widehat{X}_{(1)})$) for all $1 \leq i \leq l-1$.
\item $\widehat{X}_{(1)}$ also satisfies \hyperlink{QiS}{(QiS)} and \ref{starUstagesibullet} .
\item The dimensions of $R^{(1)}$-stabilisers on the semistable $\lambda^{(1)}$-minimal weight space drops:
\[\max_{x\in {Z}_{\min}^{(1),ss}(\widehat{X}_{(1)})} \dim \stab_{R^{(1)}}(x) < \max_{x\in Z^{(1),ss}_{\min}} \dim \stab_{R^{(1)}} (x).\]
\end{enumerate}
\end{prop}
\begin{proof}
For (a), we note that $Z^{(1),s}_{\min}$ is disjoint from $C^{(1)}(X)$, as stable points have zero dimensional stabiliser and we assumed $Z^{(1),s}_{\min} \neq Z^{(1),ss}_{\min}$. Thus $B^{(1)}(X):=p_{(1)}^{-1}(C^{[1]}(X))$ is disjoint from $X^{(1),s}_{\min}$. We then conclude (a) as $X^{s}_{\min}  \subset X^{(1),s}_{\min}$ by Remark \ref{rmk consequences of filt}. Note that $X^{s}_{\min}  \neq \emptyset$, by Assumption \eqref{non-deg} of \hyperlink{QiS}{(QiS)}.

For (b), we have seen that $Z^{s}_{\min}$ is disjoint from the centre of the blow-up in (a) and since $Z^{s}_{\min} \subset Z^{[i],s}_{\min}$ (resp. $Z^{s}_{\min} \subset Z^{(i),s}_{\min}$), we see that $Z^{[i]}_{\min}$ (resp. $Z_{\min}$) is not contained in the centre of the blow-up. Consequently, we deduce (b).

For (c), to deduce \ref{starUstagesibullet} holds for $\widehat{X}_{(1)}$, we use Lemma \ref{lem stab cant incr on blowup} for the groups $U^{[i]}$ and $U^{[i-1,i]}$, together with the fact that the dimensions of these stabilisers on $X^{[i],ss}_{\min}$ are constant. We can directly deduce that Assumptions \eqref{filt} and \eqref{split} of \hyperlink{QiS}{(QiS)} hold for $\widehat{X}_{(1)}$ as the strict transform of $Z_{\min}^{[i]}$ (resp. $Z_{\min}^{(i)}$)  is $Z_{\min}^{[i]}(\widehat{X})$ (resp. $Z_{\min}^{(i)}(\widehat{X})$) by (c).
Finally, non-degeneracy \eqref{non-deg} of \hyperlink{QiS}{(QiS)} holds for $\widehat{X}_{(1)}$ as we can identify $X^{s}_{\min}$ with a non-empty open set of $\widehat{X}_{(1),\min}^s$ by (a).

For (d), we can use (b) to identify $Z_{\min}^{(1)}(\widehat{X})$ with the strict transform of $Z_{\min}^{(1)}$. Then the result follows exactly as in \cite[Lemma 6.1 (iv)]{K2}.
\end{proof}

\begin{prop}\label{prop blowup process upstairs red fails}
Assume that the $P$-action on $X$ satisfies  \hyperlink{QiS}{(QiS)} and \ref{starUstagesibullet}. Then there exists a sequence of $P$-equivariant blow-ups resulting in a projective scheme $\widehat{X} \ra X$ such that $\widehat{X}$ satisfies \hyperlink{QiS}{(QiS)} and \hyperlink{U}{(U)}.
\end{prop}
\begin{proof}
We  iterate the above process for $R^{(1)}$: let $\widehat{X}_{(0)} = X$ and for $i =1,\dots ,n_1$, define $\widehat{X}_{(i)}$ to be the blow-up of $\widehat{X}_{(i-1)}$ along the closure of $B^{(1)}(\widehat{X}_{(i-1)})$. Since the dimensions of the $R^{(1)}$-stabilisers on the semistable locus in the $\lambda^{(1)}$-minimal weight space decreases at each stage by Proposition \ref{prop properties first blowup}, this procedure terminates with a scheme $\widehat{X}_{(n_1)}$ for which stability and semistability coincides for the $R^{(1)}$-action on the minimal $\lambda^{(1)}$-weight space. 

We then turn our attention to $R^{(2)}$ and for $i=1,\dots, n_2$, define $\widehat{X}_{(n_1+i)}$ to be the blow-up of $\widehat{X}_{(n_1 + i-1)}$ along the closure of $B^{(2)}(\widehat{X}_{(n_1 + i-1)})$. This terminates with a scheme $\widehat{X}_{(n_1 +n_2)}$ for which stability and semistability coincides for the $R^{(j)}$-action on the minimal $\lambda^{(j)}$-weight space  for $j \leq 2$. Once we have completed the blow-ups for $R^{(j)}$, we turn our attention to $R^{(j+1)}$ and we progressively construct a sequence of blow-ups
\[ \pi : \widehat{X} =\widehat{X}_{(n_1 +\dots + n_{l-1})} \ra \cdots \ra \widehat{X}_{(n_1 +n_2)} \ra \cdots \ra \widehat{X}_{(n_1)} \ra \cdots \ra X\]
such that stability and semistability coincides for the $R^{(j)}$-action on the minimal $\lambda^{(j)}$-weight space in $\widehat{X}$ for $1 \leq j \leq l-1$; that is the Upstairs Reductive Stabiliser Assumption \ref{starR0stagesibullet} holds on $\widehat{X}$. Furthermore, by Proposition \ref{prop properties first blowup} (c), \hyperlink{QiS}{(QiS)} and \ref{starUstagesibullet} also hold for $\widehat{X}$.
\end{proof}

Now we complete the proof of Theorem \ref{mainthm2} in the case where \ref{starUstagesibullet} holds but \ref{starR0stagesibullet} may fail.

\begin{prop} \label{prop qnt in stages no blow ups}
	Suppose that \hyperlink{QiS}{(QiS)} and \ref{starUstagesibullet} hold for the linearised $P$-action on $X$. 
\bnu
\item If \ref{starR0stagesibullet} holds, we have $ \dom (q)=X^{P-qs}$, which admits a projective geometric $P$-quotient.
\item If \ref{starR0stagesibullet} fails, there is an equivariant sequence of blow-ups $\pi:\hX \rar X$ such that $\hX$ satisfies the conditions of (1). Moreover, $X^{P-qs}$ is isomorphic to an open subset of $\hX^{P-qs}$ and thus $X^{P-qs}$ has a quasi-projective geometric quotient by $\hH$, with a canonical projective completion given by the geometric $P$-quotient of $\hX^{P-qs}$. \enu 
\end{prop}
\begin{proof}
In the case where \ref{starR0stagesibullet} holds, this is Proposition \ref{prop no blow ups get proj qnt}. If \ref{starR0stagesibullet} fails, we perform the blow-up process described in Proposition \ref{prop blowup process upstairs red fails} above, to obtain $\pi: \hX \rar X$ such that $\hX$ satisfies the Quotienting-in-Stages and \hyperlink{U}{(U)} (both reductive and unipotent). Hence, $\hX$ satisfies the conditions of Proposition \ref{prop no blow ups get proj qnt} and we obtain a projective geometric $\hH$ quotient of the locus \[\hX^{P-qs}= \{x\in \hX^s_{\min} \mid \hp_{[j]}(x) \notin U^{[j]}\hZ_{\min}^{[j-1]} \text{ for all } 1\leq j \leq l\} .\] We have $\pi^{-1}(Z^{s}_{\min}) \subset Z^{s}_{\min}(\hX)$ and since by Proposition \ref{prop properties first blowup} (a), $X^{s}_{\min}$ is disjoint from the centre of this blow-up, we have also $X^s_{\min} \cong \pi^{-1}(X^s_{\min})\subset \hX^s_{\min}$ . It then follows that $X^{P-qs} \cong \pi^{-1}(X^{P-qs}) \subset \hX^{P-qs}$, and hence we get a quasi-projective geometric $\hH$-quotient of $X^{P-qs}$ by restricting the quotient we obtained for $\hX^{P-qs}$. 
\end{proof} 
	
\subsection{The case when the Upstairs Unipotent Stabiliser Assumption fails}\label{sec quot stages with blowups}  We now conclude the proof of Theorem \ref{mainthm2} section by proving Theorem \ref{prop blowup process if upstairs unip fails}, which describes what to do when \ref{starUstagesibullet} fails, but \hyperlink{WUU}{(WUU)} holds. 
  
Recall from Definition \ref{defn quotienting-in-stages stable locus} that 
\[ d^{[i]}_{\min}:= d_{\min}(X^{[i],ss}_{\min},U^{[i]}) \quad \text{and} \quad d^{[i-1,i]}_{\min}:= d_{\min}(X^{[i],ss}_{\min},U^{[i-1,i]})\]
are respectively the minimal stabiliser dimensions of $U^{[i]}$ and $U^{[i-1,i]}$ on $X^{[i],ss}_{\min}$. Similarly, we will use the notation $d^{[i]}_{\max}$ and $d^{[i-1,i]}_{\max}$ for the maximal stabiliser dimensions of these subgroups on $X^{[i],ss}_{\min}$. In the case already considered, where \ref{starUstagesibullet} holds, the dimensions of the stabilisers for $U^{[i]}$ and $U^{[i-1,i]}$ on $X^{[i],ss}_{\min}$ are constant, so these maximal and minimal dimensions coincide. In the case where \ref{starUstagesibullet} fails, we will use a blow-up procedure to reduce to a situation where \ref{starUstagesibullet} holds; however, to determine the minimal weight spaces in the blow-ups we need to assume \hyperlink{WUU}{(WUU)} holds.

Using the notation of Definition \ref{def centres of blowups}, we introduce the following schemes to define the centre of each blow-up.

\begin{defn}
For a projective scheme $Y$ with an ample linearised $P$-action, we define 
\[ D^{[i]}(Y):=C(X^{[i],ss}_{\min}(Y),\hU^{[i]}) \quad \text{and} \quad D^{[i-1,i]}(Y):=C(X^{[i],ss}_{\min}(Y),\hU^{[i-1,i]})  \quad \text{for } \: 1 \leq i \leq l-1 \]
where $\hU^{[i-1,i]}:= U^{[i-1,i]} \rtimes \lambda^{[i]}(\GG_m)$.
\end{defn}

By definition, $D^{[i]}(Y)$ is the closed subscheme of $X^{[i],ss}_{\min}(Y)$ on which the dimension of the stabiliser for $\hU^{[i]}$ are maximal; equivalently this is the $U^{[i]}$-sweep of $C(Z^{[i],ss}_{\min}(Y),U^{[i]})$.

\begin{lemma}
The schemes $D^{[i]}(Y)$ and $D^{[i-1,i]}(Y)$, and thus also their closures, are $P$-invariant.
\end{lemma}
\begin{proof}
This follows as $U^{[i]}$ is a normal subgroup of $P$. Indeed, $U^{[i]}$ is the unipotent radical of $P^{[i]}:= P(\lambda^{[i]})$ and so $U^{[i]} \vartriangleleft P^{[i]}$. Since $P < P^{[i]}$, we deduce the claim. 
\end{proof}

Let us examine the first step in the blow-up procedure. The assumption \ref{cond not all of Zmin blownup} in \hyperlink{WUU}{(WUU)} ensures that not all of $Z_{\min}$ is blown up, allowing us to identify the minimal weight spaces in the blown up space as the proper transforms of their counterparts on $X$.

\begin{prop}\label{prop properties first nonred blowup}
Assume that $X$ satisfies \hyperlink{QiS}{(QiS)} and \hyperlink{WUU}{(WUU)}. Suppose that $d_{\min}^{[1]}< d_{\max}^{[1]}$. Then the blow-up $\pi_1: \widetilde{X}_{(1)} \rar X$ of $X$ along the closure $\overbar{D^{[1]}(X)}\subset X$ has the following properties.
\begin{enumerate}[label=\emph{(\alph*)}]
\item The locus $C(X^s_{\min},U^{[1]},d^{[1]}_{\min})$ is non-empty and disjoint from the centre of the blow-up. In particular, so is $X^{P-qs}$.
\item The strict transform of $Z^{[1]}_{\min}$ in $\widetilde{X}$ coincides with the blow-up $\widetilde{Z}_{(1)}$ of $Z^{[1]}_{\min}$ along the closure of $C(Z^{[1],ss}_{\min},U^{[1]})$.
\item The strict transform of $Z_{\min}^{[i]}$ (resp. $Z^{(i)}_{\min}$) is $Z_{\min}^{[i]}(\widetilde{X}_{(1)})$ (resp. $Z_{\min}^{(i)}(\widetilde{X}_{(1)})$) for $1 \leq i \leq l-1$.
\item The scheme $\widetilde{X}_{(1)}$ also satisfies \hyperlink{QiS}{(QiS)} and \hyperlink{WUU}{(WUU)}.
\item The dimensions of $U^{[1]}$-stabilisers in the $\lambda_{[1]}$-minimal weight spaces drops:
\[\max_{z\in {Z}_{\min}^{[1],ss}(\widetilde{X}_{(1)})} \dim \stab_{U^{[1]}}(z) < \max_{z\in Z^{[1],ss}_{\min}} \dim \stab_{U^{[1]}} (z).\]
\end{enumerate} 
\end{prop}
\begin{proof}  
Assumption \hyperlink{WUU}{(WUU)} implies that $C(X^s_{\min},U^{[1]},d^{[1]}_{\min}) \neq \emptyset$ and this is disjoint from the centre of the blow-up, since only points with $U^{[1]}$-stabiliser dimension equal to $d_{\max}^{[1]}$ are blown up. This proves (a). 

To prove (b), first we claim that $Z^{[1]}_{\min}$ is not contained in the centre of the blow-up. Recall that the centre of the blow up is \[D^{[1]}(X):=C(X^{[1],ss}_{\min}(X),\hU^{[1]}) = U^{[1]}C(Z^{[1],ss}_{\min}(X),U^{[1]}).\] By Lemma \ref{lem second small lemma} (see also \cite{BDHK2} Lemma 5.2), we have for $u\in U^{[1]}$ and $z\in Z^{[1]}_{\min}$ that $uz \in Z^{[1]}_{\min}$ if and only if $uz =z$. Hence, we deduce that \[D^{[1]}(X)\cap Z^{[1]}_{\min} =  C(Z^{[1],ss}_{\min}(X),U^{[1]}).\] Thus it suffices to show that $\overbar{C(Z^{[1],ss}_{\min}(X),U^{[1]})} \subsetneq Z^{[1]}_{\min}$ is a proper subscheme. But this consists only of points with $U^{[1]}$-stabiliser dimension equal to $d_{\max}^{[1]}$, by semicontinuity of stabiliser dimension. On the other hand, by \hyperlink{WUU}{(WUU)} we know that $C(Z_{\min},U^{[1]},d^{[1]}_{\min}) \neq \emptyset$. Since  $Z_{\min} \subset Z_{\min}^{[1]}$ by \hyperlink{QiS}{(QiS)} \eqref{filt}, we see that $Z_{\min}$ is not contained in the centre of the blow-up and $C(Z^{[1]}_{\min},U^{[1]},d^{[1]}_{\min}) \neq \emptyset$. It then follows that $Z^{[1]}_{\min}$ is not contained in the centre of the blow-up, which proves (b).

To prove (c), the argument is the same as (b). By \hyperlink{WUU}{(WUU)} we know that $Z_{\min}$ is not contained in the centre of the blow-up and by \hyperlink{QiS}{(QiS)} \eqref{filt} we have $Z_{\min}\subset Z^{[i]}_{\min}$ and $Z_{\min}\subset Z^{(i)}_{\min}$, which shows both $ Z^{[i]}_{\min}$ and $Z^{(i)}_{\min}$ are not contained in the centre of the blow-up and proves (c).

 For (d), as taking proper transforms preserves containment, Assumptions \eqref{filt} and \eqref{split} of  \hyperlink{QiS}{(QiS)} for $\widetilde{X}_{(1)}$ follow from the corresponding properties for $X$ as (c) holds. For the non-degeneracy assumption \eqref{non-deg} of \hyperlink{QiS}{(QiS)}, take any $x \in X^{P-qs}$; this set is non-empty by \hyperlink{WUU}{(WUU)}. By (a), $X^{P-qs}$ is disjoint from the centre of the blow-up, so we can use $\pi_1^{-1}$ to lift $x \in X^{P-qs}$ to a unique $\tilde{x} \in \widetilde{X}_{(1)}.$ We claim that $\tilde{x} \in (\widetilde{X}_{(1)})^s_{\min}$. Let us use $\tilde{p}$, $\tilde{p}_{[i]}$ to denote the corresponding maps for $\widetilde{X}_{(1)}$ rather than $X$.  To prove our claim first note that, since $p(x)$ is not blown up, we have $\tilde{p}(\tilde{x}) \in Z_{\min}(\widetilde{X}_{(1)})$. Furthermore, the linearisation on $\widehat{X}_{(1)}$ is an arbitrarily small perturbation of the one on $X$, so $\tilde{p}(\tilde{x})$ is $R$-stable, meaning that $\tilde{x} \in (\widetilde{X}_{(1)})^s_{\min}$. It follows that $\tilde{x}$ is a point witnessing non-degeneracy for $\widetilde{X}_{(1)}$, since otherwise $\pi_1(\tilde{x})$ could not lie in $X^{P-qs}$. Hence \hyperlink{QiS}{(QiS)} holds on $\widetilde{X}_{(1)}$. 
 
For \hyperlink{WUU}{(WUU)}, observe that $\tilde{p}(\tilde{x})$ has the same stabilisers as its image under $\pi_1$, since that image is not blown up. This shows that \ref{cond not all of Zmin blownup} still holds. Moreover, the same is true for each $\tilde{p}_{[i]}(\tilde{x})$ and its image under $\pi_1$, which together with the above argument gives $\tilde{x}\in (\widetilde{X}_{(1)})^{P-qs}$, so this set is non-empty.  This proves \hyperlink{WUU}{(WUU)} for  $\widetilde{X}_{(1)}$, which completes the proof of (d). 
 
The proof of statement (e) is the same as that of Lemma \ref{lemma uni stab drops blowups}.
\end{proof}

The above proposition gives the first step in the blow-up procedure for $U^{[1]}$, which we can then iterate until $d^{[1]}_{\min} = d^{[1]}_{\min}$, at which point we do the same for $U^{[2]}$, and so on until we reach $U^{[l-1]}$. We then turn our attention to the groups $U^{[i-1,i]}$, proceeding in the same way for those. At each stage in the blow-ups for the groups $U^{[i]}$ and $U^{[i-1,i]}$, we need to ensure that not all of $Z_{\min}$ is blown up so that the final blow-up satisfies \hyperlink{QiS}{(QiS)}: for this we again use condition \ref{cond not all of Zmin blownup} of \hyperlink{WUU}{(WUU)}.

\begin{thm}\label{prop blowup process if upstairs unip fails}
Assume that \hyperlink{QiS}{(QiS)} and \hyperlink{WUU}{(WUU)} hold for the linearised $P$-action on $X$. Then the following statements hold.
\begin{enumerate}
\item There exists a sequence of $P$-equivariant blow-ups resulting in a projective scheme \[\pi :\widetilde{X} \ra X\] such that $\widetilde{X}$ satisfies \hyperlink{QiS}{(QiS)} and \ref{starUstagesibullet}.
\item The Quotienting-in-Stages stable locus $ X^{P-qs}$ admits a quasi-projective geometric $P$-quotient with a natural projective completion.
\end{enumerate}
\end{thm}
\begin{proof}
We inductively consider stabilisers for the groups $U^{[i]}$ and $U^{[i-1,i]}$. The first group to consider is $U^{[1]}$ and the first blow-up in this procedure is performed in Proposition \ref{prop properties first nonred blowup}. We iterate this procedure for $U^{[1]}$:  let $\widetilde{X}_{(0)} = X$ and for $i =1,\dots n_1$, define $\widehat{X}_{(i)}$ to be the blow-up of $\widehat{X}_{(i-1)}$ along the closure of $D^{[1]}(\widehat{X}_{(i-1)})$. Since the dimensions of the $U^{[1]}$-stabilisers on the semistable locus in the $\lambda^{[1]}$-minimal weight space decreases at each stage, this procedure terminates with a scheme $\widehat{X}_{(n_1)}$ which has constant dimensional $U^{[1]}$-stabilisers on the semistable locus in the $\lambda^{[1]}$-minimal weight space equal to $d^{[1]}_{\min}$. As in the proof of Proposition \ref{prop ab uni blowups seq}, it follows that $\dim \stab_{U^{[1]}}$ is constant and equal to $d^{[1]}_{\min}$ on all of $X^{[1],ss}_{\min}(\widehat{X}_{(n_1)})$. Moreover, by repeated application of Proposition \ref{prop properties first nonred blowup}, $\widehat{X}_{(n_1)}$ satisfies \hyperlink{QiS}{(QiS)} and \hyperlink{WUU}{(WUU)}.

We then consider $U^{[2]}$ and for $i=1,\dots n_2$, define $\widetilde{X}_{(n_1+i)}$ to be the blow-up of $\widetilde{X}_{(n_1 + i-1)}$ along the closure of $D^{[2]}(\widetilde{X}_{(n_1 + i-1)})$. By running the proof of Proposition \ref{prop properties first nonred blowup} with $U^{[1]}$ replaced by $U^{[2]}$, this terminates with a scheme $\widetilde{X}_{(n_1 +n_2)}$ satisfying \hyperlink{QiS}{(QiS)} and \hyperlink{WUU}{(WUU)}, and such that for $i\leq 2$ we have $\dim \stab_{U^{[i]}}$ is constant and equal to $d^{[i]}_{\min}$ on all of $X^{[i],ss}_{\min}(\widehat{X}_{(n_1+n_2)})$. Here, for $i=2$, this follows from the blow-up construction, and for $i=1$, this is inherited from $\widehat{X}_{(n_1)}$ by Lemma \ref{lem stab cant incr on blowup}.

After repeating this procedure for each $U^{[i]}$ we obtain a blow-up $\widetilde{X}_{(N)}$ with $N = n_1 + \cdots + n_{l-1}$ satisfying \hyperlink{QiS}{(QiS)} and \hyperlink{WUU}{(WUU)}, and such that for each $1\leq i\leq l-1$ we have \[\dim \stab_{U^{[i]}}(x)=d^{[i]}_{\min}\text{ for all }x\in X^{[i],ss}_{\min}(\widehat{X}_{(n_1+n_2)});\] that is the stabiliser assumptions concerning each $U^{[i]}$ in \ref{starUstagesibullet} are satisfied.

We then consider the groups $U^{[i-1,i]}$ and inductively perform blow-ups along the closures of $D^{[i-1,i]}(\widetilde{X}_{(N)})$ and so that the dimensions of the $U^{[i-1,i]}$-stabilisers are constant on $ X^{[i],ss}_{\min}$. This procedure determines a sequence of blow-ups
\[ \pi : \widetilde{X}  \ra \cdots \ra  \widetilde{X}_{(1)} \ra X\]
such that on $\widetilde{X}$ we have \hyperlink{QiS}{(QiS)} and \ref{starUstagesibullet}, which proves (1).

For (2), since  $X^{P-qs}$ is disjoint from the centre of each blow-up, we can identify it with its preimage under $\pi$. In fact, by the argument in the proof of Proposition \ref{prop properties first nonred blowup} (d) we have $\pi^{-1}(X^{P-qs}) \subset \widetilde{X}^{P-qs}$. Now we can apply Proposition \ref{prop qnt in stages no blow ups} to construct a quasi-projective geometric quotient of the open subset $\widetilde{X}^{P-qs}$ of $\widetilde{X}$, possibly after performing a further sequence of blow-ups $\widehat{X} \ra \widetilde{X}$ by considering reductive stabilisers as in \S\ref{sec upstairs red fails}. This restricts to the desired quasi-projective geometric quotient of $X^{P-qs}$. The natural projective completion is given by the projective $P$-quotient of $\widehat{X}$, which satisfies both \hyperlink{QiS}{(QiS)} and \hyperlink{U}{(U)}.
\end{proof}

We can now draw all of these results together into a proof of Theorem \ref{mainthm2}.

\begin{proof}[Proof of Theorem \ref{mainthm2}]
	If both \hyperlink{QiS}{(QiS)} and \hyperlink{U}{(U)} hold, this is Proposition \ref{prop no blow ups get proj qnt}. If only \hyperlink{QiS}{(QiS)} and \ref{starUstagesibullet} hold, this is Proposition \ref{prop qnt in stages no blow ups}. The final case, where only \hyperlink{QiS}{(QiS)} and \hyperlink{WUU}{(WUU)} hold, is given by Theorem \ref{prop blowup process if upstairs unip fails} above.
\end{proof}

\subsection{Removing the irreducibility assumption} \label{sec remove irred assumption} The irreducibility assumption for $X$ is mainly for bookkeeping purposes, and is not necessary for our proof of Theorem \ref{mainthm2} to work.

For reducible $X$, we first replace $X$ with the closure of $X_{\min}$, so that $X_{\min}$ is still open and dense\footnote{It may be in fact that the weight space we are interested in is not the minimal one, in which case we can replace $X$ with the closure of the corresponding attracting set. This may happen, for example, when the non-degeneracy assumption \ref{non-deg} of \hyperlink{QiS}{(QiS)} fails, or when \hyperlink{UU}{(UU)} does not hold, \emph{cf.} Remark \ref{rmk QiS stable set is open}.} The proof of Theorem \ref{mainthm nrGIT quotient} is local on affines, and so goes through verbatim, where we make sure to take the affines small enough that they are irreducible. The rest of the argument, in the case where no blow-ups are necessary, is the same.
 
The only other necessary modification to the proof is required considering stabiliser dimensions for the reductive and non-reductive blow-up processes, as there may be irreducible components of $X$ for which the generic stabiliser dimension (for a particular group) is not equal to the global minimum dimension. For this, as with the weight spaces, we can simply replace $X$ with the union of those irreducible components having a point with global minimum stabiliser dimension; once we have done this, the locus of points with minimum stabiliser will be open and dense in the resulting space. Alternatively, as discussed after Definition \ref{def centres of blowups}, one might prefer to think of the \lq minimum' stabiliser dimension as being in fact a vector, recording the minimum on each irreducible component.

\section{An overview of applications} \label{sec applications}

The assumptions of Theorem \ref{mainthm2} are suited to the construction of moduli spaces of objects of fixed Harder-Narasimhan type in an abelian category (see Remark \ref{rmk on qis assumptions}). In this section, we outline how to apply this theorem and discuss when this can be applied to sheaves of fixed HN type on a polarised projective scheme.

Given a moduli problem in an abelian category with a moduli space of semistable objects constructed as a reductive GIT quotient of a parameter scheme, to construct moduli spaces of objects of fixed Harder-Narasimhan type $\tau$ one carries out the following procedure:
\begin{enumerate}
\item Compare the HKKN stratification and HN stratification on the parameter scheme appearing in the reductive GIT construction. Ideally identify the HN stratum $S_\tau$ with (a closed subscheme of) a unstable HKKN stratum. This gives a scheme $X_\tau$ (the closure of the subscheme $Y_\tau^{ss} \subset S_\tau$ in the parameter space) with an action a parabolic subgroup $P_\tau$, which we want to quotient by.
\item Determine whether \hyperlink{QiS}{(QiS)} holds: it should be straight-forward to verify the first two conditions and the third, non-degeneracy, should admit a moduli-theoretic interpretation similarly to Definition \ref{def HN type non-degen}. 
\item  Interpret the reductive stabiliser groups moduli-theoretically. For coprime HN types (where semistability coincides with stability for all HN subquotients), one may expect Assumption \ref{starR0stagesibullet} to hold similarly to the situation for sheaves.
\item Interpret the unipotent stabiliser groups moduli theoretically and check whether \ref{starUstagesibullet}, or the weaker assumption \hyperlink{WUU}{(WUU)}, holds.
\item Give a moduli-theoretic interpretation of the quotienting-in-stages stable locus $X_\tau^{P_\tau-qs}$ similarly to $\tau$-stability given in Definition \ref{def tau stable}. 
\end{enumerate}
If all steps in this procedure can be achieved, one obtains a quasi-projective moduli space of certain so-called $\tau$-stable objects of HN type $\tau$ and, in good cases, this moduli space will even be projective. 

Examples of moduli problems for which one could attempt to carry out the above include moduli of sheaves, Higgs sheaves\footnote{One difference to the standard sheaves case is that in the Higgs case there are two different notions of HN type: one can use either the HN type of the Higgs pair, or the HN type of the underlying sheaf. Both will be considered in \cite{HHJ}.} and representations of quivers, which would require rephrasing as a projective GIT set-up in order to apply results of NRGIT. In all these examples, the first of the above steps have already been carried out in \cite{GSZ,HK,Hoskins_quivers,Hoskins,Eloise}, and the second, third and fifth are expected to be routine: the real crux is the fourth. 

For example, fix a non-degenerate HN type $\tau$ of arbitrary length for sheaves over a polarised projective base scheme $B$, and consider the problem of constructing a moduli space for sheaves of HN type $\tau$. We have seen above in $\S$\ref{sec QiS Ass for sheaves} and $\S$\ref{sec U for sheaves} that, while Assumption \hyperlink{QiS}{(QiS)} holds, the Assumption  \hyperlink{WUU}{(WUU)} does not hold in general, which limits what we can presently say. However, we can give a partial generalisation of the length two case in \cite{Josh_length2} to higher lengths. First we need some definitions.

\bdf (Refined HN types and refined non-degeneracy)
\bnu\item  A refined Harder-Narasimhan type $\widetilde{\tau} = (\tau,d)$ consists of a Harder-Narasimhan type $\tau = (\tau_1,\dots \tau_l)$ and vector $d = (d^{[1]},\dots d^{[l]},d^{[1,2]},\dots d^{[l-1,l]})$ of natural numbers. We say that a sheaf $\cF$ has refined HN type $\widetilde{\tau}$ if it has HN type $\tau$ and we have \[\dim \Hom(\cF/\cF^{(i)},\cF^{(i)}) = d^{[i]}\]
\[\dim \Hom(\cF/\cF^{(i)},\cF^{(i-1)}) = d^{[j,j+1]}\]  for all $i= 1 \dots l$ and $j =1\dots l-1$.
\item By analogy with Definition \ref{def tau stable}, we say that a sheaf $\cF$ of refined HN type $\widetilde{\tau}$ is $\tilde{\tau}$-stable if each inclusion $\cF^{(i-1)}\subset \cF^{(i)}$ in the HN filtration of $\cF$ is non-split and all subquotients $\cF^{(i)}/\cF^{(i-1)}$ in the HN filtration are stable.
 \item By analogy with Definition \ref{def HN type non-degen} we say that $\widetilde{\tau}$ is non-degenerate if there exists a $\widetilde{\tau}$-stable sheaf of refined HN type $\widetilde{\tau}$.  
\enu 
\edf 

\bthm\label{thm sheaves} 
Let $B$ be a polarised projective scheme and $\widetilde{\tau}=(\tau,d)$ be a non-degenerate refined HN type for sheaves on $B$. If there exists a sheaf $\cF$ of type $\widetilde{\tau}$ such that $\gr (\cF) \cong \cF$, then there exists a quasi-projective moduli space of $\tilde{\tau}$-stable sheaves.
\ethm 
\bpf
Consider the corresponding $P_\tau$-action on $X:=\overline{Y_\tau}$ as in $\S$\ref{sec sheaves nonred git setup}. Let $Z_{\min}^d \subset Z_{\min}$ be the locus consisting of points with stabiliser vector $d$, and let $X^d_{\min} = p^{-1}(Z_{\min}^d)\subset X$. Replacing $X$ with the closure of $X^d_{\min}$, we obtain a projective scheme $W$ with a $P_\tau$-action, for which we claim we can verify the conditions of Theorem \ref{mainthm2}. By non-degeneracy of $\widetilde{\tau}$ we obtain non-degeneracy for $W$, and the first two assumptions of \hyperlink{QiS}{(QiS)} follow directly from Proposition \ref{prop QiS 1 and 2 hold for sheaves}. Using Proposition \ref{prop first part of QiS stable for sheaves} and Lemma \ref{lem sheaves p preserves stabs} we deduce that $W^{P-qs}$ consists exactly of the $\widetilde{\tau}$-stable sheaves. Hence \hyperlink{WUU}{(WUU)} holds as $W^{P-qs} \neq \emptyset$ and by assumption there exists a sheaf $\cF$ of type $\widetilde{\tau}$ such that $\gr (\cF) \cong \cF$. Then the result follows by applying Theorem \ref{mainthm2}. \epf 

In particular, to obtain a moduli space for an refined HN type, we must assume the existence of a graded sheaf of that type. A natural question is for which extended HN types this assumption is satisfied, and in particular whether it is satisfied for the minimum stabiliser vector for sheaves of a given HN type. 

In the setting of Higgs sheaves of Higgs-HN type $\tau$ on a base scheme $S$, one can show that, provided there exists a sheaf on $S$ of HN type $\tau$, it is possible to construct a Higgs bundle $\cE$ on $S$, such that $\cE \cong \gr^{\text{Higgs}-HN} (\cE)$ and this has trivial unipotent stabiliser \cite{HHJ}. Therefore, via the spectral correspondence (for example, see \cite{Simpson}), one can, at least in principle, construct examples of base schemes $B$ from $S$ and HN types for sheaves from the Higgs-HN type $\tau$ such that this assumption is satisfied for the minimal stabiliser vector. The case of Higgs bundles will be dealt with in more detail in upcoming work \cite{HHJ}.

	\bibliographystyle{amsplain}
	\bibliography{references}

\providecommand{\bysame}{\leavevmode\hbox to3em{\hrulefill}\thinspace}
\providecommand{\MR}{\relax\ifhmode\unskip\space\fi MR }
\providecommand{\MRhref}[2]{%
  \href{http://www.ams.org/mathscinet-getitem?mr=#1}{#2}
}
\providecommand{\href}[2]{#2}
\begin{thebibliography}{10}

\bibitem{BDHK2}
G.~B\'erczi, B.~Doran, T.~Hawes, and F.~Kirwan, \emph{Projective completions of
  graded unipotent quotients}, arXiv preprint 1607.04181.

\bibitem{BDHK_handbook}
G.~B\'{e}rczi, B.~Doran, T.~Hawes, and F.~Kirwan, \emph{Constructing quotients
  of algebraic varieties by linear algebraic group actions}, Handbook of group
  actions. {V}ol. {IV}, Adv. Lect. Math. (ALM), vol.~41, Int. Press,
  Somerville, MA, 2018, pp.~341--446.

\bibitem{BDHK}
G.~B\'erczi, B.~Doran, T.~Hawes, and F.~Kirwan, \emph{Geometric invariant
  theory for graded unipotent groups and applications}, J. Topology \textbf{11}
  (2018), no.~3, 826--825.

\bibitem{BDK_grad_lin}
G.~B\'erczi, B.~Doran, and F.~Kirwan, \emph{Graded linearisations}, arXiv
  preprint 1703.05226.

\bibitem{BHJK}
G.~B\'erczi, V.~Hoskins, J.~Jackson, and F.~Kirwan, \emph{Moduli of unstable
  sheaves}, In preparation.

\bibitem{BHK}
G.~B\'{e}rczi, V.~Hoskins, and F.~Kirwan, \emph{Stratifying quotient stacks and
  moduli spaces}, Geometry of Moduli: Abel Symp. \textbf{14} (2018), 1--34.

\bibitem{Berczi2018a}
G.~B\'{e}rczi, J.~Jackson, and F.~Kirwan, \emph{Variation of non-reductive
  geometric invariant theory}, Surveys in differential geometry 2017.
  {C}elebrating the 50th anniversary of the {J}ournal of {D}ifferential
  {G}eometry, Surv. Differ. Geom., vol.~22, Int. Press, Somerville, MA, 2018,
  pp.~49--69.

\bibitem{BK_hyp}
G.~B\'erczi and F.~Kirwan, \emph{Non-reductive geometric invariant theory and
  hyperbolicity}, arXiv preprint 1607.04181.

\bibitem{BB}
A.~Bia{\l}ynicki-Birula, \emph{Some properties of the decompositions of
  algebraic varieties determined by actions of a torus}, Bull. Acad. Polon.
  Sci. S\'er. Sci. Math. Astronom. Phys. \textbf{24} (1976), no.~9, 667--674.
  \MR{0453766}

\bibitem{BrambilaPaz2013}
Mata-Gutiérrez~O. Brambila-Paz~L., \emph{Un teorema de torelli para haces
  vectoriales inestables}, Aportaciones Matemáticas - Memorias de la Sociedad
  Matemática Mexicana, 46, pp. 03-16 (2013).

\bibitem{dom}
D.~Bunnett, \emph{On the moduli of hypersurfaces in toric orbifolds}, arXiv
  preprint 1906.00272.

\bibitem{BK_cohomolgy}
Gergely Bérczi and Frances Kirwan, \emph{Moment maps and cohomology of
  non-reductive quotients}, 2021.

\bibitem{Dolgachev}
I.~Dolgachev, \emph{Lectures on invariant theory}, Cambridge University Press,
  2003.

\bibitem{DH}
I.~Dolgachev and Y.~Hu, \emph{Variation of geometric invariant theory
  quotients}, Inst. Hautes {\'E}tudes Sci. Publ. Math. (1998), no.~87, 5--56,
  With an appendix by Nicolas Ressayre.

\bibitem{GSZ}
T.~Gómez, I.~Sols, and A.~Zamora, \emph{A {GIT} interpretation of the
  {H}arder--{N}arasimhan filtration}, Rev. Mat. Complut. \textbf{28} (2015),
  169--190.

\bibitem{HL_Theta}
D.~Halpern-Leistner, \emph{{$\Theta$}-stratifications, {$\Theta$}-reductive
  stacks, and applications}, Algebraic geometry: {S}alt {L}ake {C}ity 2015,
  Proc. Sympos. Pure Math., vol.~97, Amer. Math. Soc., Providence, RI, 2018,
  pp.~349--379.

\bibitem{Eloise}
E.~Hamilton, \emph{Stratifications and quasi-projective coarse moduli spaces
  for the stack of higgs bundles}, arXiv:1911.13194.

\bibitem{HHJ}
E.~Hamilton, V.~Hoskins, and J.~Jackson, \emph{Instability stratifications and
  moduli spaces for {H}iggs sheaves}, In preparation.

\bibitem{HN}
G.~Harder and M.~S. Narasimhan, \emph{On the cohomology groups of moduli spaces
  of vector bundles on curves}, Math. Ann. \textbf{212} (1975), 215--248.

\bibitem{Hesselink}
W.~H. Hesselink, \emph{Uniform instability in reductive groups}, J. Reine
  Angew. Math. \textbf{304} (1978), 74--96.

\bibitem{Hoskins_quivers}
V.~Hoskins, \emph{Stratifications associated to reductive group actions on
  affine spaces}, Quart. J. Math. Oxford \textbf{65} (2014), no.~3, 1011--1047.

\bibitem{Hoskins}
V.~Hoskins, \emph{Stratifications for moduli of sheaves and moduli of quiver
  representations}, Algebr. Geom. \textbf{5} (2018), no.~6, 650--685.

\bibitem{HK}
V.~Hoskins and F.~Kirwan, \emph{Quotients of unstable subvarieties and moduli
  spaces of sheaves of fixed {H}arder--{N}arasimhan type}, Proc. London Math.
  Soc \textbf{105} (2012), no.~4, 852--890.

\bibitem{HL}
D.~Huybrechts and M.~Lehn, \emph{The geometry of moduli spaces of sheaves},
  Aspects of Mathematics, Vieweg, 1997.

\bibitem{Josh_length2}
J.~Jackson, \emph{Moduli spaces of unstable objects: Sheaves of
  {H}arder-{N}arasimhan length 2}, arXiv preprint.

\bibitem{Kempf}
G.~R. Kempf, \emph{Instability in invariant theory}, Ann. of Math. \textbf{108}
  (1978), no.~2, 299--316.

\bibitem{Kirwan_thesis}
F.~Kirwan, \emph{Cohomology of quotients in symplectic and algebraic geometry},
  Mathematical Notes, no.~31, Princeton University Press, 1984.

\bibitem{K2}
F.~Kirwan, \emph{Partial desingularisations of quotients of nonsingular
  varieties and their betti numbers}, Ann. of Math. \textbf{122} (1985), no.~2,
  41--85.

\bibitem{Mumford}
D.~Mumford, J.~Fogarty, and F.~Kirwan, \emph{Geometric invariant theory}, third
  ed., Springer, 1993.

\bibitem{Nagata}
M.~Nagata, \emph{On the fourteenth problem of {H}ilbert}, Proc. {I}nternat.
  {C}ongress {M}ath. 1958, Cambridge Univ. Press, New York, 1960, pp.~459--462.
  \MR{0116056}

\bibitem{Ness}
L.~Ness, \emph{A stratification of the null cone via the moment map (with an
  appendix by {D}. {M}umford)}, Amer. J. Math. \textbf{106} (1984), no.~6,
  1281--1329.

\bibitem{Newstead}
P.~E. Newstead, \emph{Introduction to moduli problems and orbit spaces},
  T.I.F.R. Lecture Notes, Springer-Verlag, 1978.

\bibitem{Nitsure}
N.~Nitsure, \emph{Schematic {H}arder--{N}arasimhan filtration}, Internat. J.
  Math \textbf{22} (2011), no.~10, 1365--1373.

\bibitem{Shatz}
S.~Shatz, \emph{The decomposition and specialization of algebraic families of
  vector bundles}, Compositio Math. \textbf{35} (1977), 163--187.

\bibitem{Simpson}
C.~T. Simpson, \emph{Moduli of representations of the fundamental group of a
  smooth projective variety}, Inst. Hautes Etudes Sci. Publ. Math. \textbf{79}
  (1994), 47--129.

\bibitem{Thaddeus}
M.~Thaddeus, \emph{Geometric invariant theory and flips}, J. Amer. Math. Soc.
  \textbf{9} (1996), no.~3, 691--723.

\end{thebibliography}

\end{document}